\documentclass[draftclsnofoot,onecolumn,journal,12pt]{IEEEtran}

\usepackage{amssymb,graphicx}
\usepackage{latexsym}
\usepackage{mathtools}
\usepackage{cite}
\usepackage{color}
\usepackage{subfigure}
\usepackage{bm,bbm}
\usepackage{tabularx,booktabs}
\usepackage{ltablex}
\usepackage{mathrsfs}
\usepackage{enumitem}
\usepackage{cases}
\usepackage{dsfont}

\newcolumntype{M}[1]{>{\centering\arraybackslash}m{#1}}

\usepackage[linesnumbered,lined,boxed]{algorithm2e}
\usepackage{setspace}
\SetAlCapSkip{1em}
\SetKwComment{Comment}{$\triangleright$\ }{}
\RequirePackage{amsthm,amsmath}
\RequirePackage[colorlinks,citecolor=blue,urlcolor=blue]{hyperref}
\newcommand*{\fullref}[2]{\hyperref[{#1}]{\ref*{#1} {\it(#2)}}}
\newcommand*{\partialref}[2]{\hyperref[{#1}]{#2}}
\usepackage[deletedmarkup=xout]{changes}
\graphicspath{{./figs/}}

\renewcommand{\proof}[1]{\emph{#1}}
\renewcommand{\endproof}{{\hfill $\square$}}

\newtheorem{definition}{\bf Definition}
\newtheorem{theorem}{\bf Theorem}
\newtheorem{lemma}{\bf Lemma}
\newtheorem{proposition}{\bf Proposition}
\newtheorem{corollary}{\bf Corollary}

\newcounter{condition}
\newenvironment{condition}[1]{\refstepcounter{condition}
\par\medskip
\noindent \textbf{#1:~} \it}{\medskip \newline}

\newcommand{\norm}[1]{\left\lVert #1 \right\rVert}

\begin{document}

\title{A Restless Bandit Model for Resource Allocation, Competition and Reservation}
\author{Jing~Fu,
Bill~Moran~and~Peter~G.~Taylor
\thanks{
Jing Fu was with School of Mathematics and Statistics, The University of Melbourne, VIC 3010, Australia, and she is now with School of Engineering, RMIT University, VIC3001, Australia (e-mail: jing.fu@rmit.edu.au).
}
\thanks{Peter G. Taylor is with School of Mathematics and Statistics, The University of Melbourne, VIC 3010, Australia (e-mail: taylorpg@unimelb.edu.au).}
\thanks{Bill Moran is with Department of Electrical and Electronic Engineering, The University of Melbourne, VIC 3010, Australia (e-mail:wmoran@unimelb.edu.au).}
}
\maketitle

\begin{abstract}
We study a resource allocation problem with varying requests, and with resources of limited capacity shared by multiple requests.
It is modeled as a set of heterogeneous Restless Multi-Armed Bandit Problems (RMABPs) connected by constraints imposed by resource capacity. 
Following Whittle's relaxation idea and Weber and Weiss' asymptotic optimality proof, we propose a simple policy and prove it to be asymptotically optimal in a regime where both arrival rates and capacities increase.  
We provide a simple sufficient condition for asymptotic optimality of the policy, and in complete generality propose a method that generates a set of candidate policies for which asymptotic optimality can be checked. 
The effectiveness of these results is demonstrated by numerical experiments.
To the best of our knowledge, this is the first work providing asymptotic optimality results for such a resource allocation problem and such a combination of multiple RMABPs.
\end{abstract}



\section{Introduction}\label{sec:introduction}

\subsection{Overview and Motivation}

Modern technologies enable Internet resources such as routers, computing servers and cables to be abstracted from the physical layer to a \emph{virtual} layer, facilitating a quick response to demands for setting up communication networks or processing computing jobs. Virtual servers comprising different sets of physical resources are assigned to arriving customers who use these resources for a period of time and then return them to a pool when they depart.

Such networks are just particular examples of more general systems where users of different types arrive with a desire to be allocated resources of various kinds, to use these resources and then return them. Users are often indifferent to the precise set of resources that they are allocated, they just require allocation of some resources that will enable them to accomplish the task at hand. In such circumstances a network manager has the task of deciding whether an arriving customer should be admitted into the system and, if so, which set of resources should be assigned to satisfying their requirements.

In this paper we describe and analyze a very general model for such systems. Specifically, we study a system in which $J$ \emph{resource pools}, each made up of finite numbers of \emph{resource units} (RUs), await allocation to incoming requests of $L$ different types. We refer to the number of RUs in a resource pool as its \emph{capacity}. Each resource
pool is potentially shared and \emph{competed} for by many requests, but \emph{reservation} of RUs for still-to-arrive requests is also allowed. When a request has been accommodated by a resource pool, an appropriate number of RUs of this type are occupied by the request until it leaves the system. The released RUs can be reused by other requests. A request is permitted to occupy RUs from more than one resource pool simultaneously.  In this context, the number of requests of the same type that are accommodated by a group of resource pools varies according to a stochastic process, where the transition rates are affected by the resource allocation policy employed. Several such processes associated with the same resource pool are coupled by its capacity limitations.

By strategically assigning requests to appropriate combinations of RUs, we aim to maximize the long-run average revenue, defined as the difference between the long-run average reward earned by serving the requests and the long-run average cost incurred by using the resource pools. 
Such a resource allocation problem can be easily applied to a rich collection of classical models, such as loss networks in telecommunications, resource allocation for logistic systems, and job assignment in parallel computing. 

\cite{kelly1991loss} published a comprehensive analysis of \emph{loss network} models with and without \emph{alternative routing}. In the latter case, network traffic can be re-routed onto alternative paths when the original path fails or is full. 
In \cite{kelly1991loss}, a list of alternative paths as choices of resource pools is given for each call/request. The alternative paths are selected in turn after if preceding offered paths are unavailable. In contrast, the manager of a typical resource allocation problem described above is potentially able to change the priorities of paths dynamically. How this should be done is a key focus of this paper.

\begin{figure}[t]
\centering
\begin{minipage}{.45\textwidth}
\centering
\includegraphics[width=0.4\linewidth]{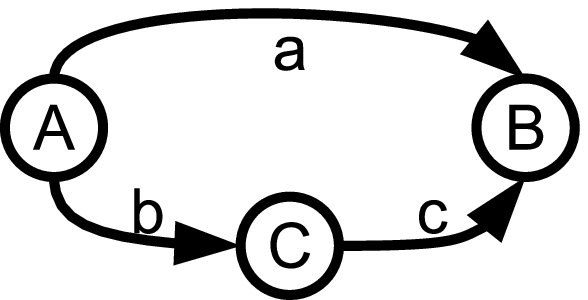}
\caption{A simple loss network.}\label{fig:loss_network}
\end{minipage}
\begin{minipage}{.45\textwidth}
\centering
\includegraphics[width=0.4\linewidth]{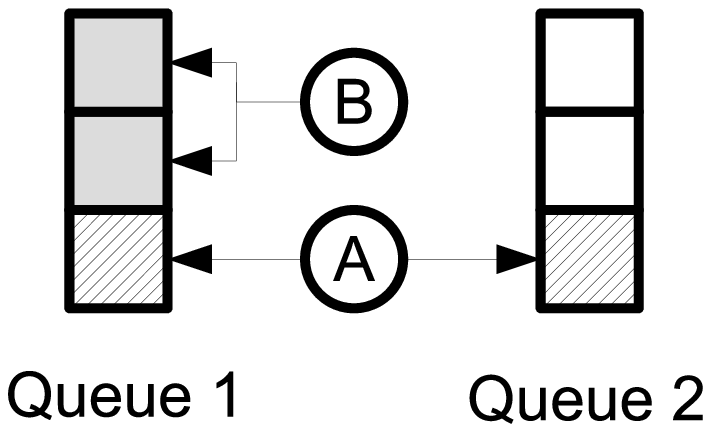}
\caption{A simple parallel queueing model.}\label{fig:queue_model}
\end{minipage}
\vspace{-0.5cm}
\end{figure}
To illustrate the kind of problem of interest here, consider the simple loss network model shown in Figure~\ref{fig:loss_network}. Links $a$, $b$ and $c$ are abstracted as resource pools with capacities equal to 1, 3 and 3, respectively: link $a$ consists of one channel as an RU, and links $b$ and $c$ each have 3 channels. Requests asking for a connection from $A$ to $B$ occupying one channel can be served by either path $\{a\}$ or $\{b,c\}$, but requests requiring two channels for each connection from $A$ to $B$ are able to be accommodated only by path $\{b,c\}$. We refer to the former and the latter as type-I and type-II requests, respectively.  An arrival of a type-I request results in one of the paths $\{a\}$ and $\{b,c\}$ being chosen by the optimizer depending on current traffic loads on the three links, where links $b$ and $c$ might be shared with existing type-II requests. Occupied channels or RUs are released immediately and simultaneously when relevant requests are completed.  

Resource allocation problems with small values of $L$ and $J$, such as the example above, can be modeled by a Markov Decision Processes (MDP), and solved through dynamic programming. However, in real-world applications, where $L$ and $J$ are large, resulting in high dimensionality of the state and action spaces, such an approach is often intractable.

In this paper we use an analysis inspired by techniques applied to Restless Multi-Armed Bandit Problems (RMABPs).
The standard RMABP consists of parallel MDPs with binary actions (they can either be ``pulled', that is activated, or not), which are competing for a limited possibility of being selected at each decision epoch. Each of the MDPs, referred to as a \emph{bandit process}, has its own individual state-dependent reward rates and transition probabilities when it is activated and when it is not. 

Attempts to solve the problem are faced with exponential growth in the size of the state space as the number of parallel bandit processes increases.  
This class of problems was described  by \cite{whittle1988restless}, who proposed a heuristic management policy that was shown to be  asymptotically optimal under non-trivial extra conditions  by~\cite{weber1990index};
this policy approaches optimality as the number of bandit processes tends to infinity.
The policy, subsequently referred to as the \emph{Whittle index policy}, always prioritizes bandit processes with higher state-dependent \emph{indices} that intuitively represent marginal rewards earned by processes if they are selected.
The Whittle indices can be computed independently for each bandit process - a process that imposes significantly reduced computational complexity. The Whittle index policy is scalable to a RMABP with a large number of bandit processes.
Also, the asymptotic optimality property, if it is satisfied, guarantees a bounded performance degradation in a large-scale system and is appropriate for large problems where optimal solutions are intractable.
The non-trivial extra conditions required by the asymptotic optimality proof in \cite{weber1990index} are related to proving the existence of a global attractor of a stochastic process.

RMABPs have been widely used in scheduling problems, such as channel detecting (see \cite{liu2012learning,wang2019whittle} ), job assignments in data centers (see \cite{fu2016asymptotic}), web crawling (see \cite{avrachenkov2016whittle}), target tracking (see \cite{krishnamurthy2007structured,le2010scheduling}) and job admission control (see \cite{nino2012admission,nino2019resource}).
Here we treat the resource allocation problem described above as a set of RMABPs coupled by linear inequalities involving random state and action variables.

\subsection{Main Contributions}

We propose a modified \emph{index policy} that takes into account the capacity constraints of the problem. The index policy prioritizes combinations of RUs with the highest indices, each of which is a real number representing the marginal revenue of using its associated RUs. The policy is simple, scalable and appropriate for a large scale resource allocation problem.

Our analysis of asymptotic optimality of the index policy proceeds through a relaxed version of the problem and study of a global attractor of a stochastic process defined in \eqref{eqn:z_process} below.
We prove that the process \eqref{eqn:z_process}  will almost surely converge to a global attractor in the asymptotic regime regardless of its initial point, and hence the index policy is asymptotically optimal if and only if this global attractor coincides with an optimal solution of the resource allocation problem. 
Following ideas similar to those of \cite{weber1990index}, optimality of the global attractor for the resource allocation problem can be deduced from its optimality for the relaxed problem, which can be analyzed with  remarkably reduced computational complexity. 

A sufficient condition for the global attractor and optimal solution to coincide is that the offered traffic for the entire system is \emph{heavy} and the resource pools in our system are \emph{weakly coupled}.
We rigorously define these concepts in Section~\ref{subsec:sufficient_condition}.
These results are enunciated in Theorems~\ref{theorem:main_second}, \ref{theorem:main} and Corollary~\ref{coro:main} in Section~\ref{sec:asym_opt}.

When the above-mentioned sufficient conditions are not satisfied, an asymptotically optimal index policy can still exist. In this case, we propose a method that can derive the parameters required by the asymptotically optimal policy.
Although asymptotic optimality is not guaranteed, 
Theorem~\ref{theorem:main} provides a verifiable sufficient condition, less stringent than the one mentioned above, to check asymptotic optimality of the index policy with adapted parameters. 
We numerically demonstrate the effectiveness of this method in Section~\ref{sec:example}.

The index policy exhibits remarkably reduced computational complexity, compared to conventional optimizers, and its potential asymptotic optimality is appropriate for large-scale systems where computational power is a scarce commodity. Furthermore, simulation studies indicate that an index policy can still be good in the pre-limit regime.
As mentioned earlier, our problem can be seen as a set of RMABPs coupled by the capacity constraints. When the capacities of all resource pools tend to infinity,  the index policy reduces to the Whittle index policy because the links between RMABPs no longer exist.

To the best of our knowledge, no existing work has proved asymptotic optimality in resource allocation problems, where resource competition and reservation are potentially permitted, nor has there been a previous analysis of such a combination of multiple, different RMABPs, resulting in a much higher dimensionality of the state space.

The remainder of the paper is organized as follows.
In Section~\ref{sec:model}, we describe the resource allocation problem.
In Section~\ref{sec:relaxation}, we apply the Whittle relaxation technique. 
In Section~\ref{sec:index_policy}, we propose an algorithm to implement an index policy.
In Section~\ref{sec:asym}, we define the asymptotic regime and we prove the asymptotic optimality of the index policy under some conditions.
To demonstrate the effectiveness of the proposed policies, numerical results are provided in Section~\ref{sec:example}.
In Section~\ref{sec:conclusions}, we present conclusions.

\subsection{Relation to the Literature}
The classical Multi-Armed Bandit Problem (MABP) is a optimization problem in which only one bandit process (BP) among $K$ BPs can be activated at any one time, while all the other $K-1$ BPs are \emph{frozen}: an active BP randomly changes its state, while state transitions will not happen to the frozen BPs.
In 1974, Gittins and Jones published the well-known \emph{index theorem} for the MABP \cite{gittins1974dynamic}, and in 1979, \cite{gittins1979bandit} proved the optimality of a simple \emph{index policy}, subsequently referred to as the \emph{Gittins index policy}.
Under the Gittins index policy, an index value, referred to as the \emph{Gittins index}, is associated with each state of each BP, and the BP with the largest index value is activated, while all the other BPs are frozen.
More details about Gittins indices can be found in \cite[Chapter 2.12]{gittins2011multiarmed} (and the references therein).

The optimality of the Gittins index policy for the conventional MABP fails for the general case where the $K-1$ BPs that are not selected can also change their states randomly; such a process is known as a Restless Multi-Armed Bandit Process (RMABP).
The RMABP was proposed by \cite{whittle1988restless}.
The RMABP allows $M=1,2,\ldots,K$ BPs to be active simultaneously.
In a similar vein to the Gittins index policy, Whittle assigned a state-dependent index value, referred to as the \emph{Whittle index}, to each BP and always activated the $M$ BPs with the highest indices. 
The Whittle indices are calculated from a \emph{relaxed} version of the original RMABP obtained by randomizing the action variables. 
\cite{whittle1988restless} defined a property of a RMABP, referred to as \emph{indexability}, under which the \emph{Whittle index policy} exists.
Whittle conjectured in \cite{whittle1988restless} that the Whittle index policy, if it exists, is \emph{asymptotically optimal}.
\cite{papadimitriou1999complexity} proved that the optimization of RMABPs is PSPACE-hard in general;
nonetheless, \cite{weber1990index} were able to establish asymptotic optimality of Whittle index policy under mild conditions.

\cite{nino2001restless} proposed a Partial Conservation Law
(PCL) for the optimality of RMABP;  this is an
extension of  the General Conservation Law (GCL) published in
\cite{bertsimas1996conservation}.  
Later, \cite{nino2002dynamic} defined a group of problems that satisfies
PCL-indexibility and proposed a new index policy that improved the
Whittle index. 
The new index policy was proved to be optimal for problems with PCL-indexibility.  
PCL-indexibility implies (and is stronger than) Whittle indexibility.
A detailed survey about the optimality of bandit problems can be found in \cite{nino2007dynamic}.

\cite{verloop2016asymptotically} proved the asymptotic optimality of the Whittle index policy in an extended version of an RMABP, where BPs randomly arrive and depart the system. She proposed an index policy that was not restricted to Whittle indexable models and numerically demonstrated its near-optimality. 
\cite{larranaga2015asymptotically} applied this extended RMABP to a queueing problem assuming convex, non-decreasing functions for both holding costs and measured values of people's impatience.
More results on  asymptotic optimality of index-like polices can be
found in \cite[Chapter IV]{fu2016thesis}. 

Asymptotically optimal policies for cost-minimization problems in network systems using a fluid approximation technique have been considered in \cite{bauerle2000asymptotic,bauerle2002optimal,stolyar2004maxweight,
nazarathy2009near} and \cite{bertsimas2015robust}.
The fluid approximation to the stochastic optimization problem can be much simpler than the original.
A key problem here is to establish an appropriate fluid problem and translate its optimal solution to a policy amenable to the stochastic problem.
Asymptotic optimality of the translated stochastic policy can be established if the fluid solution provides an upper/lower bound of the stochastic problem and the policy coincides with this bound asymptotically.
The reader is referred to \cite{meyn2008control} for a detailed description of fluid approximation across various models. 

Although the fluid approximation technique helps with asymptotic analysis in a wide range of (cost-minimization) network problems, existing results cannot be directly applied to our problem, where the arrival and departure rates of request queues are state-dependent and capacity violation over resource pools is strictly forbidden.
Our system is always stable for any offered traffic because of the strict capacity constraints.
In our case, the form of the corresponding fluid model remains unclear for generic policies. 
Even given the optimal solution of a well-established fluid model, the synthesis of an explicit policy in the stochastic model remains a challenge.

We adopt another approach, following   the ideas of \cite{whittle1988restless} and \cite{weber1990index}.
Our asymptotic optimality is derived from an optimal solution of a relaxed version of the stochastic optimization problem. The relaxed problem is still a stochastic optimization problem with a discrete state space. We propose a policy based on intuition captured by the relaxed problem, of which the optimal solution provides a performance upper bound of the original problem.
Then, we prove, under certain conditions, that this policy coincides with the upper bound asymptotically. 
The detailed analysis comprises the main content of the paper.

\section{A Resource Allocation Problem}\label{sec:model}

We use $\mathbb{N}_{+}$ and $\mathbb{N}_{0}$ to denote the sets of positive and non-negative integers, respectively, and for any $N\in\mathbb{N}_{+}$, let $[N]$ represent the set $\{1,2,\ldots,N\}$ with $[0]=\emptyset$.
Let $\mathbb{R}$, $\mathbb{R}_{+}$ and $\mathbb{R}_{0}$ be the set of all, positive and non-negative reals, respectively.

\subsection{System Model}\label{subsec:model}

Recall that there are $L$ types of requests and $J$ pools of RUs, all
potentially different, with  resource pool $j\in[J]$ having capacity $C_j$
RUs that can be dynamically allocated to and released by the $L$ types
of requests.

Each request comes with an associated list of candidate resource combinations. 
Specifically, requests from \emph{request type} $\ell \in [L]$ can be accommodated by one of a set $\mathscr{P}_{\ell}$ of candidate \emph{patterns}.
One of these candidate patterns will be selected by a policy.
Patterns are  indexed by $i\in\mathbb{N}_+$. 
If a request is accommodated by pattern $i$, $w_{j,i}$ 
RUs of pool $j\in [J]$ are occupied until the request is completed and
departs. 
We can thus identify pattern $i$ with the \emph{weight vector} $\bm{w}_i = (w_{j,i})$ that defines its requirement.
Preemption or re-allocation of requests
are not allowed.
A request is blocked if there is not enough capacity on any
of its corresponding patterns. We might also want to block a request in
other circumstances, if accepting it would be detrimental to future
performance.  In either case, we model the situation where a request
is blocked by assigning it to the dummy pattern $d(\ell)$ with the weight vector set to $\bm{0}$.



It is possible for different RTs to be satisfied by the same pattern (this occurs, in particular with the dummy pattern). In such cases, we consider there to be multiple copies of each pattern, one for each RT that it can satisfy. 
This enables us to consider the sets $\mathscr{P}_\ell$ 
to be mutually
exclusive; that is,
$\mathscr{P}_{\ell_1}\cap \mathscr{P}_{\ell_2} = \emptyset$ for any
$\ell_1\neq \ell_2$.  Given
$|\mathscr{P}_{\ell}|$ patterns for each RT $\ell$, we have in total
$I = \sum_{\ell\in[L]}|\mathscr{P}_{\ell}|$ patterns 
associated with weight vectors $\bm{w}_i\in\mathbb{N}_0^{J}$,
$i\in[I]$.
For any pattern $i$, let $\ell(i)$ be the unique RT that is satisfied by that pattern.

Let $\mathcal{W}$ be a $J\times I$ matrix with entries $w_{j,i}$.
We assume that there is  no row and exactly $L$ columns in $\mathcal{W}$ with all zero entries.
Each of these zero columns corresponds to one of the dummy patterns $d(\ell)$ where requests of type $\ell\in[L]$ are blocked.

Requests of RT $\ell$ arrive at the system sequentially, following a Poisson process, with rates $ \lambda_{\ell}$ 
and the occupation times of the requests accommodated by pattern $i\in\mathscr{P}_{\ell}$ are exponentially distributed with parameter $\mu_{i}$. 
Although there might be situations when it is reasonable to assume that the occupation time depends only on the request type $\ell$, there might also be cases where the lifetime of a request depends on the resources accommodating it, which is why we allow the occupation time distribution to depend on $i$.
The RUs used to accommodate a request are  occupied and released at the same time.
Neither the request nor the system knows the lifespan of a request until it is accomplished and departs the system.

Since there are similarities between our problem and a parallel
queueing model, we present a second example to clarify the similarities and differences. 
Consider two
resource pools corresponding to two queues as illustrated in
Figure~\ref{fig:queue_model}, where both capacities are set to three; that is, $J=2$ and $C_1=C_2=3$.  There are two types
of requests: if a type-one request is accommodated in the system, it
will simultaneously occupy one RU of both pools;
and a type-two request can be accommodated by two RUs of either
pool. In other words, $L=2$, $\mathscr{P}_{1}=\{1,2\}$,
$\mathscr{P}_2=\{3,4,5\}$, patterns $2$ and $5$ are dummy patterns
with $\bm{w}_2=\bm{w}_5=\bm{0}$, $\bm{w}_1=(1,1)$, $\bm{w}_3=(2,0)$,
$\bm{w}_4=(0,2)$ and $I=5$.

In this case, the number of occupied RUs in both resource pools may decrease or increase
by one simultaneously, or by two exclusively for  an arrival or
departure event.
The transition rates are affected by the system
controller: if the capacity constraints are not violated, there are
two choices, resource pool one or two, for accommodating a type-two request. The task of a system manager is to find a policy for deciding which of these choices to take in order to maximize some long-term objective.
Each choice will result in a parallel queueing model with dependencies
between the sizes of queues,  between the policy employed and
queue transition rates.   As mentioned in Section~\ref{sec:introduction},
conventional optimization methods cannot be  applied directly when $L$ and $J$ are large.


\subsection{A Stochastic Optimization Problem}\label{subsec:RMABPs:general_case}

We focus here an explanation of
the stochastic mechanism of the resource allocation problem.

An \emph{instantiation} is 
generated in the memory of the system controller
when a request of RT $\ell\in[L]$ is accommodated by a pattern $i\in\mathscr{P}_{\ell}$. 
Once the request departs the system, the associated instantiation will
be removed from the controller's memory.  As requests are accommodated
and completed, the number of instantiations associated with each
pattern forms a birth-and-death process,  indicating the
number of requests being served by this pattern.  As mentioned
in the second example, the birth-and-death processes for all patterns
$i\in[I]$ are coupled by capacity constraints and affected by control
decisions.

Let $N_i(t)$, $t\geq 0$,  represent the number of instantiations
for pattern $i$ at time $t$.  The
process $N_i(t)$ has state space $\mathscr{N}_i $ that  is a discrete,
finite set of possible values. The finiteness of $\mathscr{N}_i$ derives from the finite
capacities $C_j$.  
If  $N_i (t)$ is known for all $i \in [I]$, 
the number of occupied RUs in pool $j \in [J]$ at time $t$ is given by
$S_j(t)=\sum_{i \in [I]} w_{j,i}N_i(t)$, which must be less than $C_j$.
The  vector $\bm{N}(t)=(N_i(t):\ i\in[I])$ is  the state variable of the entire system taking values in
$\mathscr{N}\coloneqq\prod_{i\in[I]}\mathscr{N}_i$,
where $\prod$ represents Cartesian product.
Since the state variables are further subject to capacity constraints to be discussed in Section~\ref{subsubsec:capacity_constraints}, $\mathscr{N}$ is larger than necessary. With slightly abused notation, we still refer to $\mathscr{N}$ as the state space of the system. 

\subsubsection{Action Constraints}\label{subsubsec:action_constraints}
We associate an action variable $a_i(\bm{n})\in\{0,1\}$ with
process $i\in[I]$ when the system is
in state $\bm{n}\in\mathscr{N}$, and $\bm{a}(\bm{n})=(a_i(\bm{n}):\ i\in[I])$.  The action variable $a_i(\bm{n})$ tells us what to do with a potential new request of type $\ell(i)$. If $a_i(\bm{n}) =1$, then such a pattern will be allocated to pattern $i$.
The \emph{action constraint},
\begin{equation}\label{eqn:constraint:action}
\sum\limits_{i\in\mathscr{P}_{\ell}} a_i(\bm{n}) = 1,~
\forall \ell\in[L],~
\forall \bm{n}\in\mathscr{N},
\end{equation}
ensures that exactly one pattern (which may be the dummy pattern $d(\ell)$) is selected for each RT $\ell$ and current state $\bm{n}$.

At any time $t$, we say that the arrival process for pattern $i$ is \emph{active} or \emph{passive} according to whether $a_i(\bm{N}(t))$ is $1$ or $0$ respectively.
The birth rate of process $i\in\mathscr{P}_{\ell}$,
$\ell\in[L]$, is $\lambda_{\ell}$ if $a_i(\bm{N}(t))=1$;
and zero otherwise. The death rate of process $i$ is
$\mu_{i} N_{i}(t)$.
The time proportion that $a_{d(\ell)}(\bm{N}(t))=1$ is the \emph{blocking probability} for requests of type $\ell$.

\subsubsection{Capacity Constraints}\label{subsubsec:capacity_constraints}

To ensure feasibility of an allocation of a request of type $\ell(i)$ to pattern $i$ when the state is $\bm{n}$, we need 
\begin{equation}\label{eqn:constraint:resourcesi}
\mathcal{W}\left(\bm{n}+ \bm{e}_i\right) \leq \bm{C},~
\end{equation}
where $\bm{e}_i$ is a vector with a one in the $i$th position and zeros everywhere else and $\bm{C}\in\mathbb{N}_{+}^{J}$ is a vector with entries $C_j$. In view of the action constraint (\ref{eqn:constraint:action}), a neat way to collect together the constraints (\ref{eqn:constraint:resourcesi}) for all $i \in \mathscr{P}_{\ell}$ is to write them in the form
\begin{equation}\label{eqn:resources1}
\mathcal{W}\left(\bm{n}+\mathcal{E}_{\ell}\bm{a}(\bm{n})\right) \leq \bm{C},~\forall \bm{n}\in\mathscr{N},
\end{equation}
where $\mathcal{E}_{\ell}$ is a diagonal matrix of size $I$ with entries $e_{\ell,i,i}=1$ if $i\in\mathscr{P}_{\ell}$ and $e_{\ell,i,i}=0$ if $i\in[I]\backslash \mathscr{P}_{\ell}$.

For two different request types $\ell_1$ and $\ell_2$, a constraint of the form 
\begin{equation}\label{eqn:resources12}
\mathcal{W}\left(\bm{n}+\mathcal{E}_{\ell_1}\bm{a}(\bm{n})+\mathcal{E}_{\ell_2}\bm{a}(\bm{n})\right) \leq \bm{C},~\forall \bm{n}\in\mathscr{N},
\end{equation} 
captures the idea that the action vector $\bm{a}(\bm{n})$ must be such that the allocation decisions for $\ell_1$ and  
$\ell_2$ ensure enough capacity to implement both of them when both requests arrive simultaneously while the state is
$\bm{n}$. Another way to think about this is that, if a request of type $\ell_1$ is allocated to a non-dummy pattern $i_1$ when the state is $\bm{n}$, the decision for a request of type $\ell_2$ when the state is $\bm{n}$ must satisfy constraint \eqref{eqn:resources1} when the state is $\bm{n}+\bm{e}_{i_1}$. In particular, if there is not enough capacity to accommodate a request of type $\ell_2$ when the state is $\bm{n}+\bm{e}_{i_1}$, then a request of type $\ell_2$ must be allocated to the dummy pattern $d(\ell_2)$, when the state is $\bm{n}$.
This can be viewed as giving priority to \emph{reserving} resources for a type $\ell_1$ request over a type $\ell_2$ request
when the state is $\bm{n}$. 
As we shall see below, the decision to do this will be made in order to optimize a long-term reward function.

Observing that $\sum_{\ell \in [L]} \mathcal{E}_{\ell} = I$, we see that the constraint 
\begin{equation}\label{eqn:constraint:resources}
\mathcal{W}\left(\bm{n}+ \bm{a}(\bm{n})\right) = \mathcal{W}\Bigl(\bm{n}+ \bigl(\sum_{\ell \in [L]} \mathcal{E}_{\ell}\bigr) \bm{a}(\bm{n})\Bigr) \leq \bm{C},~
\forall \bm{n}\in\mathscr{N},
\end{equation}
can be thought of as an extended version of \eqref{eqn:resources12}.
In \eqref{eqn:constraint:resources}, requests of all types are taken into account when the state is $\bm{n}$  and allocation decisions for some types are made in order to reserve resources for other types that turn out to be more profitable in the long run. In particular, resources are reserved for those request types $\ell$ which are allocated to non-dummy patterns $i$ at the expense of those types that are allocated to less profitable patterns or the corresponding dummy patterns.
In this paper, all the results presented are based on capacity constraint~\eqref{eqn:constraint:resources}.

From \eqref{eqn:constraint:resources}, there is an upper bound, $\min_{j\in[J]}\lceil C_j/w_{j,i} \rceil$, on  the number of instantiations of pattern $i$, and this serves as a bounding state at which no further instantiation of this pattern can be added; that is, $\mathscr{N}_i=\{0,1,\ldots,\min_{j\in[J]}\lceil C_j/w_{j,i}\rceil\}$ and $|\mathscr{N}_i|=\min_{j\in[J]}\lceil C_j/w_{j,i}\rceil +1 <+\infty$.
In this context, Equation~\eqref{eqn:constraint:resources} implies the condition
\begin{equation}\label{eqn:constraint:resources:zero}
a_i(\bm{n}) = 0, ~\text{if}~i\notin \{d(\ell):\ \ell\in[L]\} \text{ and } n_i = |\mathscr{N}_i|-1.
\end{equation}


\subsubsection{Objective}
A \emph{policy} $\phi$ is defined as a mapping $\mathscr{N}\rightarrow \mathscr{A}$ where
$\mathscr{A}\coloneqq\prod_{\ell\in[L]}\{0,1\}^{|\mathscr{P}_{\ell}|}$,
determined by the action variables $\bm{a}(\bm{n})$ defined above.
When we are discussing a system operating with a given policy $\phi$, we rewrite the action and state variables as $\bm{a}^{\phi}(\cdot)$ and $\bm{N}^{\phi}(t)$, respectively.

By serving a request of type $\ell\in[L]$ and occupying an RU of pool
$j$ for one unit of time, we gain expected reward $r_{\ell}$ and pay
expected cost $\varepsilon_j$.  
The expected reward for a whole service is gained at the moment the service is completed. 
It corresponds to the situation where a request allocated to pattern $i$ earns reward at rate $r_{\ell(i)} \mu_i$ for as long as it is in the system (so that the expected revenue per customer is $(r_{\ell(i)} \mu_i) . (1/\mu_i) = r_{\ell(i)}$).
The value of $\epsilon_j$ is the cost per unit time of using a unit of capacity from resource pool $j$ in which case the expected cost of accommodating the request in pool $j$ as part of pattern $i$ is $\epsilon_j/\mu_i$.
We seek a policy that
maximizes the \emph{revenue}: the difference between expected reward
and cost, by efficiently utilizing the limited amount of
resources. 

The objective is to maximize the long-run average rate of earning revenue, which exists because, for any policy $\phi$, the process can be modeled by a finite-state Markov chain.
Let $\bm{r}=(r_{\ell}:\ \ell\in[L])$ and $\bm{\varepsilon}=(\varepsilon_{j}:\ j\in[J])$. 
For all $\ell\in[L]$ and $i\in[I]$, define a $L\times I$ matrix $\mathcal{U}$ with  entries
$u_{\ell,i}\coloneqq \mu_{i}\mathds{1}_{i\in\mathscr{P}_{\ell}}$.
By the Strong Law of Large Numbers for Continuous Time Markov Chains, see for example \cite{serfozo2009basics} Theorem 45 in Chapter 4, noting the subsequent discussion of the case where rewards are earned at jump times, the long-run average rate of earning revenue when the policy is $\phi$ is given by
\begin{equation}
R^{\phi} \coloneqq \mathbb{E}_{\pi^{\phi}} \Bigl[\bm{r}\mathcal{U} - \bm{\varepsilon}\mathcal{W}\Bigr] =\sum\limits_{i\in[I]} \sum_{n_i\in\mathscr{N}_i} \pi_i^{\phi}(n_i) \Bigl(r_{\ell(i)} \mu_i  - \sum_{j\in\mathscr{J}} w_{j,i} \varepsilon_j\Bigr)n_i,
\end{equation}
where $\pi_i^{\phi}(n_i)$ is the stationary probability that the state of process $i$ is $n_i$ when the policy is $\phi$. Then we wish to find the policy $\phi$ that maximizes $R^\phi$, that is we wish to find
\begin{equation}\label{eqn:objective}
R \coloneqq \max\limits_{\phi}R^{\phi}.
\end{equation}
Define
$\Phi$ to be the set of all policies
with  the constraints in~\eqref{eqn:constraint:action} and \eqref{eqn:constraint:resources} satisfied.
Each policy in $\Phi$ is then a \emph{feasible policy} for our resource allocation problem.

\section{Whittle Relaxation}\label{sec:relaxation}
Our resource
allocation problem with objective function defined by \eqref{eqn:objective} and constraints given by \eqref{eqn:constraint:action} and \eqref{eqn:constraint:resources} can be modeled as a set of RMABPs coupled by
capacity constraints.
We leave the
specification of the RMABPs to Appendix~\ref{app:RMABPs}.

In this section, we provide a theoretical analysis of the resource allocation problem,
following the idea of Whittle relaxation \cite{whittle1988restless}.
In the vein of a RMABP, 
we randomize the action variable $\bm{a}^{\phi}(\bm{n})$ so that its elements take values from $\{0,1\}$ with probabilities determined by the policy $\phi$ and relax constraint~\eqref{eqn:constraint:action} to require that
\begin{equation}\label{eqn:constraint:relax:action}
\lim\limits_{t\rightarrow +\infty}\mathbb{E}\biggl[\sum\limits_{i\in\mathscr{P}_{\ell}} a^{\phi}_i(\bm{N}^{\phi}(t)) \biggr]= 1,~
\forall \ell\in[L].
\end{equation}
Following similar ideas, we relax \eqref{eqn:constraint:resources} into two equations:
\begin{equation}\label{eqn:constraint:relax:resources}
\lim\limits_{t\rightarrow +\infty}\mathbb{E}\biggl[\mathcal{W}\Bigl(\bm{N}^{\phi}(t)+ \bm{a}^{\phi}(\bm{N}^{\phi}(t))\Bigr)\biggr] \leq \bm{C},
\end{equation}
and  
\begin{equation}\label{eqn:constraint:dummy}
\lim\limits_{t\rightarrow+\infty} \mathbb{E}\Bigl[a^{\phi}_i(\bm{N}^{\phi}(t))\ \mathds{1}_{N^{\phi}_i(t)=|\mathscr{N}_{i}|-1}\Bigr] = 0,~\forall i\in[I]\backslash\{d(\ell):\ \ell\in[L]\}.
\end{equation}
{\bf Remark} 
Equation~\eqref{eqn:constraint:relax:resources} 
is derived by taking expectations for both sides of Equation~\eqref{eqn:constraint:resources},
and \eqref{eqn:constraint:dummy} is a consequence of  \eqref{eqn:constraint:resources:zero}, 
so constraints described by~\eqref{eqn:constraint:relax:resources} and \eqref{eqn:constraint:dummy} are relaxed versions of the constraints described by~\eqref{eqn:constraint:resources}.
The justification for Equation~\eqref{eqn:constraint:dummy} will be discussed in Section~\ref{subsec:asym_regime}, in conjunction with the physical meanings of all variables, when we increase the scale of the entire system.
We refer to the problem with objective~\eqref{eqn:objective}, constraints~\eqref{eqn:constraint:relax:action}, \eqref{eqn:constraint:relax:resources} and \eqref{eqn:constraint:dummy} and randomized control variables $\bm{a}^{\phi}(\bm{n})$, for all $\bm{n}\in\mathscr{N}$, as the \emph{relaxed problem}.

A value $a$ in $(0,1)$ can be interpreted as a randomisation between taking $a_i^{\phi}(\bm{n}) = 0$ and $a_i^{\phi}(\bm{n}) = 1$. Specifically we take $a_i^\phi(n) = 1$ with probability $a$.
We represent the set of policies that correspond to assigning such values $a\in(0,1)$ as $\tilde{\Phi}$.
For $n_i\in \mathscr{N}_i$, $\phi\in\tilde{\Phi}$, $i\in[I]$, define
\begin{itemize}
\item $\alpha^{\phi}_i(n_i) \coloneqq \lim\limits_{t\rightarrow+\infty} \mathbb{E}\left[a^{\phi}_i(\bm{N}^{\phi}(t))\ |\ N^{\phi}_i(t)=n_i\right]$, 
which is the expectation with respect to the stationary distribution when policy $\phi$ is used,
and the vector $\bm{\alpha}^{\phi}_i \coloneqq (\alpha^{\phi}_i(n_i) :\ n_i\in\mathscr{N}_i)$;
\item the stationary probability that $N^{\phi}_i(t)=n_i$ under policy $\phi$ to be $\pi_{i,n_{i}}^{\phi}$, and the vector $\bm{\pi}^{\phi}_i \coloneqq (\pi_{i,n_i}^{\phi}:\ n_i\in\mathscr{N}_i)$.
\end{itemize}
Let $\bm{\Pi}^{\phi}_n \coloneqq \left( \bm{\pi}^{\phi}_{i}\cdot (\mathscr{N}_i):\ i\in[I]\right)^{T}$ and $\bm{\Pi}^{\phi}_a \coloneqq \left(\bm{\pi}^{\phi}_{i}\cdot \bm{\alpha}^{\phi}_i:\ i\in[I]\right)^{T}$, where 
$(\mathscr{N}_i)$ represents the vector $(0,1,\ldots,|\mathscr{N}_i|-1)$.
The Lagrangian function for the optimization problem with objective function (\ref{eqn:objective}) and constraints \eqref{eqn:constraint:relax:action}, \eqref{eqn:constraint:relax:resources} and \eqref{eqn:constraint:dummy} is
\begin{multline}\label{eqn:dual_func}
g(\pmb{\gamma},\bm{\nu},\bm{\eta})\coloneqq
\max\limits_{\phi\in\tilde{\Phi}}
(\bm{r}\mathcal{U}-\bm{\varepsilon}\mathcal{W})\bm{\Pi}^{\phi}_n-\sum\limits_{\ell=1}^{L}\nu_{\ell}\Bigl(\sum\limits_{i\in\mathscr{P}_{\ell}}\bm{\pi}^{\phi}_i\cdot \bm{\alpha}^{\phi}_i -1 \Bigr)\\
- \pmb{\gamma}\cdot \Bigl(\mathcal{W}(\bm{\Pi}^{\phi}_n+\bm{\Pi}^{\phi}_a)-\bm{C}\Bigr)
-\sum\limits_{i\in[I]\backslash\{d(\ell):\ \ell\in[L]\}}\eta_i \pi^{\phi}_{i,|\mathscr{N}_i|-1} \alpha^{\phi}_{i}(|\mathscr{N}_i|-1),
\end{multline}
where 
$\bm{\nu}\in\mathbb{R}^{L}$, $\pmb{\gamma}\in\mathbb{R}_{0}^{J}$ and $\bm{\eta}\in\mathbb{R}^{I-L}$ are Lagrange multiplier vectors for constraints~\eqref{eqn:constraint:relax:action}, \eqref{eqn:constraint:relax:resources} and \eqref{eqn:constraint:dummy}, respectively. 
In~\eqref{eqn:dual_func}, the constraints no longer apply to variables $\bm{\alpha}^{\phi}_i$ ($i\in[I]$) but appear in the maximization as cost items weighted by their Lagrange multipliers.
For $i\in[I]\backslash\{d(\ell):\ \ell\in[L]\}$,
define functions 
\begin{multline}\label{eqn:dual_func_i}
\Lambda_{i}(\phi,\pmb{\gamma},\nu_{\ell(i)},\eta_i)
\coloneqq (r_{\ell(i)}\mu_i-\bm{\varepsilon}\cdot\bm{w}_i)\bm{\pi}^{\phi}_i\cdot (\mathscr{N}_i)-\nu_{\ell(i)}\bm{\pi}^{\phi}_i\cdot \bm{\alpha}^{\phi}_i 
- \pmb{\gamma}\cdot \bigl(\bm{w}_i(\bm{\pi}^{\phi}_i\cdot (\mathscr{N}_i)+\bm{\pi}^{\phi}_i\cdot \bm{\alpha}^{\phi}_i)\bigr)\\
-\eta_i \pi^{\phi}_{i,|\mathscr{N}_i|-1}\alpha^{\phi}_{i}(|\mathscr{N}_i|-1),
\end{multline}
where we recall that $\bm{w}_i$ is the weight vector of pattern $i$ given by the $i$th column vector of $\mathcal{W}$; similarly, for $\ell\in[L]$, $\pmb{\gamma}\in\mathbb{R}_0^J$ and $\eta\in\mathbb{R}$, define
$\Lambda_{d(\ell)}(\phi,\pmb{\gamma},\nu_{\ell},\eta) \coloneqq -\nu_{\ell}\alpha^{\phi}_{d(\ell)}(n)$, where $n$ is the only state in $\mathscr{N}_{d(\ell)}$.
From Equation \eqref{eqn:dual_func}, for $\pmb{\gamma}\in\mathbb{R}_0^J$, $\bm{\nu}\in\mathbb{R}^L$ and $\bm{\eta}\in\mathbb{R}^{I-L}$,
\begin{equation}\label{eqn:dual_func2}
g(\pmb{\gamma},\bm{\nu},\bm{\eta})=\max\limits_{\phi\in\tilde{\Phi}}\sum\limits_{i\in[I]}\Lambda_{i}(\phi,\pmb{\gamma},\nu_{\ell(i)},\eta_i)
+\sum\limits_{\ell\in[L]}\nu_{\ell}
+ \pmb{\gamma}\cdot \bm{C}.
\end{equation}
where $\eta_{d(\ell)}$ $(\ell\in[L])$ are unconstrained real numbers that are used for notational convenience. 


In the maximization problem on the right hand side of~\eqref{eqn:dual_func2}, there is no constraint that restricts the value of one $\Lambda_i(\phi,\pmb{\gamma},\nu_{\ell(i)},\eta_i)$ once the others are known. 
As a result, we can maximize the sum in \eqref{eqn:dual_func2} by maximizing each of the summands independently. 
We can thus write \eqref{eqn:dual_func2} as
\begin{equation}
g(\pmb{\gamma},\bm{\nu},\bm{\eta})=\sum\limits_{i\in[I]}\max\limits_{\phi\in\tilde{\Phi}}\Lambda_{i}(\phi,\pmb{\gamma},\nu_{\ell(i)},\eta_i)
+\sum\limits_{\ell\in[L]}\nu_{\ell}
+ \pmb{\gamma}\cdot \bm{C},
\end{equation}
but with the maximum over $\phi\in\tilde{\Phi}$. 
Observe now that maximizing $\Lambda_i$ over $\phi$ is equivalent to choosing $\bm{\alpha}^{\phi}_i(n_i)$ from $[0,1]^{|\mathscr{N}_i|}$, by interpreting $\alpha^{\phi}_{i,n}\in[0,1]$ as the probability that process $i$ is activated under policy $\phi$ when it is in state $n$. Thus,
\begin{equation}\label{eqn:dual_problem}
g(\pmb{\gamma},\bm{\nu},\bm{\eta})
=\sum\limits_{i\in[I]}\max\limits_{\bm{\alpha}^{\phi}_i\in[0,1]^{|\mathscr{N}_i|}}\Lambda_{i}(\phi,\pmb{\gamma},\nu_{\ell},\eta_i)
+\sum\limits_{\ell\in[L]}\nu_{\ell} + \pmb{\gamma}\cdot \bm{C}.
\end{equation}

By slightly abusing notation, we refer to the policy $\phi$ determined by an action vector $\bm{\alpha}^{\phi}_i$ as the policy for pattern $i$,
and define $\Phi_{i}$ as the set of all policies for pattern $i$.



\begin{definition}\label{def:sub-problem}
The maximization of $\Lambda_{i}(\phi,\pmb{\gamma},\nu_{\ell},\eta_i)$ over $\bm{\alpha}^{\phi}_{i}\in[0,1]^{|\mathscr{N}_{i}|}$ is the \emph{sub-problem} for pattern $i\in[I]$. 
\end{definition}

For given $\pmb{\gamma}$, $\bm{\nu}$ and $\bm{\eta}$, the sub-problem for any pattern is an MDP, so that it can be numerically solved by dynamic programming. By solving the sub-problems for all patterns $i\in[I]$, we obtain $g(\pmb{\gamma},\bm{\nu},\bm{\eta})$.
For any $\pmb{\gamma}$, $\bm{\nu}$ and $\bm{\eta}$, the Lagrangian function $g(\pmb{\gamma},\bm{\nu},\bm{\eta})$ is a performance upper bound for the primal problem described in \eqref{eqn:objective}, \eqref{eqn:constraint:relax:action}, \eqref{eqn:constraint:relax:resources} and \eqref{eqn:constraint:dummy}, which is a relaxed version of the original resource allocation problem. 
Thus there will be a non-negative gap between this upper bound and the maximized performance of the original problem.

\subsection{Analytical Solutions}\label{sec:indexability}
\begin{proposition}\label{prop:a_opt}
For given $\bm{\nu}$ and  $\pmb{\gamma}$, there exists $\bm{E}\in\mathbb{R}^{I-L}$ such that, for any $\bm{\eta}>\bm{E}$, a policy of the sub-problem  for pattern $i$, referred to as $\bar{\varphi}\in\Phi_i$, determined by action vector $\bm{\alpha}^{\bar{\varphi}}_i\in[0,1]^{|\mathscr{N}_i|}$ is optimal for this sub-problem, if, for $n\in\mathscr{N}_i$,
\begin{numcases}{\alpha^{\bar{\varphi}}_{i}(n)}
=1 &\text{if }$0<\lambda_{\ell}(r_{\ell}-\frac{1}{\mu_{i}}\sum\limits_{j\in \mathscr{J}_i}\varepsilon_{j}w_{j,i})
- (1+\frac{\lambda_{\ell}}{\mu_{i}})\sum\limits_{j\in \mathscr{J}_i}w_{j,i}\gamma_{j}-\nu_{\ell} \text{ and } n< |\mathscr{N}_i|-1$,\label{eqn:a_opt:a}\\
\in [0,1]
&\text{if }$0=\lambda_{\ell}(r_{\ell}-\frac{1}{\mu_{i}}\sum\limits_{j\in \mathscr{J}_i}\varepsilon_{j}w_{j,i})
- (1+\frac{\lambda_{\ell}}{\mu_{i}})\sum\limits_{j\in \mathscr{J}_i}w_{j,i}\gamma_{j}-\nu_{\ell} \text{ and } n< |\mathscr{N}_i|-1$,\label{eqn:a_opt:b}\\
=0&\text{otherwise},\label{eqn:a_opt:c}
\end{numcases}
where $\ell=\ell(i)$.
\end{proposition}
The proof will be given in Appendix~\ref{app:prop:a_opt} in the e-companion to this paper.
In the maximization of $\Lambda_i(\phi,\pmb{\gamma},\nu_{\ell(i)},\eta_i)$), the only term of $\Lambda_i$ dependent on $\bm{\eta}$ is $-\eta_i\pi^{\phi}_{i,|\mathscr{N}_i|-1}\alpha^{\phi}_i(|\mathscr{N}_i|-1)$.
The choice of a sufficiently large $\eta_i$ guarantees that $\alpha^{\phi}_i(|\mathscr{N}_i|-1)$ is $0$ for an optimal policy of the sub-problem, so that constraints~\eqref{eqn:constraint:dummy} of the relaxed problem are satisfied. 
For convenience, in what follows we fix $\bm{\eta}$ to be one of these large values 
so that $\alpha^{\phi}_i(|\mathscr{N}_i|-1)$ is also fixed to be $0$ for any optimal policy $\phi$ of the sub-problem for pattern $i$.
By slightly abusing notation, in all subsequent equations and discussions, we directly require $\alpha^{\phi}_i(|\mathscr{N}_i|-1)=0$ ($i\in[I]\backslash\{d(\ell):\ell\in[L]\}$) unless specified otherwise.

{\bf Remark}
Recall that the action variables $\bm{\alpha}^{\phi}_i$ for any pattern $i\in[I]$ and policy $\phi\in\Phi_i$ are potentially state-dependent. However, the right hand sides of equations~\eqref{eqn:a_opt:a}-\eqref{eqn:a_opt:c} are independent of the state variable $n$ which appears on their left hand side, provided that this is less than $|\mathscr{N}_i|-1$. This state-independence phenomenon is a consequence of the linearity of the reward and cost rates in the state variable, $N^{\phi}_i(t)$, for pattern $i\in[I]\backslash\{d(\ell):\ \ell\in[L]\}$: from our definition in Section~\ref{sec:model}, the reward and cost rates of process $i$ in state $N^{\phi}_i(t)$ are $r_{\ell(i)}\mu_{i}N^{\phi}_i(t)$ and $\sum_{j\in\mathscr{J}_i}\varepsilon_j w_{j,i}N^{\phi}_i(t)$, respectively.
A detailed analysis is provided in the proof of Proposition~\ref{prop:a_opt}.


Using an argument similar to that in \cite{whittle1988restless}, we can derive from \eqref{eqn:a_opt:a}-\eqref{eqn:a_opt:c} an abstracted \emph{priority} for each \emph{pattern-state pair} (PS pair) $(i,n)$ with $n\in\mathscr{N}_i\backslash\{|\mathscr{N}_i|-1\}$ and $i\in[I]$; unlike in \cite{whittle1988restless}, here, this priority is $(\pmb{\gamma},\bm{\nu})$-dependent. 
The priority of a PS pair $(i,n)$ with $n\in\mathscr{N}_i\backslash\{|\mathscr{N}_i|-1\}$
is determined by the \emph{index}
\begin{equation}\label{eqn:index_value}
\Xi_{i}(\pmb{\gamma},\bm{\nu})\coloneqq\lambda_{\ell(i)}\Bigl(r_{\ell(i)}-\frac{1}{\mu_{i}}\sum\limits_{j\in \mathscr{J}_i}\varepsilon_{j}w_{j,i}\Bigr)
- \Bigl(1+\frac{\lambda_{\ell(i)}}{\mu_{i}}\Bigr)\sum\limits_{j\in \mathscr{J}_i}w_{j,i}\gamma_{j}-\nu_{\ell(i)},
\end{equation}
and \eqref{eqn:a_opt:a}-\eqref{eqn:a_opt:c} can be characterized as comparing $\Xi_i(\pmb{\gamma},\bm{\nu})$ with $0$. 
When there is strict inequality in the comparison (that is, the cases described in \eqref{eqn:a_opt:a} and \eqref{eqn:a_opt:c}), the value of $\alpha^\phi_i(n)$ is specified, since PS pairs $(i,n)$ for all $n\in\mathscr{N}_i\backslash\{|\mathscr{N}_i|-1\}$ correspond to the same $\Xi_i(\pmb{\gamma},\bm{\nu})$ value. 
However, there is still freedom to decide different values of $\alpha^\phi_i(n)$,
when $\Xi_i(\pmb{\gamma},\bm{\nu})=0$ (the case described in \eqref{eqn:a_opt:b}).
A detailed discussion about priorities of PS pairs corresponding to the same $\Xi_i(\pmb{\gamma},\bm{\nu})$ will be provided in Section~\ref{subsec:decomposable}.
By solving the sub-problem of dummy pattern $d(\ell)$ $(\ell\in[L])$ which involves only one state $n\in\mathscr{N}_{d(\ell)}$, we obtain an optimal policy $\varphi$ determined by
\vspace{-0.1cm}
\begin{equation}\label{eqn:a_opt:dummy}
\alpha^{\varphi}_{d(\ell)}(n) \left\{\begin{array}{ll}
=1, &\text{if } 0 < -\nu_{\ell}, \vspace{-0.1cm}\\
\in[0,1], &\text{if }0 = -\nu_{\ell},\vspace{-0.1cm}\\
=0, &\text{otherwise}. 
\end{array}\right.
\end{equation}
\vspace{-0.1cm}
The priority of the state of a dummy pattern is then $\Xi_{d(\ell)}(\pmb{\gamma},\bm{\nu}) \equiv -\bm{\nu}$ for any $\pmb{\gamma}$. 

In addition, from Equation~\eqref{eqn:a_opt:c} in Proposition~\ref{prop:a_opt}, for any given $\bm{\nu}\in\mathbb{R}^I$ and $\pmb{\gamma}\in\mathbb{R}_0^J$, there exists $\bm{\eta}\in\mathbb{R}^{I-L}$ such that it is optimal to make states $|\mathscr{N}_i|-1$ passive (that is, $\alpha^{\bar{\varphi}}_i (|\mathscr{N}_i| - 1) = 0$) for all $i\in[I]\backslash\{d(\ell):\ \ell\in[L]\}$. 
Among all PS pairs $(n,i)$ ($n\in\mathscr{N}_i$, $i\in[I]$), we assign, without loss of generality, the least  priority to  those PS pairs $(i,|\mathscr{N}_i|-1)$ for which $i\in[I]\backslash\{d(\ell):\ \ell\in[L]\}$.

The policy $\bar{\varphi}$ determined by \eqref{eqn:a_opt:a}-\eqref{eqn:a_opt:c} and \eqref{eqn:a_opt:dummy} is optimal for the relaxed problem described by \eqref{eqn:objective}, \eqref{eqn:constraint:relax:action}, \eqref{eqn:constraint:relax:resources} and \eqref{eqn:constraint:dummy}, if the given multipliers $\bm{\nu}$ and $\pmb{\gamma}$ that appear in \eqref{eqn:a_opt:a}-\eqref{eqn:a_opt:c} and \eqref{eqn:a_opt:dummy}  satisfy the \emph{complementary slackness conditions} of this relaxed problem, defined by \vspace{-0.3cm}
\begin{condition}{Complementary Slackness}\label{cond:complementary_slackness}
\vspace{-0.3cm}
\begin{equation}
\label{eqn:relaxaction:slack}
\nu_{\ell}\Bigl(\sum\limits_{i\in\mathscr{P}_l} \bm{\pi}^{\phi}_i\cdot\bm{\alpha}^{\phi}_i -1\Bigr)=0,~
\forall l\in[L],\vspace{-0.3cm}
\end{equation}
and \vspace{-0.3cm}
\begin{equation}
\label{eqn:relaxconstraint:slack}
\gamma_{j} \Bigl(\bm{\omega}_j\cdot \left(\bm{\Pi}^{\phi}_n+\bm{\Pi}^{\phi}_a\right)-C_j\Bigr) = 0,~\forall j\in[J], \vspace{-0.3cm}
\end{equation}
where $\bm{\omega}_j=(w_{j,i}:\ i\in[I])$ is the $j$th row of matrix $\mathcal{W}$.\vspace{-0.3cm}
\end{condition}
In this context, if resource pool $j\in[J]$ is very popular so that the capacity constraint corresponding to the $j$th row in  \eqref{eqn:constraint:relax:resources} achieves equality, then $\gamma_j$ is allowed to be positive, leading to a lower value of $\Xi_i(\pmb{\gamma},\bm{\nu})$ than for $\gamma_j=0$.
On the other hand, if resource pool $j\in[J]$ cannot be fully occupied and the $j$th capacity constraint in \eqref{eqn:constraint:relax:resources} is satisfied with a strict inequality, then the 
complementary slackness condition described in \eqref{eqn:relaxconstraint:slack} forces $\gamma_j$ to be zero.
Following this mechanism, when resource pool $j\in[J]$ is overloaded and its priority is reduced by increasing $\gamma_j$, the offered traffic to this resource pool will be reduced in line with its priority.

If
there exist multipliers $\bm{\nu}$, $\pmb{\gamma}$ and a policy $\bar{\varphi}$ determined by \eqref{eqn:a_opt:a}-\eqref{eqn:a_opt:c}, such that  the complementary slackness conditions~\eqref{eqn:relaxaction:slack} and \eqref{eqn:relaxconstraint:slack} are satisfied by taking $\phi=\bar{\varphi}$, then, by the strong duality theorem, this policy $\bar{\varphi}$ is optimal for the relaxed problem;
this observation, together with Theorem~\ref{theorem:main_second} in Section~\ref{sec:asym_opt}, leads to the existence of an asymptotically optimal policy feasible for the original problem, derived with priorities of patterns induced by the descending order of $\Xi_i(\pmb{\gamma},\bm{\nu})$. More details about the analysis in the asymptotic regime will be provided in Section~\ref{sec:asym}.
Here we focus on the non-asymptotic regime, and specifically on the choice and computation of $\pmb{\gamma}$ and $\bm{\nu}$.

\subsection{Decomposable Capacity Constraints}\label{subsec:decomposable}

In the general case, it is not clear whether the 
the complementary slackness conditions~\eqref{eqn:relaxaction:slack} and \eqref{eqn:relaxconstraint:slack}
can be satisfied and, 
even if they are, what the values of $\pmb{\gamma}$ and the corresponding $\bm{\nu}$ are.
More important is the question of how the multipliers help with proposing the asymptotically optimal policy applicable to the original problem.

In Sections~\ref{subsec:decomposable} and \ref{subsec:sufficient_condition}, we concentrate on the complementary slackness conditions and the existence of an optimal policy (for the relaxed problem) satisfying \eqref{eqn:a_opt:a}-\eqref{eqn:a_opt:c}.
Recall that \eqref{eqn:a_opt:a}-\eqref{eqn:a_opt:c} intuitively suggest priorities of patterns induced by $\Xi_i(\pmb{\gamma},\bm{\nu})$. Later in Section~\ref{sec:index_policy}, a policy feasible for the original problem will be proposed based on given priorities of patterns, and its asymptotic optimality will be discussed in Section~\ref{sec:asym_opt}.

\subsubsection{Priorities of PS Pairs}\label{subsubsec:priorities_of_pairs}

As described in Section~\ref{sec:indexability}, 
the priorities of PS pairs are determined by the descending order of $\Xi_i(\pmb{\gamma},\bm{\nu})$, with higher priorities given by higher values of $\Xi_i(\pmb{\gamma},\bm{\nu})$. 
It may happen that, because of different tie-breaking rules, the same $\pmb{\gamma}$ and $\bm{\nu}$ lead to different priorities. For given $\pmb{\gamma}\in\mathbb{R}_0^J$ and $\bm{\nu}\in\mathbb{R}^L$, let $\mathscr{O}(\pmb{\gamma},\bm{\nu})$ represent the set of all rankings of PS pairs according to the descending order of $\Xi_i(\pmb{\gamma},\bm{\nu})$ ($i\in[I]$). 
Also, for notational convenience, let $\mathscr{O}$ represent the set of all PS pair rankings.

To emphasize the priorities of these PS pairs, according to a given ranking $o\in\mathscr{O}$,  we label all these pairs by their order $\iota^{o}\in[N]$ with $N\coloneqq \sum_{i\in[I]} |\mathscr{N}_i|$ and $(i_{\iota^o},n_{\iota^o})$ giving the pattern and the state of the $\iota^{o}$th PS pair.
We will omit the superscript $o$ and use $\iota$ for notational simplicity unless it is necessary to specify the underlying ranking.
There exists one and only one $\ell\in[L]$ satisfying $i_{\iota} \in \mathscr{P}_{\ell}$ for any PS pair labeled by $\iota$. Such an $\ell$ is denoted by $\ell_{\iota}$.

\IncMargin{1em}
\begin{algorithm}
\small 
\linespread{0.4}\selectfont

\SetKwFunction{FPriorityPolicy}{PriorityPolicy}
\SetKwProg{Fn}{Function}{:}{\KwRet}
\SetKwInOut{Input}{Input}\SetKwInOut{Output}{Output}
\SetAlgoLined
\DontPrintSemicolon
\Input{a vector of non-negative reals $\pmb{\gamma}\in\mathbb{R}_0^J$ and a ranking of PS pairs $o\in\mathscr{O}$.}
\Output{a policy $\bar{\varphi}(o)\in\tilde{\Phi}$ determined by action variables $\bm{\alpha}^{\bar{\varphi}(o)}_i\in[0,1]^{|\mathscr{N}_i|}$ for all $i\in[I]$ and a vector of reals $\bm{\nu}(o,\pmb{\gamma})$.}

\Fn{\FPriorityPolicy{$o,\pmb{\gamma}$}}{

    $\bm{\alpha}^{\bar{\varphi}}_i\gets \bm{0}$ for all $i\in[I]$ \tcc*{Variables $\bm{\alpha}^{\bar{\varphi}}_i$ determine a policy $\bar{\varphi}$}
    Initializing the list of candidate PS pairs as the list of all PS pairs\;
    $\iota \gets 0$ \tcc*{Iteration variable} 
    \While {$\iota<N $ {\bf and} the list of candidate PS pairs is not empty}{	
	 	    $\iota\gets \iota+1$\;
			\If {PS pair $\iota$ is not in the list of candidate PS pairs}{
			    {\bf continue}\;
			}
	    $a_1 \gets \inf \left\{\left\{\alpha^{\bar{\varphi}}_{i_{\iota}}(n_{\iota}) \in [0,1] \left|~\sum_{i\in\mathscr{P}_{\ell_{\iota}}}\bm{\pi}^{\bar{\varphi}}_i\cdot \bm{\alpha}^{\bar{\varphi}}_i = 1\right.\right\}\cup\{1\}\right\}$\;

	    \tcc*{The maximal probability of activating PS pair $\iota$ such that }
		\tcc*{the action constraint for RT $\ell_{\iota}$ is not violated.}

  	    $a_2 \gets \inf \left\{\left\{\alpha^{\bar{\varphi}}_{i_{\iota}}(n_{\iota}) \in [0,1] \left|~\exists j\in[J], \bm{\omega}_j\cdot\left(\bm{\Pi}^{\bar{\varphi}}_n+\bm{\Pi}^{\bar{\varphi}}_a\right)=C_j\right.\right\}\cup\{1\}\right\}$\;

		\tcc*{The maximal probability of activating PS pair $\iota$ such that}
	 	\tcc*{the capacity constraints are not violated.}

		$\alpha^{\bar{\varphi}}_{i_{\iota}}(n_{\iota})\gets \min\{a_1,a_2\}$ \;
		\tcc*{Update $\alpha^{\bar{\varphi}}_{i_{\iota}}(n_{\iota})$ with the maximal activating probability}
		\tcc*{without violating any constraint.}

		\uIf{$\sum_{i\in\mathscr{P}_{\ell_{\iota}}}\bm{\pi}^{\bar{\varphi}}_i\cdot \bm{\alpha}^{\bar{\varphi}}_i = 1$}{
			\tcc*{If the action constraint achieves equality under policy $\bar{\varphi}$}
			\tcc*{determined by updated $\bm{\alpha}^{\bar{\varphi}}_i$, $i\in[I]$.}
		    $\nu_{\ell_{\iota}}(o,\pmb{\gamma})\gets\Xi_{i_{\iota}}(\pmb{\gamma},\bm{0})$\;

	       remove all PS pairs $\iota'>\iota$ with $\ell_{\iota'}=\ell_{\iota}$ from the list of candidate PS pairs\;
	   	}
   		\ElseIf{$\exists j\in[J],\ \bm{\omega}_j\cdot\left(\bm{\Pi}^{\bar{\varphi}}_n+\bm{\Pi}^{\bar{\varphi}}_a\right)=C_j$}{
			\tcc*{If a capacity constraint achieves equality under policy $\bar{\varphi}$}
			\tcc*{determined by updated $\bm{\alpha}^{\bar{\varphi}}_i$, $i\in[I]$.}

		    remove all PS pairs $\iota'>\iota$ with $w_{j,i_{\iota'}}>0$ from the list of candidate PS pairs\;
		}

}

$\bm{\alpha}^{\bar{\varphi}(o)}_i\gets \bm{\alpha}^{\bar{\varphi}}_{i}$ for all $i\in[I]$\;

}
\caption{Priority-style policy for the relaxed problem}\label{algo:varphi_gamma}
\end{algorithm}
 \DecMargin{1em}

For any given ranking of PS pairs $o\in\mathscr{O}$, we can generate a policy $\bar{\varphi}(o)$ with priorities of PS pairs defined by $o$, such that 
\eqref{eqn:constraint:relax:action}, \eqref{eqn:constraint:relax:resources} and \eqref{eqn:constraint:dummy} are satisfied:
the policy $\bar{\varphi}(o)$ is feasible for the relaxed problem but not necessarily feasible for the original problem.
The pseudo-code for generating $\bar{\varphi}(o)$ is presented in Algorithm~\ref{algo:varphi_gamma}.
The key idea for generating such a $\bar{\varphi}(o)$ is to initialize $\bm{\alpha}_i^{\bar{\varphi}(o)}$ to $\bm{0}$ for all $i\in[I]$, and sequentially activate the PS pairs according to their priorities defined by $o$ until either a relaxed action or capacity constraint described in \eqref{eqn:constraint:relax:action} and \eqref{eqn:constraint:relax:resources}, respectively, achieves equality. In particular, 
\begin{enumerate}[label=(\Roman*)]
\item
if a relaxed action constraint described in \eqref{eqn:constraint:relax:action} achieves equality by activating PS pairs less than or equal to $\iota$, then the multiplier $\nu_{\ell_{\iota}}$ is set to $\Xi_{i_{\iota}}(\pmb{\gamma},\bm{0})$, and all later PS pairs $\iota'>\iota$ with $\ell_{\iota'}=\ell_{\iota}$ are \emph{disabled} from being activated and are removed from the \emph{list of candidate pairs} awaiting later activation; \label{case:equality:action}
\item similarly, if a relaxed capacity constraint described in \eqref{eqn:constraint:relax:resources} associated with resource pool $j\in[J]$ achieves equality by activating PS pairs less than or equal to $\iota$, then all later PS pairs $\iota'>\iota$ with $w_{j,i_{\iota'}}>0$ are disabled and removed from the list of candidate states.\label{case:equality:capacity}
\end{enumerate}
Maintaining an iteratively updated list of candidate pairs in this way continues until all 
action constraints in \eqref{eqn:constraint:relax:action} achieve equality: the policy $\bar{\varphi}(o)$ is determined by the resulting $\bm{\alpha}^{\bar{\varphi}(o)}_i$ ($i\in[I]$), and the multipliers $\bm{\nu}$ are updated in \ref{case:equality:action}.
The vector of these multipliers is denoted by $\bm{\nu}(o,\pmb{\gamma})$.
The PS pair labeled by $\iota$ satisfying the condition described in \ref{case:equality:capacity} is called the \emph{critical pair}, with the corresponding resource pool $j$ referred to as the \emph{critical pool} of PS pair $\iota$, denoted by $j_{\iota}(o)$. Note that, from the description in \ref{case:equality:capacity}, there might be more than one resource pool for which the capacity constraints achieve equality simultaneously while activating PS pair $\iota$; we choose one of them to be $j_{\iota}(o)$ and refer to this resource pool as the critical pool of $\iota$.
Let $\mathscr{I}(o)$ represent the set of all critical pairs with respect to the policy $\bar{\varphi}(o)$. 
\begin{lemma}\label{lemma:critical_pattern}
For any $o\in\mathscr{O}$ and $\iota,\iota'\in\mathscr{I}(o)$, if $\iota \neq \iota'$ then $i_{\iota}\neq i_{\iota'}$.
\end{lemma}
\proof{Proof.}
Consider critical pairs $\iota,\iota'\in\mathscr{I}(o)$ with $\iota\neq \iota'$, and assume $\iota < \iota'$ without loss of generality. Since $\iota$ is a critical pair, there is a critical resource pool $j_{\iota}$ which is fully occupied. In this case, if $i_{\iota} = i_{\iota'}$, then pair $\iota'$ must require some resource units from pool $j_{\iota}$ and so $\alpha^{\bar{\varphi}(o)}_{\iota'}=0$.
PS pair $\iota'$ cannot be critical, which violates the condition $\iota'\in\mathscr{I}(o)$. 
Hence, $i_{\iota} \neq i_{\iota'}$. This proves the lemma.
\endproof

Recall, for any ranking $o$, the policy $\bar{\varphi}(o)$ must satisfy the action and capacity constraints~\eqref{eqn:constraint:relax:action}, \eqref{eqn:constraint:relax:resources} and \eqref{eqn:constraint:dummy}. Also, since \eqref{eqn:constraint:relax:action} holds, the 
complementary slackness conditions 
corresponding to the action constraints~\eqref{eqn:relaxaction:slack} are satisfied by taking $\phi=\bar{\varphi}(o)$. However, the complementary slackness conditions 
corresponding to the capacity constraints~\eqref{eqn:relaxconstraint:slack} and equations~\eqref{eqn:a_opt:a}-\eqref{eqn:a_opt:c} are not necessarily satisfied if we plug in $\phi=\bar{\varphi}(o)$ and $\pmb{\gamma}$: the policy $\bar{\varphi}(o)$ is a heuristic policy  applicable for the relaxed problem defined by \eqref{eqn:objective}, \eqref{eqn:constraint:relax:action}, \eqref{eqn:constraint:relax:resources} and \eqref{eqn:constraint:dummy} derived by intuitively prioritizing PS pairs according to their ranking $o\in\mathscr{O}$.

In Section~\ref{subsec:sufficient_condition} we shall define a particular class of resource allocation models, for which we can show the complementary slackness conditions are indeed satisfied.

\begin{definition}\label{define:decomposable}
The system said to be \emph{decomposable} if there exist multipliers $\pmb{\gamma}\in\mathbb{R}_0^J$, $\bm{\nu}\in\mathbb{R}^L$ and a ranking $o\in\mathscr{O}(\pmb{\gamma},\bm{\nu})$ such that $\bm{\nu}=\bm{\nu}(o,\pmb{\gamma})$ and the complementary slackness conditions~\eqref{eqn:relaxaction:slack} and \eqref{eqn:relaxconstraint:slack} are achieved by taking $\phi=\bar{\varphi}(o)$. In this case the optimal values of the dual variables are called \emph{decomposable values}. 
\end{definition}

Recall that, in the general case, for $\pmb{\gamma}\in\mathbb{R}_0^J$ and $\bm{\nu}\in\mathbb{R}^L$, even if $o\in\mathscr{O}(\pmb{\gamma},\bm{\nu})$, the policy $\bar{\varphi}(o)$ is not necessarily optimal (because it does not necessarily satisfy \eqref{eqn:a_opt:a}-\eqref{eqn:a_opt:c}). 
When the policy $\bar{\varphi}(o)$ is optimal for the relaxed problem, the ranking $o$ can be used to construct an index policy applicable to the original problem (detailed steps are provided in Section~\ref{sec:index_policy}).  Theorem~\ref{theorem:main_second} (in Section~\ref{sec:asym_opt}) then ensures that such an index policy is asymptotically optimal.

\subsubsection{Derivation of the Pair Ranking} \label{subsubsec:derive_ranking}

We start with a proposition that shows how the values of the Lagrange multipliers $\bm{\nu}$ and $\pmb{\gamma}$ can be derived from a knowledge of the critical pair and critical resource pool corresponding to a given order $o\in \mathscr{O}$.

\begin{proposition}\label{prop:solution_existence}
For any given $\pmb{\gamma}_0\in\mathbb{R}_0^J$ and $o\in\mathscr{O}$, the linear equations
\begin{equation}\label{eqn:necessary_gamma}
\nu_{\ell_{\iota}}(o,\pmb{\gamma}_0)=\Xi_{i_{\iota}}(\pmb{\gamma},\bm{0}),~\forall \iota\in\mathscr{I}(o)
\vspace{-0.3cm}
\end{equation}
and \vspace{-0.3cm}
\begin{equation}\label{eqn:necessary_gamma:zero}
\gamma_{j} = 0,~\forall j \notin \{j_{\iota}(o)\in[J] ~|~\iota\in\mathscr{I}(o)\}
\end{equation}
have a unique solution $\pmb{\gamma}\in\mathbb{R}^J$.
\end{proposition}
The proof of Proposition~\ref{prop:solution_existence} will be given in Appendix~\ref{app:prop:solution_existence} in the e-companion. 
For a ranking $o\in\mathscr{O}$, define an function $\mathcal{T}^o$ of $\pmb{\gamma}_0\in\mathbb{R}_0^J$ with respect to $o\in\mathscr{O}$: $\mathcal{T}^o(\pmb{\gamma}_0)\coloneqq \pmb{\gamma}$ where $\pmb{\gamma}$ is the unique solution of \eqref{eqn:necessary_gamma} and \eqref{eqn:necessary_gamma:zero}. 
Let $\mathcal{T}^o_j(\pmb{\gamma}_0)$ represent the $j$th element of $\mathcal{T}^o(\pmb{\gamma}_0)$.

\begin{proposition}\label{prop:converge_gamma}
If there exist $\pmb{\gamma}_0\in\mathbb{R}_0^J$ and $o\in\mathscr{O}(\pmb{\gamma}_0,\bm{0})$ such that $\mathcal{T}^o(\pmb{\gamma}_0)=\pmb{\gamma}_0$, then  
$\pmb{\gamma}_0$ is a vector of decomposable multipliers and the policy $\bar{\varphi}(o)$ based on  ranking $o$ is optimal for the relaxed problem defined by \eqref{eqn:objective}, \eqref{eqn:constraint:relax:action}, \eqref{eqn:constraint:relax:resources} and \eqref{eqn:constraint:dummy}.
\end{proposition}
The proof of Proposition~\ref{prop:converge_gamma} will be given in Appendix~\ref{app:prop:converge_gamma} in the e-companion.  Recall that $\mathscr{I}(o)$ is the set of critical pairs with respect to the policy $\bar{\varphi}(o)$, $j_{\iota}(o)$ is the critical resource pool corresponding to critical pair $\iota\in\mathscr{I}(o)$ according to ranking $o$, and $\nu_{\ell_{\iota}}(o,\pmb{\gamma}_0)$ is an output of Algorithm~\ref{algo:varphi_gamma} when the inputs are $o$ and $\pmb{\gamma}=\pmb{\gamma}_0$.

{\bf Remark} Proposition~\ref{prop:converge_gamma} provides a way of checking decomposability of $\pmb{\gamma}_0$ and optimality of $\bar{\varphi}(o)$.
By Proposition~\ref{prop:converge_gamma}, any fixed point $\pmb{\gamma}_0\in\mathbb{R}_0^J$ of the function $\mathcal{T}^o$ with respect to a ranking $o\in\mathscr{O}(\pmb{\gamma}_0,\bm{0})$ is a decomposable vector. 
The decomposability of $\pmb{\gamma}_0$ can be checked without requiring knowledge of any $\bm{\nu}\in\mathbb{R}^L$. 
Also, we present the following corollary of Proposition~\ref{prop:converge_gamma}.
\begin{corollary}\label{coro:converge_gamma}
For $\pmb{\gamma}_0\in\mathbb{R}_0^J$ and $o\in\mathscr{O}(\pmb{\gamma}_0,\bm{0})$, if $\mathcal{T}^o(\pmb{\gamma}_0)\neq \pmb{\gamma}_0$, $\mathcal{T}^o(\pmb{\gamma}_0)\in \mathbb{R}_0^J$ and $o\in\mathscr{O}(\mathcal{T}^o(\pmb{\gamma}_0),\bm{0})$, then $\mathcal{T}^o(\mathcal{T}^o(\pmb{\gamma}_0))=\mathcal{T}^o(\pmb{\gamma}_0)$. 
\end{corollary}
Note that the hypothesis of Corollary~\ref{coro:converge_gamma} requires all components of $\mathcal{T}^o(\pmb{\gamma}_0)$ to be nonnegative. This is not such an easy condition to satisfy. 
The proof of Corollary~\ref{coro:converge_gamma} will be given in Appendix~\ref{app:coro:converge_gamma} in the e-companion.

In this context, consider a given $\pmb{\gamma}_0\in\mathbb{R}_0^J$ and a ranking $o\in\mathscr{O}(\pmb{\gamma}_0,\bm{0})$. If $\pmb{\gamma}_0$ is a fixed point of $\mathcal{T}^o$, then it is the vector of decomposable multipliers; if it is not but $\mathcal{T}^o(\pmb{\gamma}_0)$ is a nonnegative fixed point of $\mathcal{T}^o$, then $\mathcal{T}^o(\pmb{\gamma}_0)$ represents the decomposable multipliers. 
However, in both cases we need to propose a specific $\pmb{\gamma}_0$; it requires prior knowledge of the multipliers, which is, in general, not available.
Section~\ref{subsec:sufficient_condition} will discuss a special case where the decomposability is provable and we have a method of deriving the decomposable multipliers.
Here,  to make a reasonably good choice of the Lagrangian multipliers in a general system, we embark on a \emph{fixed point iteration method}.

Since Proposition~\ref{prop:converge_gamma} requires a fixed point $\pmb{\gamma}$ of the function $\mathcal{T}^o$ with $o\in\mathscr{O}(\pmb{\gamma},\bm{0})$, we need to iterate not only the value of $\pmb{\gamma}$ but also the corresponding ranking $o$ which affects the function $\mathcal{T}^o$ and should be an element of $\mathscr{O}(\pmb{\gamma},\bm{0})$.
Following the idea of conventional
fixed point interation methods, for $k\in\mathbb{N}_0$, 
let $\pmb{\gamma}_{k+1}= \bigl(\mathcal{T}^{o_k}(\pmb{\gamma}_k)\bigr)^+$ with initial $\pmb{\gamma}_0$ and $o_0\in\mathscr{O}(\pmb{\gamma}_0,\bm{0})$, 
where $(\bm{v})^+\coloneqq (\max\{0,v_i\}: i\in[N])$ for a vector $\bm{v}\in \mathbb{R}^N$ ($N\in\mathbb{N}_+$). 
Construct a ranking $o_{k+1}\in\mathscr{O}(\pmb{\gamma}_{k+1},\bm{0})$ according to $o_k$:
for any two different PS pairs $(i,n)$ and $(i',n')$ with $\Xi_i(\pmb{\gamma}_{k+1},\bm{0})=\Xi_{i'}(\pmb{\gamma}_{k+1},\bm{0})$, $(i,n)$ precedes $(i',n')$ in the ranking $o_{k+1}$ if and only if $(i,n)$ precedes $(i',n')$ in the ranking $o_k$.
Here, the operation $(\cdot)^+$ is used to make all the elements of $\pmb{\gamma}_{k+1}$ non-negative, so that $\pmb{\gamma}_{k+1}$ is feasible for the function $\mathcal{T}^{o_{k+1}}$.
Thus the ranking $o_{k+1}$ inherits the tie-breaking rule used for $o_k$ so that the difference between $o_k$ and $o_{k+1}$, which must satisfy $o_k\in\mathscr{O}(\pmb{\gamma}_k,\bm{0})$ and $o_{k+1}\in\mathscr{O}(\pmb{\gamma}_{k+1},\bm{0})$, is minimized.
Corollary~\ref{coro:converge_gamma} can be used to check whether the $\pmb{\gamma}_{k+1}$ is a fixed point of the function $\mathcal{T}^{o_k}$.
Also, $\pmb{\gamma}_{k+1}$ and $o_{k+1}$ are uniquely determined by $\pmb{\gamma}_k$ and $o_k$. We can consider $(\pmb{\gamma}_k,o_k)$ as an entity which is an argument delivered to the function $\mathcal{T}^{o_k}(\pmb{\gamma}_k)$, and wish to find a fixed point in this sense.

In the general case, the function $\mathcal{T}^{o_k}(\pmb{\gamma}_k)$ is discontinuous in $\pmb{\gamma}_k$ and the sequence $\{\pmb{\gamma}_k\}_{k=0}^{\infty}$ is heuristically generated with no proof of convergence to a fixed point. 
In fact, the choice of $\pmb{\gamma}_{k+1}= \bigl(\mathcal{T}^{o_k}(\pmb{\gamma}_k)\bigr)^+$ may result in the sequence $\{\pmb{\gamma}_k\}_{k=0}^{\infty}$ being trapped in oscillations. 
To avoid this,  with slight abuse of notation, we modify the iteration as $\pmb{\gamma}_{k+1}= \bigl(c\mathcal{T}^{o_k}(\pmb{\gamma}_k)+(1-c)\pmb{\gamma}_k\bigr)^+$ with a parameter $c\in [0,1]$, which captures the effects of exploring the new point $\mathcal{T}^{o_k}(\pmb{\gamma}_k)$.
Numerical examples of iterating $\pmb{\gamma}_k$ will be provided in Section~\ref{sec:example}.

With an upper bound, $U\in\mathbb{N}_+$, we take $k^*\coloneqq \arg\min_{k=1,2,\ldots,U} \lVert \pmb{\gamma}_{k-1}-\pmb{\gamma}_k\rVert$ and consider $o_{k^*}$ as a reasonably good ranking of PS pairs. 
Such $o_{k^*}$ is pre-computable with computational complexity no worse than $O(U(N^2+J^2))$, where $N^2$ and $J^2$ result from ordering the $N$ pairs and solving the $J$ linear equations, respectively.
In Section~\ref{sec:index_policy}, we show that an index policy feasible for the original problem can always be generated with such an $o_{k^*}$, and the implementation complexity is $O(I)$  in terms of computation and storage.

\vspace{-0.3cm}
\subsection{Weakly Coupled Constraints}\label{subsec:sufficient_condition}
\vspace{-0.3cm}
Here, we discuss a sufficient condition under which the sequence $\{\pmb{\gamma}_k\}_{k=0}^{\infty}$ is provably convergent; and, in Section~\ref{sec:example}, when this condition fails, we show via numerical examples that the sequence might still converge.

\begin{definition}
Recall the matrix $\mathcal{W}=(w_{j,i})$ defined in Section~\ref{subsec:model}.
We say that row $j\in[J]$ is 
\begin{enumerate}
\item a \emph{type-1} row if there is at most one $i\in[I]$ with $w_{j,i}>0$; \vspace{-0.15cm}
\item a \emph{type-2} row if there is more than one $i\in[I]$ with $w_{j,i}>0$.
\end{enumerate} \vspace{-0.3cm}
\end{definition}
That is, row $j$ is a type-1 row if resource pool $j$ is not shared by patterns of different types; and is a type-2 row, otherwise.
Denote by $\mathscr{J}_i=\{j\in[J]\ |\ w_{j,i}>0\}$ the set of resource pools used by pattern $i$. 
We then define a condition. \vspace{-0.3cm}
\begin{condition}{Weak Coupling}\label{cond:hypothesis}
A system is weakly coupled if, for any pattern $i$, there is at most one $j\in\mathscr{J}_i$ with row $j$ of $\mathcal{W}$ being a type-2 row. \vspace{-0.3cm}
\end{condition}
This condition implies that there is at most one shared resource pool associated with each pattern. 
In a weakly coupled system, if pattern $i_1$ shares a resource pool $j_{12}$ with pattern $i_2$ and pattern $i_1$ shares a resource pool $j_{13}$ with pattern $i_3$ then $j_{12}=j_{13}$. 
A system where each of the patterns requires only one resource pool is clearly weakly coupled.
Note that, in a  weakly coupled system, dependencies between state variables of different patterns still exist, because 
each resource pool can be shared by requests of multiple RTs.

\begin{definition}\label{define:w_star}
For a weakly coupled system define, for each $i\in[I]\backslash\{d(\ell):\ell\in[L]\}$, $w^*_i=w_{j,i}$ where $j$ is the only resource pool in $\mathscr{J}_i$ shared with other patterns, if there is one; or any member of the set $\arg\min\limits_{j'\in\mathscr{J}_i}\frac{C_{j'}}{w_{j',i}}$, otherwise.
\end{definition}
\begin{definition}\label{define:o_star}

For a weakly coupled system define, for $\bm{\nu}\in\mathbb{R}^L$, a set of PS rankings $\mathscr{O}^*(\bm{\nu})\subset \mathscr{O}$ such that, for any $o\in\mathscr{O}^*(\bm{\nu})$, PS pairs $\iota\in[N]$ are ranked according to the descending order of 
\begin{equation}\label{eqn:main:index}
\Xi_{\iota}^*= \left\{
\begin{array}{ll}
\frac{\Xi_{i_{\iota}}(\bm{0},\bm{0})-\nu_{\ell_{\iota}}}{w^*_{i_{\iota}}(1+\lambda_{\ell_{\iota}}/\mu_{i_{\iota}})}, &\text{if } \nexists \ell\in[L],~i_{\iota}=d(\ell),\\
0, & \text{otherwise},
\end{array}\right.
\end{equation} 
\end{definition}

\begin{proposition}\label{prop:equal_opt}
If the system is weakly coupled and there exists a ranking $o\in\mathscr{O}^*(\bm{0})$ satisfying $\nu(o,\bm{0})=\bm{0}$, then the capacity constraints described in \eqref{eqn:constraint:relax:resources} are decomposable and the policy $\bar{\varphi}(o)$ is optimal for the relaxed problem defined by \eqref{eqn:objective} and \eqref{eqn:constraint:relax:action}-\eqref{eqn:constraint:dummy}.
In particular, there exist decomposable multipliers $\pmb{\gamma}\in\mathbb{R}_0^{J}$ satisfying, for $j\in[J]$,
\begin{enumerate}[label=\roman*)]
\item if there is a critical PS pair $\iota\in\mathscr{I}(o)$ with critical resource pool $j=j_{\iota}(o)$, and no 
$j'\neq j$ with $j'\in\mathscr{J}_{i_{\iota}}$ is critical for any other PS pair $\iota' \in\mathscr{I}(o)$, then
\vspace{-0.3cm}
\begin{equation}\label{eqn:equal_opt:a}
\gamma_j = \frac{\Xi_{i_{\iota}}(\bm{0},\bm{0})-\nu_{\ell_i}}{w_{j,i_{\iota}}\left(1+\lambda_{\ell_{\iota}}/\mu_{i_{\iota}}\right)};
\vspace{-0.3cm}
\end{equation}
\item if there are critical PS pairs $\iota$ and $\iota'$ in $\mathscr{I}(o)$ with critical resource pools $j=j_{\iota}(o)\neq j_{\iota'}(o)$ and $ j_{\iota'}(o)\in\mathscr{J}_{i_{\iota}}$, then  \vspace{-0.3cm}
\begin{equation}\label{eqn:equal_opt:b}
\gamma_j=\frac{w_{j_{\iota'}(o),i_{\iota}}}{w_{j,i_{\iota}}}\left(\frac{\Xi_{i_{\iota}}(\bm{0},\bm{0})-\nu_{\ell_{\iota}}}{w_{j_{\iota'}(o),i_{\iota}}\left(1+\lambda_{\ell_{\iota}}/\mu_{i_{\iota}}\right)}
-\frac{\Xi_{i_{\iota'}}(\bm{0},\bm{0})-\nu_{\ell_{\iota'}}}{w_{j_{\iota'}(o),i_{\iota'}}\left(1+\lambda_{\ell_{\iota'}}/\mu_{i_{\iota'}}\right)}\right);\vspace{-0.3cm}
\end{equation}
\item otherwise, \vspace{-0.3cm}
\begin{equation}\label{eqn:equal_opt:c}
\gamma_j = 0. \vspace{-0.3cm}
\end{equation}
\end{enumerate}
\end{proposition}
The proof is given in Appendix~\ref{app:prop:equal_opt} in the e-companion. 
Note that, from Lemma~\ref{lemma:critical_pattern}, for any critical PS pairs $\iota,\iota'\in\mathscr{I}(o)$ with $\iota\neq \iota'$, it follows that $i_{\iota}\neq i_{\iota'}$.
If the system is weakly coupled, for any $j\in[J]$, there exist at most two different critical pairs $\iota\in\mathscr{I}(o)$ satisfying $j\in\mathscr{J}_{i_{\iota}}$. Also, in a weakly coupled system, for the second case stated in Proposition~\ref{prop:equal_opt},  if there are critical PS pairs $\iota$ and $\iota'$ in $\mathscr{I}(o)$ with critical resource pools $j=j_{\iota}(o)\neq j_{\iota'}(o)$ and $ j_{\iota'}(o)\in\mathscr{J}_{i_{\iota}}$, then $ j_{\iota}(o)\notin \mathscr{J}_{i_{\iota'}}$ because there is at most one resource pool in $\mathscr{J}_{i_{\iota}}$ shared with other patterns.

In Proposition~\ref{prop:equal_opt}, 
the assumption that the system is weakly coupled constrains the way in which resource pools are shared by different requests. 
The case where there is an $o\in\mathscr{O}^*(\bm{0})$ with $\bm{\nu}(o,\bm{0})=\bm{0}$ will occur when the relaxed action constraint~\eqref{eqn:constraint:relax:action} is satisfied with $\alpha^{\bar{\varphi}(o)}_{d(\ell)}(n)>0$ for the only $n\in\mathscr{N}_{d(\ell)}$ and for all $\ell\in[L]$. 
To see this, note that
the construction of the policy $\bar{\varphi}(o)$ guarantees that the resulting multipliers $\bm{\nu}(o,\bm{0})$ will be non-negative, and so
it follows from~\eqref{eqn:a_opt:dummy} that $\alpha^{\bar{\varphi}(o)}_{d(\ell)}(n) > 0$ only if $\nu_{\ell}(o,\bm{0}) = 0$. That is, having $\nu_{\ell}(o,\bm{0}) = 0$ is associated with there being a positive probability that the dummy pattern $d(\ell)$ is selected in the relaxed system.
Furthermore, if there is a PS pair $\iota$ (for a non-dummy pattern $i_{\iota}\in\mathscr{P}_{\ell}$) which satisfies the condition described in \ref{case:equality:action}, that is PS pair $\iota$ causes the relaxed action constraint~\eqref{eqn:constraint:relax:action} to bite, Algorithm~\ref{algo:varphi_gamma} will ensure that $\alpha^{\bar{\varphi}(o)}_{i_{\iota'}}(n_{\iota'}) = 0$ for all PS pairs $\iota'$ 
ranked lower than $\iota$ according to the order $o$. In particular, this will cause  $\alpha^{\bar{\varphi}(o)}_{d(\ell)}(n) = 0$ for the only $n\in\mathscr{N}_{d(\ell)}$.

So if $\alpha^{\bar{\varphi}(o)}_{d(\ell)}(n) > 0$, it is because the relaxed capacity constraints~\eqref{eqn:constraint:relax:resources} bite before the relaxed action constraints~\eqref{eqn:constraint:relax:action}.
If this is true for all $\ell$, then the capacity constraints are biting for every request type, and so we refer to the case where there is an $o\in\mathscr{O}^*(\bm{0})$ with $\bm{\nu}(o,\bm{0})=\bm{0}$ as a \emph{heavy traffic} condition.

\begin{condition}{Heavy Traffic}\label{cond:heavy_traffic}
The system is in heavy traffic if there is a ranking $o\in\mathscr{O}^*(\bm{0})$ such that $\bm{\nu}(o,\bm{0})=\bm{0}$. \vspace{-0.5cm}
\end{condition}

{\bf Remark} 
The property of being weakly coupled and in heavy traffic simplifies the analysis of the complementary slackness condition of the relaxed problem. In particular, the index related to a pattern, described in Equation~\eqref{eqn:index_value}, is affected only by the multipliers of resource pools $j \in[J]$ with $w_{j,i}>0$. Weak coupling helps reduce the number of such multipliers $\gamma_j$, so that the index of a pattern is affected by at most one $\gamma_ j$, which in turn affects other pattern indices. When the system is weakly coupled and in heavy traffic, we can analytically solve the $I$ linear equations~\eqref{eqn:necessary_gamma} and \eqref{eqn:necessary_gamma:zero} and derive the $\phi$ and $\pmb{\gamma}$ that satisfy the complementary slackness condition described in equaitons~\eqref{eqn:relaxaction:slack} and \eqref{eqn:relaxconstraint:slack}.  A detailed discussion is provided in the proof of Proposition~\ref{prop:equal_opt}.

Proposition~\ref{prop:equal_opt} guarantees the decomposability of a system when it is weakly coupled and in heavy traffic.
The property of being weakly coupled and in heavy traffic is stronger than necessary for decomposability, but it is simple to check and is satisfied in a number of common resource allocation problems. We consider examples about how to easily define such a system.

As explained above, the heavy traffic property is usually satisfied when the service capacity is not enough (or just enough) to address its high traffic load.
On the other hand, the weak coupling specifies the structure of the weight matrix $\mathcal{W}$. 
For instance, if each pattern involves only one resource pool (that is, for all $i\in[I]\backslash\{d(\ell)|\ell\in[L]\}$, $|\mathscr{J}_i|=1$), then the system is weakly coupled as each resource pool is still potentially shared by requests of different types.

Within the above framework, we can model skill-based resource pooling in call centers (see \cite{wallace2004resource,cezik2008staffing}) as a weakly coupled resource allocation problem; and when its traffic load is also heavy, the system is decomposable.
In each call center, agents are trained for several skills, such as two or three languages, and are able to handle some but not all of the incoming calls. We classify these agents into multiple call centers according to their trained skills; that is, all agents in the same call center have the same skills and are able to serve the same types of calls. In this context, a call corresponds to a request, an agent corresponds to a resource unit, a call center is a resource pool and a call type is a request type.

Since each call is served by an agent with appropriate skills, each pattern consists of only one call center (resource pool) and selecting a pattern means selecting an agent (a resource unit) from the corresponding call center: this problem is weakly coupled. Note that agents of each call center are potentially serving different types of calls simultaneously, and the capacity constraints \eqref{eqn:constraint:resources} still restrict the system because of the limited number of agents in each call center.

In particular, the trained skills are used to establish the availability of an agent to serve a call, and do not relate to any concept defined in the resource allocation problem. When an agent of call center $j\in[J]$  is able to serve calls of type $\ell\in[L]$, regardless of the skills needed for this service, there is a pattern $i$ in $\mathscr{P}_{\ell}$ with $w_{j,i}=1$ and $w_{j',i}=0$ for all $j'\in[J]\backslash\{j\}$. 
For instance, a call center has agents who can speak English and Chinese, and there are two types of calls: one requires English or French-speaking agents and the other Chinese or Japanese speakers. A call of either type can be served by an agent of this call center, and many calls of both types can be served by this call center simultaneously.

Other problems with similar features, such as health-care task scheduling for agents with different  qualifications (see \cite{lieder2015task}) and home health-care scheduling (see \cite{fikar2017home}), can also be modeled as weakly coupled systems. And, of course, when the systems are also in heavy traffic, they are decomposable.

A virtual machine (VM) replacement problem can be modeled as a resource allocation problem (see \cite{stolyar2013infinite,stolyar2017large}). VM replacement is about consolidating multiple VMs onto a set of physical machines/servers, where each physical server can usually accommodate more than one VM simultaneously. To consolidate a VM, certain numbers of physical units, such as  CPU cycles, memory, disk, or I/O ports, located on a server will be occupied by this VM until it is completed. The VMs and servers are potentially different, and, because of different server profiles or user preference, a server is not necessarily able to accommodate every VM. 
Consider a simple version, for which the capacity of a server is determined by the total amount of only one type of physical unit: this server has a plentiful supply of the others or is not aware of other physical units.
In this case, regarding a VM as a request, a server as a resource pool and a physical unit of the shortage type as a resource unit, we obtain a resource allocation problem that is weakly coupled. 
Similar problems in computer networks, such as the virtual node embedding (see \cite{esposito2016distributed}), server provisioning in distributed cloud environments (see \cite{wei2017data}), and wireless resource scheduling (see \cite{chen2017wireless}), can potentially be modeled as weakly coupled resource allocation problems. And as before, when the weakly coupled systems are in heavy traffic, the decomposability property holds.

As in \cite{stolyar2013infinite,stolyar2017large}, for general VM replacement problems,  each server capacity is not necessarily constrained by physical units of just one type. As above, we model a VM as a request, a physical unit as a resource unit, and the set of all physical units of the same type located on the same server as a resource pool. In this context, the capacity of each resource pool is determined by the total number of its associated physical units of a given type on a server and the weak coupling property cannot hold in general. 
It follows that, unlike the preceding examples, the system is not necessarily decomposable.
However, as discussed in Section~\ref{subsubsec:derive_ranking}, a decomposable system that is not weakly coupled or in heavy traffic can be found by finding a fixed point $\pmb{\gamma}\in\mathbb{R}_0^J$ of the function $\mathcal{T}^o$ ($o\in\mathscr{O}(\pmb{\gamma},\bm{0})$). Numerical examples of such systems will be provided in Section~\ref{sec:example}.

\section{The Index Policy: Its Implementation in the Non-Asymptotic Regime}\label{sec:index_policy}
In Section~\ref{sec:relaxation}, we considered the relaxed problem with constraints \eqref{eqn:constraint:relax:action}-\eqref{eqn:constraint:dummy}. Here, we return to the original problem with constraints \eqref{eqn:constraint:action} and \eqref{eqn:constraint:resources}.

For each RT $\ell\in[L]$, we refer to a policy $\varphi\in\Phi$ as an \emph{index policy} according to PS-pair ranking $o\in\mathscr{O}$, if it always prioritizes a candidate process in 
a PS pair with a ranking equal or higher than those of all the other candidate processes.
This  policy $\varphi$ is applicable to the original problem while, the policy $\bar{\varphi}(o)$ proposed in Section~\ref{subsubsec:priorities_of_pairs} is not in general.
The method of implementing such a $\varphi$ is not unique; for instance, the computation of the ranking of the PS pairs can vary. Here we propose one possible implementation.

For $t>0$, we maintain a sequence of $I$ ordered PS pairs $(i,N^{\varphi}_i(t))$ ($i\in[I]$) that are associated with the $I$ patterns, according to the given ranking $o$ and the state vector $\bm{N}^{\varphi}(t)$: PS pair $(i,N^{\varphi}_i(t))$ is placed ahead of $(i',N^{\varphi}_{i'}(t))$ if and only if the former precedes the latter in the ranking $o$.
Let $i^o_{\sigma}(\bm{N}^{\varphi}(t))$ ($\sigma\in[I]$) represent the pattern associated with the $\sigma$th PS pair in this ordered sequence.

For a general ranking $o\in\mathscr{O}$, the variables $i^o_{\sigma}(\bm{N}^{\varphi}(t))$ are potentially updated at each state transition. Nonetheless, for the purpose of this paper, we mainly focus on the rankings $o\in\mathscr{O}(\pmb{\gamma},\bm{\nu})$ (for some $\pmb{\gamma}\in\mathbb{R}^J_0$ and $\bm{\nu}\in\mathbb{R}^L$) that follow the descending order of $\Xi_i(\pmb{\gamma},\bm{\nu})$.
In this case, the variables $i^o_{\sigma}(\bm{N}^{\varphi}(t))$ are updated only if a pattern $i\in[I]\backslash\{d(\ell)|\ell\in[L]\}$ transitions into or out of its boundary state $|\mathscr{N}_i|-1$.

Consider the capacity constraints
\begin{equation}\label{eqn:constraint:resources:epsilon:1}
\sum\limits_{i'\in[I]}w_{j,i'}N^{\varphi}_{i'}(t) + \sum\limits_{i'\in[I],i'\neq i}w_{j,i'}a^{\varphi}_{i'}(\bm{N}^{\varphi}(t))+ a^{\varphi}_i(\bm{N}^{\varphi}(t)) \leq \Bigl\lceil C_j\bigl(1-\bar{\epsilon}_{j,\iota(i,N^{\varphi}_i(t))}\bigr)\Bigr\rceil,~\forall j\in[J], i\in[I],
\end{equation}
where $\iota(i,N^{\varphi}_i(t))\in[N]$ represents the order of PS pair $(i,N^{\varphi}_i(t))$ in the ranking $o$ and $\bar{\bm{\epsilon}}\in [0,1]^{J\times [N]}$ is a given matrix of parameters.
Apart from this matrix of parameters, constraints~\eqref{eqn:constraint:resources:epsilon:1}  are the same as constraints~\eqref{eqn:constraint:resources}.
As we shall discuss in Section~\ref{subsec:policies}, we choose the $\bar{\epsilon}_{j,\iota}$ such that $\bar{\epsilon}_{j,\iota}C_j \geq w_{j,i_{\iota}}-1$ and, for any $j\in[J]$ and PS pairs $\iota<\iota'$ with respect to the given ranking $o$, if $w_{j,i_{\iota}},w_{j,i_{\iota'}}>0$, then $\bar{\epsilon}_{j,\iota} < \bar{\epsilon}_{j,\iota'}$.
In this context, if $\bar{\epsilon}_{j,\iota}C_j \in [ w_{j,i_{\iota}}-1,w_{j,i_{\iota}})$ for all $\iota\in[N]$ and $j\in[J]$, then constraints~\eqref{eqn:constraint:resources:epsilon:1} reduce to \eqref{eqn:constraint:resources}; otherwise, they are more stringent than \eqref{eqn:constraint:resources}.
The parameter $\bar{\bm{\epsilon}}$ is used to specify the trajectory of the underlying process $\bm{N}^{\varphi}(t)$ when the system is scaled to the asymptotic regime. 
This specification is required 
for proving the asymptotic optimality of the index policy $\varphi$.

In the interests of notational consistency, we shall use the form \eqref{eqn:constraint:resources:epsilon:1} throughout
but, here, since we do not worry about the asymptotic behavior, we consider the case with $\bar{\epsilon}_{j,\iota}C_j \in [ w_{j,i_{\iota}}-1,w_{j,i_{\iota}})$ for all $\iota\in[N]$ and $j\in[J]$ so that \eqref{eqn:constraint:resources:epsilon:1} reduces to \eqref{eqn:constraint:resources}.
A detailed discussion about the scaling procedure and the role of $\bar{\bm{\epsilon}}$ in the asymptotic case will be provided in Section~\ref{sec:asym}.

Under the index policy $\varphi$, we select $L$ patterns to accept new arrivals of $L$ types according to their orders in sequence $i^o_{\sigma}(\bm{N}^{\varphi}(t))$ ($\sigma\in[I]$).
In particular, at a decision making epoch $t>0$, we initialize $a^{\varphi}_i(\bm{N}^{\varphi}(t))=0$ for all $i\in[I]$ and a set of \emph{available patterns} to be $[I]$. 
If, for $i=i^o_1(\bm{N}^{\varphi}(t))$, constraints~\eqref{eqn:constraint:resources:epsilon:1}
will not be violated by setting $a^{\varphi}_i(\bm{N}^{\varphi}(t))=1$, then set $a^{\varphi}_i(\bm{N}^{\varphi}(t))=1$ and remove all patterns associated with request type $\ell(i)$ from the set of available patterns.

The other $L-1$ patterns are selected similarly and iteratively.
That is, we look for the smallest $\sigma\in\{2,3,\ldots,I\}$ such that 
\begin{itemize}
\item $i^o_{\sigma}(\bm{N}^{\varphi}(t))$ remains in the set of available patterns; and
\item the capacity constraints~\eqref{eqn:constraint:resources:epsilon:1} will not be violated by setting $a^{\varphi}_i(\bm{N}^{\varphi}(t))=1$ where $i=i^o_{\sigma}(\bm{N}^{\varphi}(t))$.
\end{itemize}
If there is such a $\sigma$, set $a^{\varphi}_i(\bm{N}^{\varphi}(t))=1$ for $i=i^o_{\sigma}(\bm{N}^{\varphi}(t))$, remove all patterns associated with request type $\ell(i)$ from the set of available patterns and continue selecting the remaining $L-2$ patterns in the same manner.
When all of the $L$ patterns have been selected we can stop. 
Detailed steps are provided in Algorithm~\ref{algo:index_policy}, 
which has a computational complexity that is linear in $I$.

\IncMargin{1em}
\begin{algorithm}
\small 
\linespread{0.4}\selectfont

\SetKwFunction{FIndexPolicy}{IndexPolicy}
\SetKwProg{Fn}{Function}{:}{\KwRet $\bm{a}^{\varphi}(\bm{n})$}
\SetKwInOut{Input}{Input}\SetKwInOut{Output}{Output}
\SetAlgoLined
\DontPrintSemicolon
\Input{a ranking of PS pairs $o\in\mathscr{O}$ and a given state $\bm{n}\in\mathscr{N}$.}
\Output{the action variables $\bm{a}^{\varphi}(\bm{n})$ under the index policy $\varphi\in\Phi$ with respect to ranking $o$ when the system is in state $\bm{n}$.}
\Fn{\FIndexPolicy{$o,\bm{n}$}}{
	$\bm{a}^{\varphi}(\bm{n}) \gets \bm{0}$ \tcc*{Initializing the action variables}
	$\mathscr{P}\gets [I]$ \tcc*{Initializing the set of available patterns}
	$\sigma \gets 1$ \tcc*{Starting with the pattern with the highest priority}
    \While {$\mathscr{P}\neq \emptyset $}{	
	 	    $i\gets i^o_{\sigma}(\bm{n})$\;
			\If {$i\in\mathscr{P}$ {\bf and} Constraints~\eqref{eqn:constraint:resources:epsilon:1} are not violated by setting $a^{\varphi}_i(\bm{n}) =1$ and $\bm{N}^{\varphi}(t)=\bm{n}$}{
			    $a^{\varphi}_i(\bm{n})\gets 1$\;
				Remove all patterns $i'\in\mathscr{P}$ with $\ell(i')=\ell(i)$ from $\mathscr{P}$\;
			}
	    	$\sigma\gets \sigma +1$\;
	}
}
\caption{Implementing the index policy $\varphi$ with respect to ranking $o$.}\label{algo:index_policy}
\end{algorithm}
 \DecMargin{1em}

The performance of $\varphi$ is mainly determined by the given order $o\in\mathscr{O}$. 
Based on later discussion of the asymptotic regime, if the policy $\bar{\varphi}(o)$ is optimal for the relaxed problem in the asymptotic regime, then $\varphi$ is asymptotically optimal for the original problem.
Even without the proved asymptotic optimality, the ranking $o$ should ensure good performance of $\varphi$ as it is always rational to prioritize patterns according to their potential profits.
As long as there are reasonably good $\pmb{\gamma}$ and $\bm{\nu}$ leading to a $o\in\mathscr{O}(\pmb{\gamma},\bm{\nu})$, which correctly reflects the potential profits of patterns, the performance degradation of $\bar{\varphi}(o)$ is likely to be limited for the relaxed problem and close to the optimal solution of the original problem; and the index policy $\varphi$ derived from $o$ is a promising choice for managing resources.

The selection of $\pmb{\gamma}$, $\bm{\nu}$  and $o\in\mathscr{O}(\pmb{\gamma},\bm{\nu})$ is discussed in Section~\ref{sec:relaxation}. The key point is to guarantee good performance of $\bar{\varphi}(o)$: the policy that is guaranteed to be optimal for the relaxed problem when the system is decomposable.

\section{Stochastic Optimization in a Scaled System} \label{sec:asym}

In this section, we establish asymptotic optimality of $\varphi$.

\vspace{-0.3cm}
\subsection{Scaling Parameter}\label{subsec:asym_regime}
\vspace{-0.3cm}
With a parameter $h\in\mathbb{N}_+$, let $\bm{C}\coloneqq h\bm{C}^{0}$, $\bm{C}^{0}\in\mathbb{N}_{+}^{J}$, and the arrival rates scale as $\bm{\lambda}\coloneqq h\bm{\lambda}^{0}$, $\bm{\lambda}^{0}\in\mathbb{R}_{+}^{L}$.
We refer to the parameter $h$ as the \emph{scaling parameter}, and the \emph{asymptotic regime} as the limiting case with $h\rightarrow +\infty$.

We split the process associated with pattern $i$ into $h$ identical \emph{sub-processes}
$(i,k)$, $k\in[h]$, and divide $N^{\phi}_{i}(t)$, the number of instantiations for pattern $i$ under policy $\phi$ at time $t$, into $h$ pieces.
The number of instantiations of the $k$th piece is $N^{\phi}_{i,k}(t)$, so that $N^{\phi}_{i}(t) = \sum_{k=1}^{h}N^{\phi}_{i,k}(t)$.
We refer to $N^{\phi}_{i,k}(t)$ as the number of instantiations for \emph{sub-pattern} $(i,k)$.
The counting process given by $N^{\phi}_{i,k}(t)$ ($k\in [h], i\in [I]$) has state space 
$\mathscr{N}_{i}^{0} \coloneqq \{0,1,\ldots,\min_{j\in \mathscr{J}_i}\lceil C_{j}^{0}/w_{j,i}\rceil\}$.
For any dummy  pattern $d(\ell)$, we take $\mathscr{N}^0_{d(\ell)}= \mathscr{N}_{d(\ell)} = \{0\}$.

The objective and constraints defined by \eqref{eqn:objective}, \eqref{eqn:constraint:action} and \eqref{eqn:constraint:resources} still apply to the sums of variables $ \sum_{k=1}^{h}N^{\phi}_{i,k}(t)\coloneqq  N^{\phi}_i(t)$, $i\in[I]$.
We say the process associated with pattern $i$ is \emph{replaced} by the $h$ sub-processes associated with sub-patterns $(i,k)$, $k\in[h]$.
Each sub-pattern $(i,k)$ earns reward $r_{\ell(i)}$ per each served request and the cost rate that a request accommodated by this sub-pattern imposes on resource pool $j\in[J]$ is $\varepsilon_jw_{j,i}$; that is,
the sub-process $(i,k)$ maintains the same reward and cost rates in states $n\in\mathscr{N}_i^0$ as process $i$.
Let $\bm{N}_{h}^{\phi}(t) = (N^{\phi}_{i,k}(t):\ i\in[I], k\in[h])$ be the state variable after this replacement,
and $a^{\phi}_{i,k}(\bm{N}^{\phi}_h(t))\in\{0,1\}$ ($i\in[I],k\in[h]$) be the action variables  with respect to the process $\bm{N}^{\phi}_h(t)$.
To clarify, we rewrite the objective described in \eqref{eqn:objective} as
\begin{equation}\label{eqn:objective:h}
\max\limits_{\phi}  \frac{1}{h}\sum\limits_{i\in[I]} \sum\limits_{k\in[h]}\sum_{n_i\in\mathscr{N}_i} \pi_{i,k}^{\phi,h}(n_i) \Bigl(r_{\ell(i)} \mu_i  - \sum_{j\in\mathscr{J}} w_{j,i} \varepsilon_j\Bigr)n_i,
\end{equation}
where $\pi^{\phi,h}_{i,k}(n_i)$ represents the stationary probability that the state of sub-process $(i,k)$ is $n_i$ under policy $\phi$ with scaling parameter $h$.
We divide the total revenue earned by all sub-patterns by $h\in\mathbb{N}_+$ so that the objective function is always bounded when $h\rightarrow +\infty$. 
The policy $\phi$ in \eqref{eqn:objective:h}  is determined by the action variables $a^{\phi}_{i,k}(\bm{N}^{\phi}_h(t))$ ($i\in[I],k\in[h]$)
subject to \vspace{-0.3cm}
\begin{equation}\label{eqn:constraint:action:h}
\sum\limits_{i\in\mathscr{P}_{\ell}}\sum\limits_{k\in[h]}a^{\phi}_{i,k}(\bm{N}^{\phi}_h(t))=1,~\forall \ell\in[L],~\forall t \geq 0, \vspace{-0.3cm}
\end{equation}
and  \vspace{-0.3cm}
\begin{equation}\label{eqn:constraint:resources:h}
\sum\limits_{i\in[I]}\frac{w_{j,i}}{h}\sum\limits_{k\in[h]}\Bigl(N^{\phi}_{i,k}(t)+a^{\phi}_{i,k}(\bm{N}^{\phi}_h(t))\Bigr)\leq C_j^0,~\forall j\in[J],~\forall  t \geq 0.
\end{equation}
The constraints in \eqref{eqn:constraint:resources:h} are obtained by substituting $C_j^0h$ for $C_j$ in the constraints stated in \eqref{eqn:constraint:resources}, and thus \eqref{eqn:constraint:resources:h} is equivalent to \eqref{eqn:constraint:resources}.
Also, to guarantee that the maximal value of each $N^{\phi}_{i,k}(t)$ ($k\in[h]$), $\min_{j\in\mathscr{J}_i}\lceil C_j^0/w_{j,i}\rceil$, is not exceeded, define, for $k\in[h]$ and $i\in[I]\backslash\{d(\ell):\ \ell\in[L]\}$,
\vspace{-0.4cm}
\begin{equation}\label{eqn:constraint:resources:zero:h}
a^{\phi}_{i,k}(\bm{N}^{\phi}_h(t))=0,~\text{if }  N^{\phi}_{i,k}(t)=|\mathscr{N}_i^0|-1,
\vspace{-0.4cm}
\end{equation}
which corresponds to the redundant constraints described in \eqref{eqn:constraint:resources:zero}.

{\bf Remark}
As in Section~\ref{sec:introduction},  we activate exactly one sub-process $(i,k)$ ($i\in\mathscr{P}_{\ell}$, $k\in[h]$) for RT $\ell\in[L]$ regardless of the scaling parameter $h\in\mathbb{N}_+$.
The birth and death rates of this active sub-process are $h\lambda_{\ell}^0$ and $N^{\phi}_{i,k}(t)\mu_{i}^0$, respectively, so that if $h\lambda^0_{\ell} \gg (|\mathscr{N}^0_i|-1)\mu_{i}^0$, the number of instantiations of pattern $i$ will increase rapidly until it is restricted by the capacity constraints.

A model with a single active sub-process at any time has different stochastic properties compared to the case where the number of active sub-processes is proportional to $h$ (which was discussed in \cite{weber1990index}). 
To illustrate the difference, we present an example in Appendix~\ref{app:example3} of the e-companion.

The optimization problem consisting of the $hI$ sub-processes associated with $hI$ sub-patterns, coupled through constraints~\eqref{eqn:constraint:action:h}-\eqref{eqn:constraint:resources:zero:h} can be analyzed and relaxed along the same lines as in Section~\ref{sec:relaxation}. 
Let $\alpha^{\phi}_{i,k}(n)\coloneqq \lim_{t\rightarrow +\infty}\mathbb{E}\{a^{\phi}_{i,k}(\bm{N}^{\phi}_h(t))| N^{\phi}_{i,k}(t)=n\}$ ($n\in\mathscr{N}^0_{i}$, $i\in[I]$, $k\in[h]$) represent the action variables of the $hI$ sub-problems for the relaxed problem scaled by $h$. 


All the sub-processes corresponding to a given pattern $i\in[I]$ in the same state $n\in\mathscr{N}_i^0$ are equivalent.
The controller then is concerned only  with the total number of active sub-processes of a given pattern in a given state. 
Define the random variable $Z^{\phi,h}_{\iota}(t)$ to be the proportion of sub-processes in PS pair $\iota$ at time $t$ under policy $\phi$ where $h$ is the scaling parameter; that is, \vspace{-0.3cm}
\begin{equation}\label{eqn:z_process}
Z^{\phi,h}_{\iota}(t)=\frac{1}{hI}\Bigl|\bigl\{(i,k)\in[I]\times[h]\ \left|\ N^{\phi}_{i,k}(t)=n_{\iota},\ i_{\iota}=i\right.\bigr\}\Bigr|. \vspace{-0.3cm}
\end{equation} 

Let $\bm{Z}^{\phi,h}(t)=(Z^{\phi,h}_{\iota}(t):\ \iota\in[N])$ and $\mathscr{Z}$ be the probability simplex $\{\bm{z}\in[0,1]^{N}\ |\ \sum_{\iota\in[N]}z_{\iota}=1\}$.
In this model, the process $\bm{Z}^{\phi,h}(t)$ is analogous to the counting process $\bm{N}^{\phi}_h(t)$ in the original process.
When the process $\bm{Z}^{\phi,h}(t)$ takes value $\bm{z}\in\mathscr{Z}$, it can transition only  to a state of the form $\bm{z}+\bm{e}_{\iota,\iota'}\in\mathscr{Z}$ with $i_{\iota}=i_{\iota'}$, where $\bm{e}_{\iota,\iota'}\in \mathbb{R}^{N}$ is a vector with $\iota$th element $+1/hI$, $\iota'$th element $-1/hI$ and all the other elements set to zero. 
For our birth-and-death process, a transition will happen only with $n_{\iota'} = n_{\iota}\pm 1$ corresponding to the arrival and departure of a request, respectively.
For any given $h\in\mathbb{N}_+$, 
the state space of the process $\bm{Z}^{\phi,h}(t)$ is a subset of $\mathscr{Z}$ and thus the system is always stable.
We refer to the system with $h\rightarrow +\infty$ as the \emph{asymptotic regime}.

Note that any resource allocation problem in the non-asymptotic regime coincides with a scaled problem described in \eqref{eqn:objective:h}-\eqref{eqn:constraint:resources:zero:h} with given $h<+\infty$. The scaling parameter $h$ is introduced to rigorously specify the trajectory of the entire system going from a non-asymptotic regime to an asymptotic regime.

\vspace{-0.5cm}

\subsection{Index Policies in a Scaled System}\label{subsec:policies}
\vspace{-0.3cm}

In Section~\ref{sec:index_policy}, for a ranking $o\in\mathscr{O}$, we proposed an index policy $\varphi\in\Phi$ for the resource allocation problem in the non-asymptotic regime; this coincides with the problem described in \eqref{eqn:objective:h}-\eqref{eqn:constraint:resources:zero:h} with given $h<+\infty$. 
For clarity, we translate the description of $\varphi$ to a policy used for a scaled system with the notation described in Section~\ref{subsec:asym_regime}.

For a ranking $o\in\mathscr{O}$, the index policy $\varphi$ \emph{activates} a sub-process in the first PS pair $\iota\in[N]$ in ranking $o$ with $Z^{\varphi,h}_{\iota}(t)>0$ and the action and capacity constraints holding; that is, $\varphi$ selects a sub-process $(i_{\iota},k)$ ($k\in[h]$) satisfying $N^{\varphi}_{i_{\iota},k}(t)=n_{\iota}$ and sets $a^{\varphi}_{i_{\iota},k }(\bm{N}^{\varphi}_h(t))$ to $1$. The condition $Z^{\varphi,h}_{\iota}(t)>0$ is required because there has to be some sub-processes in PS pair $\iota$ for us  to be able to activate.
Once a sub-process in PS pair $\iota$ is selected for activation, the action constraint~\eqref{eqn:constraint:action:h} for RT $\ell_{\iota}$ achieves equality: there is exactly one active sub-process for a specific RT $\ell\in[L]$.
Resource units in associated resource pools are reserved for this activated sub-process in PS pair $\iota$. 
In this way, $L$ sub-processes in $L$ different PS pairs will be activated iteratively, according to the ranking $o$, for the $L$ RTs.

Under the index policy $\varphi$,  the transition matrix of process $\bm{Z}^{\varphi,h}(t)$ is determined by the value of $\sum_{k\in[h],N^{\varphi}_{i_{\iota},k}(t)=n_{\iota}}a^{\varphi}_{i_{\iota},k }(\bm{N}^{\varphi}_h(t))$ for each PS pair $\iota\in[N]$, which is either $0$ or $1$  and is dependent on $\bm{N}^{\varphi}_h(t)$ through only $\bm{Z}^{\varphi,h}(t)$.
Define $\upsilon^{\varphi,h}_{\iota}(\bm{z})$, $\iota\in[N]$, $\bm{z}\in \mathscr{Z}$, to be the ratio of the number of active sub-processes in PS pair $\iota$, for which the corresponding sub-patterns are prepared to accept a request,  to the total number of sub-processes in this PS pair under $\varphi$, when the proportions of sub-processes in all PS pairs are currently specified by $\bm{z}$. 
That is, at time $t$, for $\iota\in[N]$,
\begin{equation}
\upsilon^{\varphi,h}_{\iota}\bigl(\bm{Z}^{\varphi,h}(t)\bigr)= \frac{\sum_{k\in[h],N^{\phi}_{i_{\iota},k}(t)=n_{\iota}}a^{\phi}_{i_{\iota},k}(\bm{N}^{\varphi}_h(t))}{IhZ^{\varphi,h}_{\iota}(t)},
\end{equation}
where we recall that 
the numerator on the right hand side relies on $\bm{N}^{\varphi}_h(t)$ through $\bm{Z}^{\varphi,h}(t)\in\mathscr{Z}$.
Note that, for arbitrarily large $h$, the value of $\upsilon^{\varphi,h}_{\iota}(\bm{z})Ihz_{\iota}$ ($\bm{z}\in\mathscr{Z}$), representing the number of active sub-processes in PS pair $\iota$, is never more than $1$ because the policy $\varphi$ must satisfy the action constraints \eqref{eqn:constraint:action:h}.
Let $\bm{\upsilon}^{\varphi,h}(\bm{z})=(\upsilon^{\varphi,h}_{\iota}(\bm{z}):\ \iota\in[N])$.
Although different tie-breaking rules lead to the same process $\bm{Z}^{\varphi,h}(t)$, 
we shall stipulate that, when there is more than one sub-process $(i,k)$ ($i\in [I], k\in[h]$)  in the same PS pair available for activation, we prioritize the one with the smaller value of $k$.
In this context, the variables $\bm{\upsilon}^{\varphi,h}(\bm{z})$, $\bm{z}\in\mathscr{Z}$, provide sufficient information for the index policy $\varphi$ to make decisions on the counting process $\bm{N}^{\varphi}_h(t)$.

Let $\zeta_{\iota}^{\varphi,h}(\bm{z})$ represent the maximal proportion of sub-processes in PS pair $\iota$ that can be active if we consider only the capacity constraints defined by \eqref{eqn:constraint:resources:h} (neglecting the action constraints defined by \eqref{eqn:constraint:action:h}) with proportions of sub-processes in all PS pairs specified by $\bm{z}$ under policy $\varphi$.
We obtain that, for $\iota\in[N]\backslash\{d(\ell):\ \ell\in[L]\}$,
\begin{equation}\label{eqn:tilde_z}
\zeta^{\varphi,h}_{\iota}(\bm{z}) =
\min\Biggl\{z_{\iota},\max\biggl\{0,\min\limits_{j\in\mathscr{J}_i} \frac{1}{w_{j,i}Ih}\Bigl\lceil hC_{j}^{0}(1-\epsilon^h_{j,\iota})-\sum\limits_{\iota'=1}^{N}w_{j,i_{\iota'}}n_{\iota'}z_{\iota'}Ih-\sum\limits_{\iota'\in \mathscr{N}_{\iota}^{+}}w_{j,i_{\iota'}}\upsilon^{\varphi,h}_{\iota'}(\bm{z})z_{\iota'}Ih\Bigr\rceil\biggr\}\Biggr\},
\end{equation}
where $\mathscr{N}_{\iota}^{+}$, $\iota\in[N]$, is the set of PS pairs $\iota'\in[N]$ with $\iota' <\iota$ (with higher priorities than pair $\iota$) with respect to ranking $o$, 
and $\epsilon^h_{j,\iota}\in [0,1]$ corresponds to $\bar{\epsilon}_{j,\iota}$ in \eqref{eqn:constraint:resources:epsilon:1}.

Here, the parameter $\epsilon^h_{j,\iota}$ is defined so that
\begin{equation}\label{eqn:epsilon_condition}
0<\lim\limits_{h\rightarrow +\infty}\epsilon^h_{j,\iota} < \lim\limits_{h\rightarrow +\infty}\epsilon^h_{j,\iota'} \leq 1
\end{equation}
for any $\iota < \iota'$, $w_{j,i_{\iota}}>0$ and $w_{j,i_{\iota'}}>0$.
We need $\epsilon^h_{j,\iota}$ to indicate the priorities of PS pairs in resource pool $j$, because
the last term in \eqref{eqn:tilde_z}, 
$\sum_{\iota'\in \mathscr{N}_{\iota}^{+}}w_{j,i_{\iota'}}\upsilon^{\varphi,h}_{\iota'}(\bm{z})z_{\iota'}Ih$,
is $o(h)$.
In particular, 
in order to follow the strict capacity constraints described in \eqref{eqn:constraint:resources:h}, we need to define the $\epsilon^h_{j,\iota}$ so that $\epsilon^h_{j,\iota}C_{j}^0 h \geq w_{j,i_{\iota}}-1$
for all $j\in[J]$, $\iota\in[N]$ and $h\in\mathbb{N}_+$ and $\lim\limits_{h\rightarrow+\infty}\epsilon^h_{j,\iota}$ exists.
Let $\bm{\epsilon}^h \coloneqq (\epsilon^h_{j,\iota}:\ j\in[J],\iota\in[N])$ and $\bm{\epsilon} \coloneqq \lim\limits_{h\rightarrow +\infty}\bm{\epsilon}^h$. 
Define $\mathscr{E}^h$, $h\in\mathbb{N}_+ \cup \{+\infty\}$, and $\Psi$ as the sets of all such vectors $\bm{\epsilon}^h$ and sequences of such vectors $\psi\coloneqq (\bm{\epsilon}^1,\bm{\epsilon}^2,\ldots,)$, respectively.

Equation \eqref{eqn:epsilon_condition} specifies possible trajectories of $\bm{\epsilon}^h$ as $h\rightarrow +\infty$, and is required for subsequent proofs of asymptotic optimality.
Note that, although the asymptotic regime is a limiting situation,  using an asymptotically-optimal policy is likely to be appropriate for systems with finite but large $h$.

In \eqref{eqn:tilde_z}, the value of $\zeta^{\varphi,h}_{\iota}(\bm{z})$ is constrained by the remaining capacities of relevant resource pools, the proportion of sub-processes currently in  PS pair $\iota$ and the proportions of active sub-processes in PS pairs with higher priorities.
Recall that  $\upsilon^{\varphi,h}_{\iota}(\bm{z})$ represents  the proportion of active sub-processes in PS pair $\iota$, for which the corresponding sub-patterns are prepared to accept a request, when the proportions of sub-processes in all PS pairs are currently specified by $\bm{z}$.
Together with the action constraints described in \eqref{eqn:constraint:action:h}, under an index policy $\varphi$, 
for $z_{\iota}>0$,\vspace{-0.15cm}
\begin{equation}\label{eqn:probability_active}
\upsilon^{\varphi,h}_{\iota}(\bm{z}) = \frac{1}{z_{\iota}hI}
\min\biggl\{\zeta^{\varphi,h}_{\iota}(\bm{z})hI,\  \max\Bigl\{0,1- \sum\limits_{\begin{subarray}\ \iota'\in \mathscr{N}^{+}_{\iota}\\ l_{\iota'}=l_{\iota}\end{subarray}}\zeta^{\varphi,h}_{\iota'}(\bm{z})hI\Bigr\}\biggr\}.
\end{equation}
If $z_\iota = 0$, then there are no sub-processes in PS pair $\iota$ and so $\upsilon^{\varphi,h}_{\iota}(\bm{z})$ can take any value in $[0,1]$ without making a difference to the evolution of the process.
For completeness, define, for $\bm{z}$ with $z_{\iota}=0$ and $\bm{z}_{\iota}^{x} 
\coloneqq (z_{1},z_{2},\ldots,z_{\iota-1},x,z_{\iota+1},\ldots,z_{N})$, 
$\upsilon^{\varphi,h}_{\iota}(\bm{z})=\lim\limits_{x\downarrow 0}\upsilon^{\varphi,h}_{\iota}(\bm{z}^{x}_{\iota})$.
For any given $\bm{z}\in\mathscr{Z}$, $\zeta^{\varphi,h}_{\iota}(\bm{z})$ and $\upsilon^{\varphi,h}_{\iota}(\bm{z})$ can be obtained iteratively using equations \eqref{eqn:tilde_z} and \eqref{eqn:probability_active} from $\iota=1$ to $N$.

{\bf Remark}
Although capacity constraints were not considered in \cite{weber1990index}, the construction of $\bm{\upsilon}^{\varphi,h}(\bm{z})$ ($\bm{z}\in\mathscr{Z}$) follows ideas similar to those used in that paper. 
Recall that, for a given ranking, $o\in\mathscr{O}$,  the policy $\bar{\varphi}(o)$ is generated by Algorithm~\ref{algo:varphi_gamma} and is infeasible for the original problem. 
This gives rise to the interesting property that values of $\upsilon^{\varphi,h}_{\iota}(\bm{z})$ and $\alpha^{\bar{\varphi}(o)}_{i_{\iota},k}(n_{\iota})$ ($k\in[h]$) for all $h\in\mathbb{N}_+\cup\{+\infty\}$, are always independent of those of PS pairs $\iota'$ with $\iota' > \iota$: the PS pairs with lower priorities than $\iota$.
The property is important 
for the proofs of Theorems~\ref{theorem:main_second} and \ref{theorem:main}.

\vspace{-0.3cm}

\subsection{Asymptotic Optimality}\label{sec:asym_opt}
\vspace{-0.3cm}

For given $h\in\mathbb{N}_+$, define the long-run average revenue normalized by $h$ of the resource allocation problem under policy $\phi$ to be $R^{\phi,h}$; that is, \vspace{-0.2cm}
\begin{equation}
R^{\phi,h} \coloneqq \frac{1}{h}\sum\limits_{i\in[I]} \sum\limits_{k\in[h]}\sum_{n_i\in\mathscr{N}_i} \pi_{i,k}^{\phi,h}(n_i) \Bigl(r_{\ell(i)} \mu_i  - \sum_{j\in\mathscr{J}} w_{j,i} \varepsilon_j\Bigr)n_i, \vspace{-0.2cm}
\end{equation}
the objective function described in \eqref{eqn:objective:h}.
\vspace{-0.3cm}
\begin{definition}
We say that the index policy $\varphi$ derived from PS ranking $o$ by iterating \eqref{eqn:tilde_z} and \eqref{eqn:probability_active} is \emph{asymptotically optimal} if 
$
\lim\limits_{\lVert \bm{\epsilon}\rVert \to \bm{0}}\lim\limits_{h\rightarrow+\infty }|R^{\varphi,h}-\max\limits_{\phi\in\Phi}R^{\phi,h}| = 0$.
\end{definition}
Recall that the index policy $\varphi$ described in Section~\ref{subsec:policies} is dependent on the parameter $\bm{\epsilon}^h$ with $\bm{\epsilon} \coloneqq \lim\limits_{h\rightarrow +\infty}\bm{\epsilon}^h$. The $\bm{\epsilon}$ is used to guarantee strict priorities of the sub-processes in the asymptotic regime as discussed in Section~\ref{subsec:policies}. 
The policy $\bar{\varphi}(o)$, proposed in Section~\ref{subsubsec:priorities_of_pairs} for the relaxed problem, is generally not applicable to the original resource allocation problem.
Although the policies $\bar{\varphi}(o)$ and $\varphi$ both rely on the same ranking $o\in\mathscr{O}$, they are different policies. 
\begin{theorem}\label{theorem:main_second}
For given $o\in\mathscr{O}$, 
$\varphi$ derived from $o$ by iterating \eqref{eqn:tilde_z} and \eqref{eqn:probability_active} is asymptotically optimal if and only if 
\vspace{-0.5cm}
\begin{equation}\label{eqn:main_second}
\lim\limits_{h\rightarrow +\infty}|R^{\bar{\varphi}(o),h} -\max\limits_{\phi\in\Phi}R^{\phi,h}|=0.
\end{equation}
\end{theorem}
The proof is given  in Appendix~\ref{app:theorem:main_second}.
Theorem~\ref{theorem:main_second} indicates that asymptotic optimality of $\varphi$ is equivalent to the convergence between $R^{\bar{\varphi}(o),h}$ and the maximized long-run average revenue of the original problem as $h\rightarrow +\infty$.
It is proved by showing the existence of $\lim_{h\rightarrow +\infty} \lim_{t \rightarrow +\infty}\mathbb{E}[\bm{Z}^{\bar{\varphi}(o),h}(t)]$ and
a global attractor of the process $\bm{Z}^{\varphi,h}(t)$ as $t,h\rightarrow +\infty$ and $\lVert\bm{\epsilon}\rVert\rightarrow 0$, and specifically that they coincide with each other.
The long-run average revenue $R^{\varphi,h}$ then coincides with $R^{\bar{\varphi}(o),h}$ as $h\rightarrow +\infty$ and $\lVert \bm{\epsilon}\rVert \rightarrow 0$.

A similar condition relevant to the global attractor was required in \cite{weber1990index} for asymptotic optimality of the Whittle index policy in a general RMABP. It does not necessarily hold. However, in our problem, each sub-process is a queueing process with the departure rate increasing in its queue size.
Such a sub-process is a special case of a general bandit process.
We prove in general that the underlying process $\bm{Z}^{\varphi,h}(t)$, regardless of its initial point, will converge to any specified neighborhood of a fixed point almost surely as $t,h\rightarrow +\infty$ and $\lVert \bm{\epsilon}\rVert \rightarrow 0$.
Detailed proofs are provided in Appendix~\ref{app:theorem:main_second}.

Theorem~\ref{theorem:main_second}, in itself,  does not provide a verifiable condition.
This is given in our next theorem. If there exists $H\in\mathbb{R}$ such that, for all $h>H$, the system is decomposable, we say the system is decomposable in the asymptotic regime.
\begin{theorem}\label{theorem:main}
If the capacity constraints described in \eqref{eqn:constraint:resources:h} (or equivalently \eqref{eqn:constraint:resources}) are decomposable with decomposable multipliers $\pmb{\gamma}\in\mathbb{R}^J_0$ in the asymptotic regime, then there exist $\bm{\nu}\in\mathbb{R}^L$ and a PS pair ranking $o\in\mathscr{O}(\pmb{\gamma},\bm{\nu})$ such that the index policy $\varphi$ derived from $o$ is asymptotically optimal.
\end{theorem}
The proof, based on Theorem~\ref{theorem:main_second}, is given in Appendix~\ref{app:theorem:main}. 
Recall that we discussed decomposability of multipliers in Section~\ref{sec:relaxation}
and provided examples of provably decomposable systems.

{\bf Remark}
Theorem~\ref{theorem:main} binds asymptotic optimality of $\varphi$ to the decomposability of the relaxed problem. 
For a decomposable system (see Definition~\ref{define:decomposable}), there always exists a ranking $o\in\mathscr{O}$ such that $\bar{\varphi}(o)$ is optimal for the relaxed problem.
If a system is decomposable in the asymptotic regime 
then \eqref{eqn:main_second} is satisfied.
This follows because, for any $h\in\mathbb{N}_+\cup\{+\infty\}$, $R^{\varphi,h} \leq \max_{\phi\in\Phi}R^{\phi,h} \leq \max_{\phi\in\tilde{\Phi}}R^{\phi,h}$ where $\Phi$ and $\tilde{\Phi}$ are the sets of feasible policies for the original and relaxed problem, and $R^{\varphi,h}$ coincides with $R^{\bar{\varphi}(o),h}$ as $h\rightarrow +\infty$ and $\lVert\bm{\epsilon}\rVert \rightarrow 0$.

Similarly, we say the system is in heavy traffic in the asymptotic regime if there exists $H\in\mathbb{R}$ such that, for all $h>H$, the system is in heavy traffic.
\vspace{-0.3cm} 
\begin{corollary}\label{coro:main}
If the system is in heavy traffic in the asymptotic regime and is weakly coupled,
then there exist decomposable multipliers $\pmb{\gamma}\in\mathbb{R}^J_0$, satisfying \eqref{eqn:equal_opt:a}-\eqref{eqn:equal_opt:c}, and a PS pair  ranking $o\in\mathscr{O}^*(\bm{0})$,
so that  the index policy $\varphi$ derived from $o$ is asymptotically optimal. \vspace{-0.3cm}
\end{corollary}
The proof, invoking Theorem~\ref{theorem:main} and Proposition~\ref{prop:equal_opt}, is given in Appendix~\ref{app:coro:main}. 
In particular, the PS pair ranking $o$, described in Corollary~\ref{coro:main}, exists in closed form: its PS pairs are ranked according to the descending order of $\Xi_{\iota}^*$ with $\bm{\nu}=\bm{0}$.

\section{Numerical Results}\label{sec:example}

We demonstrate by simulation the performance of the index policy $\varphi$, defined in Section~\ref{subsec:policies} (or equivalently, defined in Section~\ref{sec:index_policy} for a given $h<+\infty$), in systems that are not weakly coupled or in heavy traffic in comparison with baseline policies.

In this section, the confidence intervals of all the simulated average revenues at the $95\%$ level based on the Student's t-distribution are maintained within $\pm 3\%$ of the observed mean.
We recall that the capacities $\bm{C}$ and arrival rates $\bm{\lambda}$ are scaled by the scaling parameter $h$.




Along with the fixed point iteration method proposed in Section~\ref{subsubsec:derive_ranking},
we have been able to find systems which are not weakly coupled or in heavy traffic but are decomposable.
Here, we provide two examples, where $L$ and $J$ are sampled uniformly from the sets $\{2,3,4,5\}$ and $\{10,11,\ldots,20\}$, respectively.
Let
$\epsilon_M = \max_{j\in[J],\iota\in[N]}\epsilon_{j,\iota}$.
We refer to an index policy $\varphi$ with specific $\epsilon_M\in [0,1]$ as INDEX($\epsilon_M$).

\begin{figure}[t]
\centering
\subfigure[]{\includegraphics[width=0.45\linewidth]{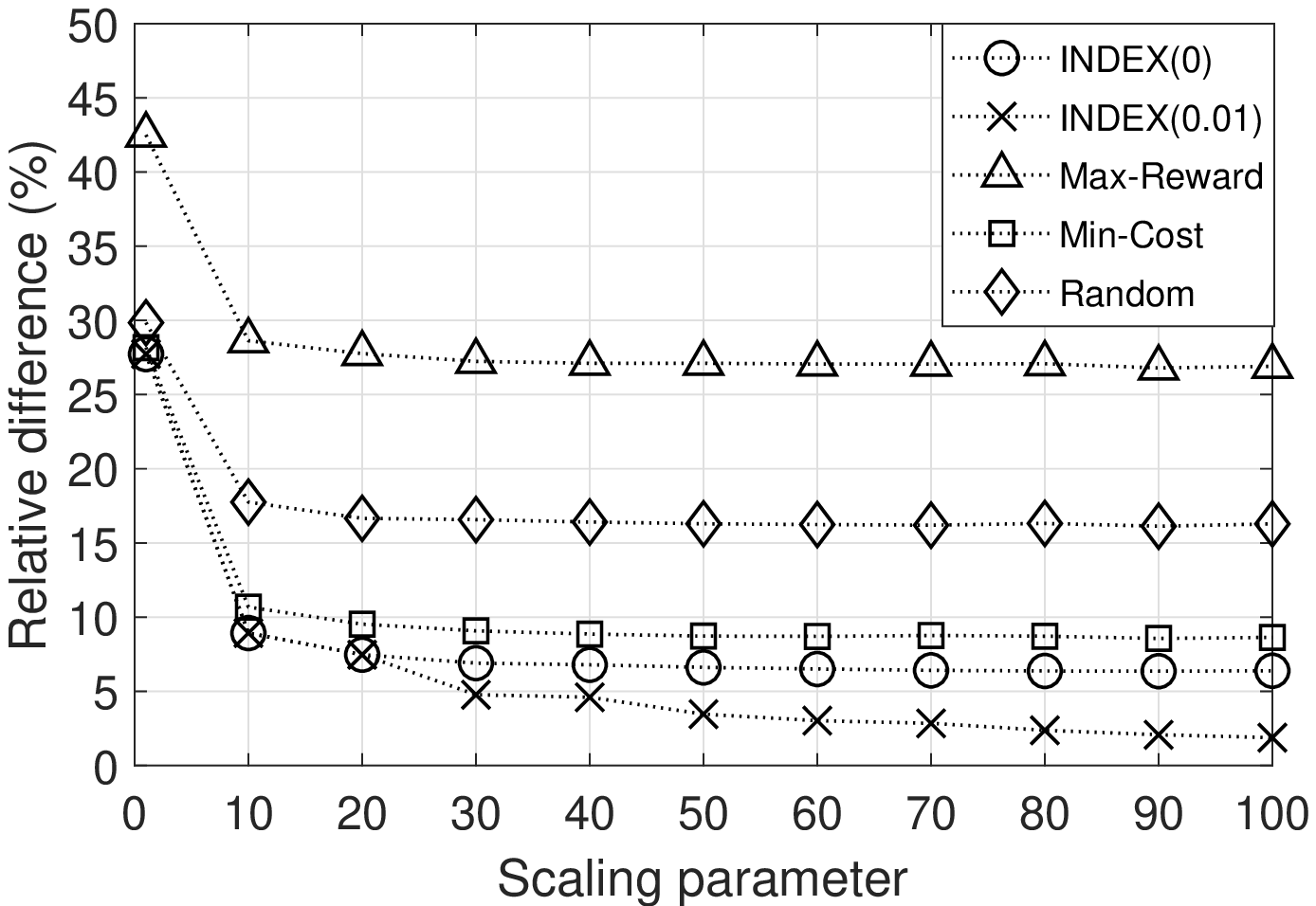}\label{fig1:opt_v1}}
\subfigure[]{\includegraphics[width=0.45\linewidth]{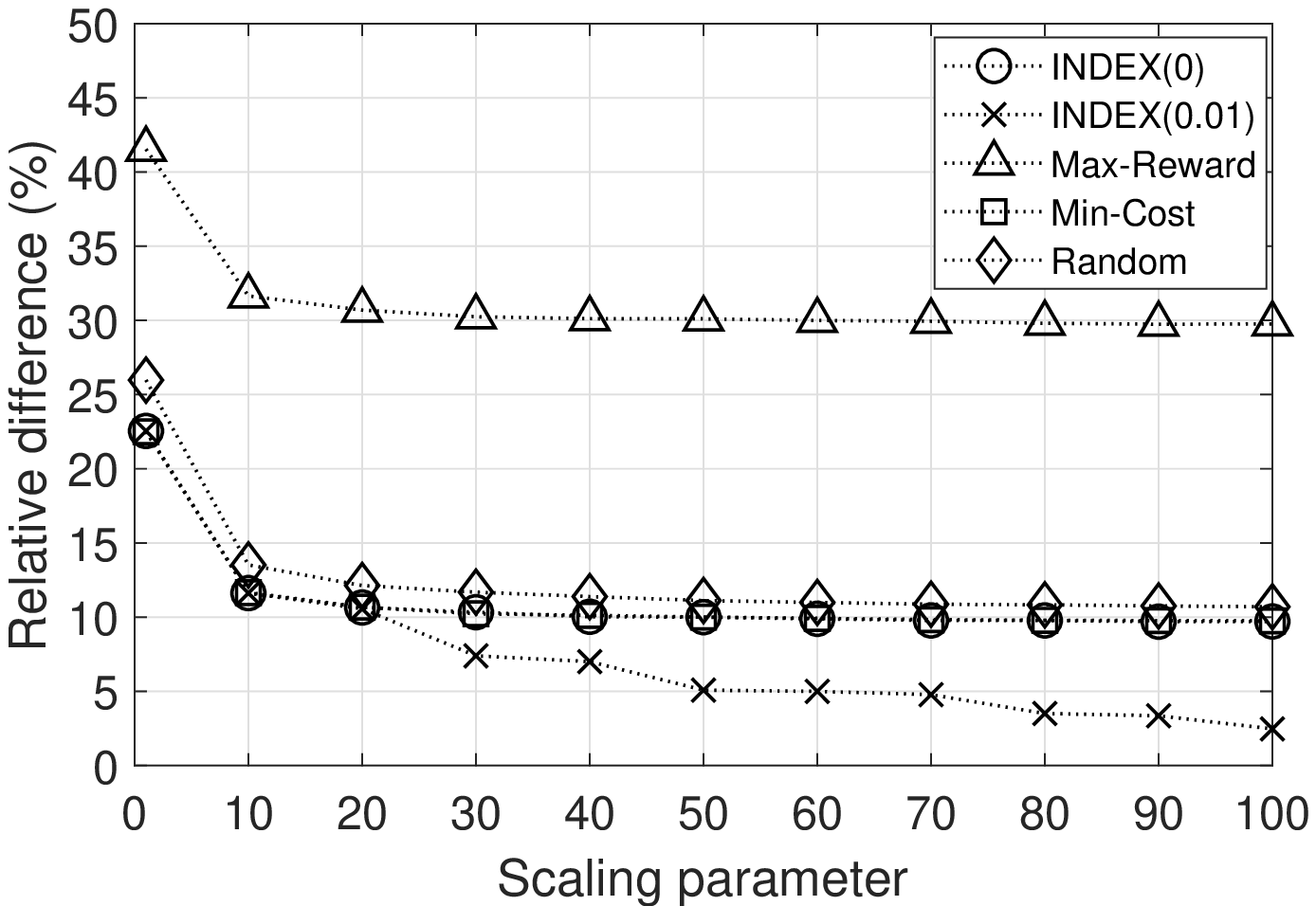}\label{fig1:opt_v3}}
\caption{Relative difference of a specific policy to $R(o_{k^*})$ against the scaling parameter of the system: (a) diverse performance and non-zero decomposable multipliers; (b) similar performance and non-zero decomposable multipliers; and (c) zero decomposable multipliers.}\label{fig:fig1}
\vspace{-0.5cm}
\end{figure}

We consider three baseline policies: two greedy policies that prioritize patterns with maximal reward rates and minimal cost rates, and one policy randomly uniformly selecting an available pattern for each request type. We refer to the three policies as \emph{Max-Reward}, \emph{Min-Cost} and \emph{Random}.
The Max-Reward and Min-Cost policies are in fact index policies with PS pairs ranked according to the descending order of their reward rates and the ascending order of their cost rates, respectively. 
The Random policy was proposed by \cite{stolyar2017large} for a VM replacement problem, aiming to minimize the system blocking probabilities in the case with finite capacities.
It is not a feasible policy of the original problem with capacity constraints \eqref{eqn:constraint:resources} because it does not reserve resource units for a specific pattern that is more profitable than the others.
When there are not enough resource units in a pool to accommodate multiple request types that have chosen their patterns involving this pool, the Random policy will always assign the resource units to the request that arrives first.

In Figure~\ref{fig:fig1}, we compare the performance of INDEX(0), INDEX($0.01$), the baseline policies and $\bar{\varphi}(o_{k^*})$, where $o_{k^*}$ is the ranking of the multipliers $\pmb{\gamma}_{k^*}$ resulting from the fixed point iteration method (described in Section~\ref{subsubsec:derive_ranking}) with parameter $c=0.5$ and initial point $\pmb{\gamma}_0=\bm{0}$. 
The system parameters are listed in Appendix~\ref{app:simulation:opt} and are generated by pseudo-random functions.  
The discovered multipliers $\pmb{\gamma}_{k^*}$ for simulations in Figures~\ref{fig1:opt_v1} and \ref{fig1:opt_v3} are 
$(269.555,0,0,0,0,273.11,0,$ $347.995,0,0,0,8.323\times 10^{-7},9.726\times 10^{-5},0)$  and $\bm{0}$,
respectively, satisfying $\mathcal{T}^{o_{k^*}}(\pmb{\gamma}_{k^*})=\pmb{\gamma}_{k^*}$ in the asymptotic regime.
By Proposition~\ref{prop:converge_gamma}, these $\pmb{\gamma}_{k^*}$ are decomposable multipliers and, by Theorem~\ref{theorem:main}, the index policies derived from the rankings $o_{k^*}$ are asymptotically optimal.
Let $R(o)\coloneqq \lim\limits_{h\rightarrow +\infty}R^{\bar{\varphi}(o),h}$ ($o\in\mathscr{O}$) of which the existence is guaranteed in the proof of Theorem~\ref{theorem:main_second}.
For the decomposable systems with $h<+\infty$ and $\bar{\varphi}(o_{k^*})$ optimal for the relaxed problem in the asymptotic regime, the asymptotic long-run average revenue, $R(o_{k^*})$, is no less than the optimum of the original problem:  $R(o_{k^*})$ is an upper bound of $R^{\phi,h}$ for any $\phi\in\Phi$.

Figure~\ref{fig:fig1} illustrates the relative difference of average revenues,
$\bigl(R(o_{k^*}) - R^{\phi,h}\bigr)/R(o_{k^*})$ for $\phi=\text{INDEX}(0),\text{INDEX}(0.01)$, Max-Reward, Min-Cost and Random,
against the scaling parameter $h$.

In this context, there are two aspects of performance evaluation presented in Figure~\ref{fig:fig1}.
First, we see the performance of the index policies in the non-asymptotic regime by comparing their long-run average revenues with an upper bound on the optimum. 
In particular, Figures~\ref{fig1:opt_v1} and \ref{fig1:opt_v3} show that INDEX($0.01$) significantly outperforms INDEX(0) for large $h$: the small but positive $\bm{\epsilon}$ does affects the performance of $\varphi$.
The performance of INDEX($0.01$) is close to the upper bound of the optimal solution with  relative difference less than $5\%$ for $h$ greater than $50$ in all three examples: its performance degradation against the optimal solution is limited in the non-asymptotic regime.

On the other hand, by comparing to $R(o_{k^*})$, a trend of coincidence between $R^{\text{INDEX}(0.01),h}$ and $R(o_{k^*})$ is observed in Figure~\ref{fig:fig1} as $h$ increases from $1$ to $100$, consistent with the proved asymptotic optimality of $\varphi$.
Recall that the examples presented in Figure~\ref{fig:fig1} are not for systems with weak coupling or heavy traffic but the index policy $\varphi$ is still proved to be asymptotically optimal here.
Also, the performance of $\varphi$ is close to the optimum without requiring extremely large $h$.

\begin{figure}[t]
\vspace{-0.5cm}
\centering
\subfigure[]{\includegraphics[width=0.45\linewidth]{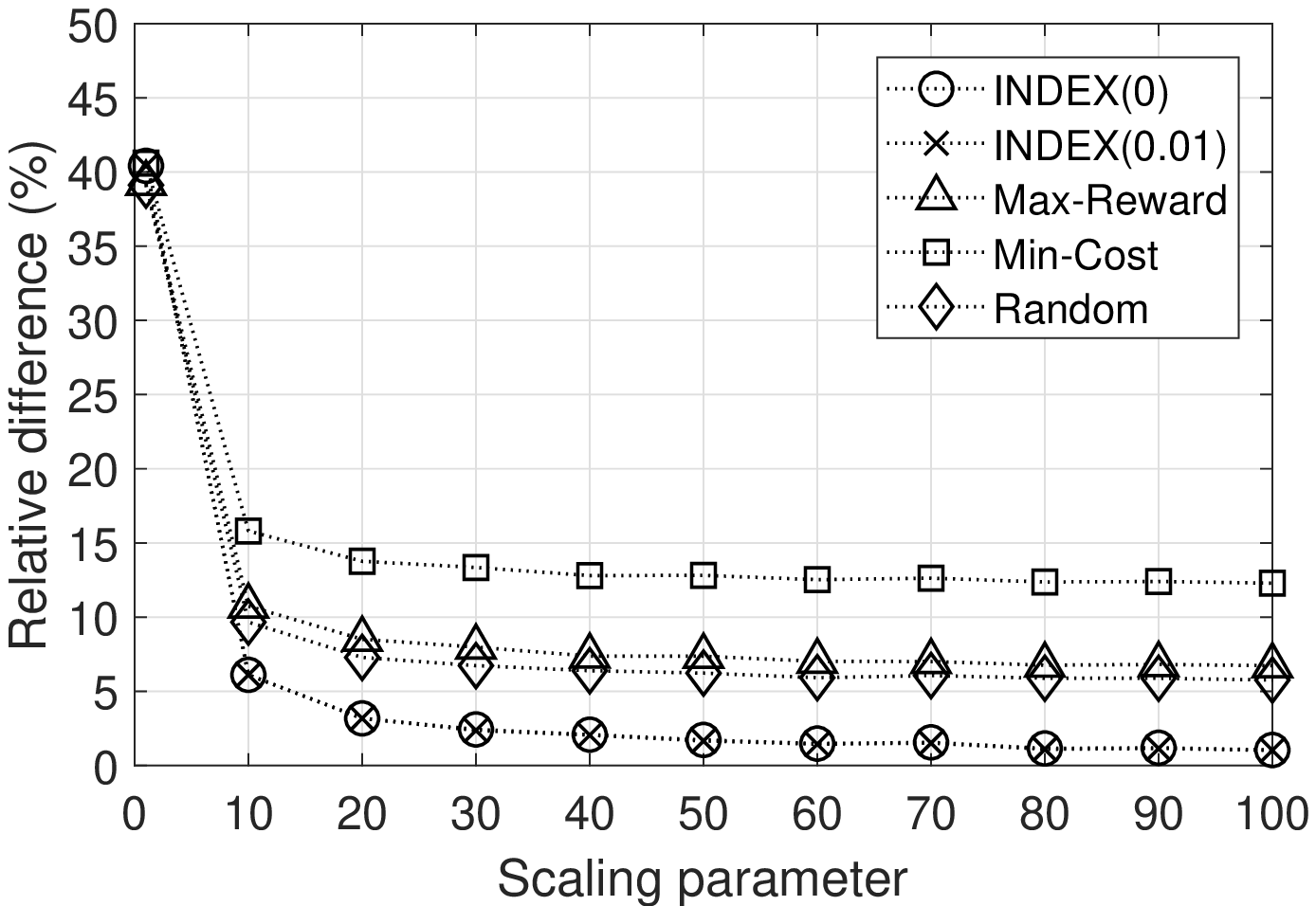}\label{fig2:nonopt_v1}}
\subfigure[]{\includegraphics[width=0.45\linewidth]{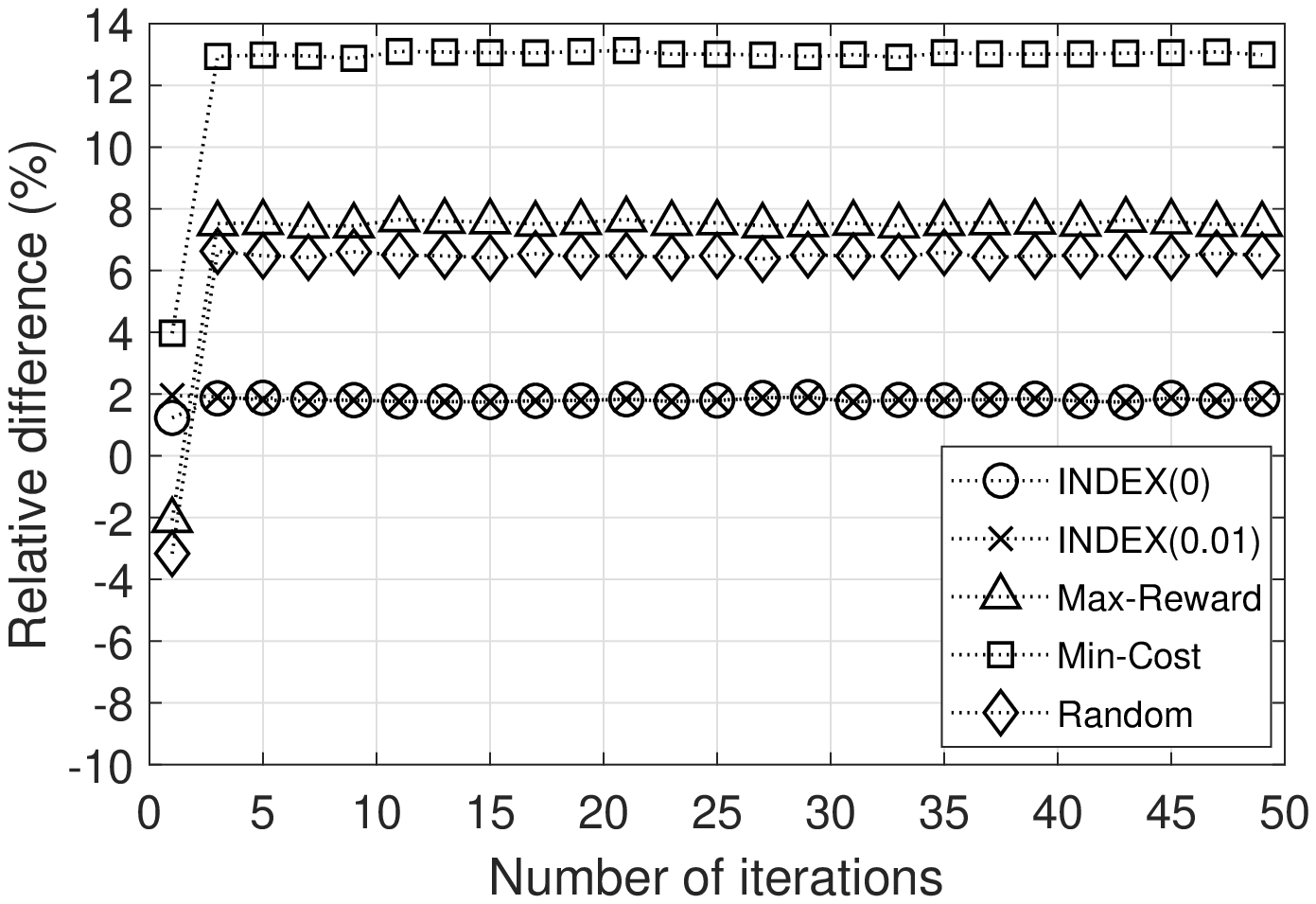}\label{fig2:nonopt_v2}}
\caption{(a) Relative difference of a specific policy to $R(o_{k^*})$ against scaling parameter of the system; (b) Relative difference of a specific policy to $R(o_k)$  against $k$.}\label{fig:fig2}
\vspace{-0.5cm}
\end{figure}

In Figure~\ref{fig:fig2}, we consider another example with multipliers that are not decomposable (that is, $\mathcal{T}^{o_{k^*}}(\pmb{\gamma}_{k^*})\neq \pmb{\gamma}_{k^*}$).
Similar to Figure~\ref{fig:fig1}, in Figure~\ref{fig2:nonopt_v1}, we plot the relative difference of revenue of INDEX(0),  INDEX($0.01$) and the baseline policies to $R(o_{k^*})$ against the scaling parameter; while, in Figure~\ref{fig2:nonopt_v2}, fixing the scaling parameter $h=50$, we illustrate curves of the relative differences, $\bigl(R(o_k)-R^{\phi,h}\bigr)/R(o_k)$ ($\phi=\text{INDEX}(0),\text{INDEX}(0.01), \text{Max-Reward}, \text{Min-Cost}, \text{Random}$), against the number of iterations $k$ for the fixed point iteration method. 
Note that the rankings $o_k$ are potentially different as $k$ varies, so as $R(o_k)$.
In Figure~\ref{fig2:nonopt_v1}, the INDEX($0$) and INDEX($0.01$) are proposed based on the ranking $o_{k^*}$, while, with slightly abused notation, in Figure~\ref{fig2:nonopt_v2},  INDEX($0$) and INDEX($0.01$) represent the index policies $\varphi$, which are derived from the rankings $o_k$ associated with the varying $k$, with $\epsilon_M=0$ and $0.01$, respectively. 
The system parameters for the simulations in Figure~\ref{fig:fig2} are listed in Appendix~\ref{app:simulation:opt}.

Figure~\ref{fig2:nonopt_v1} can be read in a similar way to Figure~\ref{fig:fig1} except that $R(o_{k^*})$ is not a proved upper bound for the average revenue for the original problem.
Here, INDEX($0$) and INDEX($0.01$) perform similarly and numerically converge to $R(o_{k^*})$ as $h$ increases although the system is not necessarily decomposable. 
The convergence is consistent with Theorem~\ref{theorem:main_second} which generally holds without requiring decomposability.
On the other hand, for each finite $h$ (which corresponds to the non-asymptotic regime), 
INDEX($0$) and INDEX($0.01$) significantly outperform all the other baseline policies, although the system is not proved to be decomposable, and their performance advantages are likely to maintain as $h$ continues increasing.

Figure~\ref{fig2:nonopt_v2} illustrates the performance trajectory as the iteration number $k$ (the x-axis) for the fixed point iteration method increases for a system with $h=50$ (in the non-asymptotic regime).
Recall that, for the simulations presented here, the average revenues of INDEX($0$) and INDEX($0.01$) and $R(o_k)$ are varying with $k$ while all of the baseline policies are independent of $k$. 
We observe a shape jump on the curves between $k=1$ and $5$.
This is caused by the initial setting, $\pmb{\gamma}_0=\bm{0}$, which is not a good choice of multipliers.
After several steps of the iteration method, the curves in Figure~\ref{fig2:nonopt_v2} become almost flat; 
that is, the values of $R(o_k)$, $R^{\text{INDEX}(0)}$ and $R^{\text{INDEX}(0.01)}$ become relatively stable for $k=5$ to $50$. 
Also, in Figure~\ref{fig2:nonopt_v2}, after the performance becomes stable, INDEX($0$) and INDEX($0.01$) achieve clearly higher long-run average revenues than those of the baseline policies: given the poor setting at the beginning, the fixed point iteration method can still lead to a reasonably good ranking $o_{k^*}$ and its associated multipliers $\pmb{\gamma}_{k^*}$.

\vspace{-0.3cm}
\section{Conclusions}\label{sec:conclusions}
\vspace{-0.3cm}

We have modeled a resource allocation problem, described by \eqref{eqn:objective}, \eqref{eqn:constraint:action} and \eqref{eqn:constraint:resources}, as a combination of various RMABPs coupled by limited capacities
of the shared resource pools, which are shared, competed for, and reserved by requests.  
This presents us with an optimization problem for a stochastic system, aimed at maximizing the long-run average revenue by dynamically accommodating requests into appropriate resource pools.

Using the ideas of Whittle relaxation \cite{whittle1988restless} and the asymptotic optimality proof of \cite{weber1990index},
we have proved the asymptotic optimality of an index policy (referred to as $\varphi$) if the capacity constraints are decomposable with multipliers $\pmb{\gamma}\in\mathbb{R}_0^J$ (Theorem~\ref{theorem:main}).
The asymptotic optimality is proved based on the existence of a global attractor $\bm{z}\in\mathscr{Z}$ for the underlying process $\bm{Z}^{\varphi,h}(t)$ as $h\rightarrow +\infty$ and $\bm{\epsilon}\rightarrow \bm{0}$. 
We have proved in general that such a global attractor exists, and then proposed a necessary and sufficient condition for asymptotic optimality in Theorem~\ref{theorem:main_second}: the performance of the attractor $\bm{z}$ approaches the optimum of the original problem in the asymptotic regime.
This condition holds if the system is decomposable.

We have proved a sufficient condition, described as the property of being weakly coupled and in heavy traffic, 
for the existence of such decomposable multipliers as well as the asymptotic optimality of policy $\varphi$ (Corollary~\ref{coro:main}).
The property is not necessary, but is easy to verify and covers a significant class of resource allocation problems. We have listed examples of systems with the property satisfied in Section~\ref{subsec:sufficient_condition}.

In a general system, we have proposed a fixed point method to fine tune the multipliers $\pmb{\gamma}\in\mathbb{R}_0^J$ and a ranking $o\in\mathscr{O}(\pmb{\gamma},\bm{0})$. 
We have proved that, if there exists a fixed point $\pmb{\gamma}\in\mathbb{R}_0^J$ of the function $\mathcal{T}^o$ satisfying $o\in\mathscr{O}(\pmb{\gamma},\bm{0})$, then this $\pmb{\gamma}$ is a vector of decomposable multipliers.
We have successfully discovered the decomposable multipliers in some situations without assuming weak coupling or heavy traffic by applying the fixed point method.
Also, in Section~\ref{sec:example}, we have compared the index policy $\varphi$ with different parameter $\bm{\epsilon}$ to baseline policies through simulations for systems that are not weakly coupled or in heavy traffic in the non-asymptotic regime. The index policy achieves clearly higher performance than the baseline policies.
To the best of our knowledge, no existing work provides rigorous asymptotic optimality for a resource allocation problem where dynamic allocation, competition and reservation are permitted.

\vspace{-0.3cm}

\section*{Acknowledgment}
Jing Fu's and Peter Taylor's research is supported by the Australian Research Council (ARC) Centre of Excellence for the Mathematical and Statistical Frontiers (ACEMS) and ARC Laureate Fellowship FL130100039.
\vspace{-0.1cm}



%
%
%



\bibliographystyle{apalike}
\bibliography{references/IEEEabrv,references/bib1}






\appendices

\section{Multiple Restless Multi-Armed Bandit Problems}\label{app:RMABPs}
Consider the special case 
\begin{enumerate}[label=\roman*)]
\item $L=1$; \label{premise:L}
\item $w_{j,i}\in\{0,1\}$ for all $i\in[I]\backslash\{d(\ell):\ \ell\in[L]\}$ and $j\in[J]$; and \label{premise:w}
\item $\sum_{j\in[J]}w_{j,i}w_{j,i'}=0 $ for any $i\neq i'$, $i,i'\in[I]$, \label{premise:j}
(so that the intersection of the sets of resource types used by different patterns is empty). 
\end{enumerate}


Here, the stochastic optimization problem defined by
\eqref{eqn:objective}, \eqref{eqn:constraint:action} and
\eqref{eqn:constraint:resources} becomes
\begin{equation}
\max\limits_{\phi} \lim\limits_{t\rightarrow +\infty}\mathbb{E}\biggl[\sum\limits_{i\in\mathscr{P}_1}\Bigl(r_{\ell(i)}\mu_{i}-\sum\limits_{j\in[J]}\varepsilon_{j}w_{j,i}\Bigr)N^{\phi}_{i}(t)\biggr]
\end{equation}
subject to
\begin{equation}
\sum\limits_{i\in\mathscr{P}_1} a^{\phi}_i(\bm{N}^{\phi}(t)) = 1,~\forall~\bm{N}^{\phi}(t)\in\mathscr{N},
\end{equation}
This is a typical RMABP as defined in \cite{whittle1988restless}, and
the process for each pattern $i\in[I]$ is a
restless Bandit Process (BP).  

Similarly, if Premises~\ref{premise:w} and \ref{premise:j} hold but Premise~\ref{premise:L} does not then the resource allocation problem
can be modeled by $L$ independent RMABPs.  If Premises~\ref{premise:w} and
\ref{premise:j} hold then the capacity constraints
defined in \eqref{eqn:constraint:resources} are no longer necessary.


We can also obtain $L$ independent RMABPs by simply dropping the
capacity constraints in~\eqref{eqn:constraint:resources} without
assuming any of the three premises, but the physical meaning of the
entire problem is changed: the capacity constraints are then
reflected by just  the definition of state space
$\mathscr{N}_i=\{0,1,\ldots, \min_{j\in\mathscr{J}_i}\lceil
C_j/w_{j,i} \rceil \}$, for each pattern $i\in[I]$, rather than be
limited by the sum of RUs simultaneously occupied by different patterns. 
In the general case where capacity constraints are in place, the $L$ RMABPs, each of which corresponds to the $|\mathscr{P}_{\ell}|$ BPs coupled by the action constraint described in~\eqref{eqn:constraint:action}, are further linked by these capacity constraints.

\section{Proof of Proposition~\ref{prop:a_opt}}\label{app:prop:a_opt}
We start with a lemma.
\begin{lemma}\label{lemma:deterministic}
 For all $n=0,1,\ldots,|\mathscr{N}_i|-2$ , there exists a policy $\varphi$ 
determined by $\alpha_i^{\varphi}(n)\in\{0,1\}$ that maximizes the long-run average revenue of the underlying MDP associated with the sub-problem for pattern $i$.
\end{lemma}
\proof{Proof of Lemma~\ref{lemma:deterministic}.}

The underlying MDP associated with the sub-problem for pattern $i$ is a birth-and-death process, and its state variable $n$ represents the number of instantiations generated for pattern $i$.

Defining  state $0$ as the absorbing state,  for all $n\in [|\mathscr{N}_i|-2]$,  Bellman's equation for the value function of this MDP is
\begin{multline}\label{eqn:add_1}
V(n) = \max\limits_{\alpha^{\varphi}_i(n) \in[0,1]} \biggl\{\frac{1}{\alpha^{\varphi}_i(n)\lambda_{\ell(i)}+n\mu_i}\Bigl(r_{\ell(i)}\mu_{i}n - \sum\limits_{j\in\mathscr{J}_i}\varepsilon_j w_{j,i}n-g\Bigr)\\
+ \frac{\alpha^{\varphi}_i(n)\lambda_{\ell(i)}}{\alpha^{\varphi}_i(n)\lambda_{\ell(i)}+n\mu_i}V(n+1) + \frac{n\mu_i}{\alpha^{\varphi}_i(n)\lambda_{\ell(i)}+n\mu_i}V(n-1)\biggr\},
\end{multline} 
where $g\in\mathscr{R}$ is a given parameter that is equal to the maximized long-run average revenue of this MDP and $V(0)=0$.
The parameter $g$ acts as an attached cost, such that the action variable $\alpha^{\varphi}_i(n)$ under an optimal solution $\varphi$ that maximizes the long-run average revenue of this MDP will also maximize the right hand side of \eqref{eqn:add_1} (see \cite{ross1992applied}). 
The value function $V(\cdot)$, which is a solution of the Bellman equation \eqref{eqn:add_1}, can be computed through the \emph{value iteration} technique.
In this context, the expression on the right hand side of \eqref{eqn:add_1} is maximized by exploring $\alpha^{\varphi}_i(n)\in[0,1]$ with all the other parameters given.
Now note that the right hand side of \eqref{eqn:add_1} is of the form $A+\frac{B}{C\alpha^{\varphi}_i(n) + D}$ ($A,B,C,D\in \mathbb{R}$), which
is either increasing or decreasing in $\alpha^{\varphi}_i(n)$ or remains constant.
Hence,  for all $n\in[|\mathscr{N}_i|-2]$, the right hand side of \eqref{eqn:add_1} can be optimized either by taking $\alpha^{\varphi}_i(n) =0$ or by taking $\alpha^{\varphi}_i(n) =1$.

We then consider the action variable for  state $0$.
For any policy $\varphi$ with $\alpha^{\varphi}_i(0)= 0$, $g=0$, because the MDP will stay in state $0$ all the time.
Thus the average revenue $g$ for the optimal policy $\varphi$ must be non-negative. 
From Bellman's equation, $\alpha^{\varphi}_i(0)$ under this optimal policy $\varphi$ is either $0$ or maximizes
\begin{equation}\label{eqn:add_2}
\max\limits_{\alpha^{\varphi}_i(0) \in (0,1]} \left\{\frac{-g}{\alpha_i^{\varphi}(0)\lambda_{\ell(i)}}+ V(1)\right\}.
\end{equation}
Since $g\geq 0$, $\alpha^{\varphi}_i(0) = 1$ will always maximize the bracketed term in \eqref{eqn:add_2}.
It follows that there exists an optimal policy $\varphi$ with $\alpha^{\varphi}_i(0) \in \{0,1\}$ that maximizes the long-run average revenue of this MDP.
This proves the lemma.

\endproof

\proof{Proof of Proposition~\ref{prop:a_opt}.}

From Lemma~\ref{lemma:deterministic},   for all $n=0,1,\ldots,|\mathscr{N}_i|-2$,
there exists a policy $\varphi$ 
determined by $\alpha_i^{\varphi}(n)\in\{0,1\}$ that maximizes the long-run average revenue of this MDP.

For any $\varphi\in\Phi_i$ and $n=0,1,\ldots,|\mathscr{N}_i|-2$, if $\alpha^{\varphi}_i(n) = 1$, there is a positive transition rate from state $n$ to state $n+1$; otherwise, that transition rate is $0$.
The transition rate from $n$ to $n-1$, if $n\geq 1$, is always positive and independent of the policy employed.

For any $\varphi\in\Phi_{i}$ and $0\leq n \leq |\mathscr{N}_{i}|-1$, if $\alpha^{\varphi}_{i}(n)=0$ then there is no difference between making $\alpha^{\varphi}_{i}(n')=0$ or $1$ when $n'>n$, because we cannot reach state $n'$ if we start below state $n$.
For any $\varphi\in\Phi_{i}$ and $0\leq n < n' \leq |\mathscr{N}_{i}|-1$, this allows us to assume that if $\alpha^{\varphi}_{i}(n)=0$ then $\alpha^{\varphi}_{i}(n') = 0$, which makes $\varphi$ a \emph{threshold policy}.
For $\varphi\in \Phi_{i}$, we define
\begin{displaymath}
m^{\varphi} =\begin{cases}
0 ,& \text{if }\alpha^{\varphi}_{i}(0) = 0,\\
m,& \text{if } \alpha^{\varphi}_{i}(m-1) = 1,\ \alpha^{\varphi}_{i}(m) = 0,\ m\in[|\mathscr{N}_i|-1],
\end{cases}
\end{displaymath}
so that $m-1$ is the maximum value of $n$ for which $\alpha_i^{\phi}(n)=1$.

Let $\pi^{m}_{n}$ represent the steady state probability of state $n\in \mathscr{N}_{i}$ under policy $\varphi$ with $m^{\varphi}=m$.
We maximize the right hand side of \eqref{eqn:dual_func_i} 
with respect to all threshold policies defined by $m$ taking into account the specific form of the stationary distribution.
For an optimal solution $\varphi^{*}$, we obtain 
\begin{equation}\label{app:solve_opt_lr_3:1}
m^{\varphi^{*}}\in \arg\max\limits_{m\in[|\mathscr{N}_i|-1]}\Biggl\{\pi^{m}_{0}\sum\limits_{n=1}^{m}\biggl[\frac{(\lambda_{\ell})^{n}}{n! (\mu_{i})^{n}}\Bigl(r_{i}(n,\pmb{\gamma}) - \frac{n\mu_{i}}{\lambda_{\ell}}\bigl(\nu_{\ell}+\sum\limits_{j\in \mathscr{J}_i}w_{j,i}\gamma_{j}\bigr)\Bigr)\biggr]\Biggr\}
\end{equation}
where 
$r_{i}(n,\pmb{\gamma}) = n\mu_{i}r_{\ell}-\sum\limits_{j\in\mathscr{J}_i}\varepsilon_{j}w_{j,i}n - \sum\limits_{j\in\mathscr{J}_i}\gamma_{j}w_{j,i}n$.
Equation \eqref{app:solve_opt_lr_3:1} can be rewritten as
\begin{equation}\label{app:solve_opt_lr_3:2}
m^{\varphi^{*}}\in \arg\max\limits_{m} \Biggl\{(\lambda_{\ell} \tilde{r}_{i} - \nu_{\ell})\frac{\sum\limits_{n=0}^{m-1}\frac{(\lambda_{\ell})^n}{n! (\mu_{i})^{n}}}{\sum\limits_{n=0}^{m}\frac{(\lambda_{\ell})^n}{n! (\mu_{i})^n}}\Biggr\}
\end{equation}
where
$\tilde{r}_{i} =r_{\ell}-\frac{1}{\mu_{i}}\sum\limits_{j\in\mathscr{J}_i}\varepsilon_{j}w_{j,i}-\frac{1}{\mu_{i}}\sum\limits_{j\in\mathscr{J}_i}\gamma_{j}w_{j,i} - \frac{1}{\lambda_{\ell}}\sum\limits_{j\in\mathscr{J}_i}\gamma_{j}w_{j,i}$.
Note that the right hand side of \eqref{app:solve_opt_lr_3:2} may return a set with more than one element: the optimal value for $m^{\varphi^{*}}$ is not necessarily unique. 
If this is the case, we choose any of the possible maxima. 
From \eqref{app:solve_opt_lr_3:2}, if $\nu_{\ell}< \lambda_{\ell}\tilde{r}_{i}$, then $m^{\varphi^{*}} = |\mathscr{N}_{i}|-1$; if $\nu_{\ell}= \lambda_{\ell}\tilde{r}_{i}$, then $m^{\varphi^{*}}$ can be any value in $\mathscr{N}_i$; otherwise, $m^{\varphi^{*}} = 0$. We thus interpret $\varphi^{*}$ as an index policy where indices for states $n\in\mathscr{N}_{i}\backslash\{|\mathscr{N}_i|-1\}$ are given by $\lambda_{\ell}\tilde{r}_{i}$.
This proves Proposition~\ref{prop:a_opt}.

\endproof

\section{Proof of Proposition~\ref{prop:solution_existence}}\label{app:prop:solution_existence}

\proof{Proof.}
For given $\pmb{\gamma}_0\in\mathbb{R}_0^{J}$, $o\in\mathscr{O}(\pmb{\gamma}_0)$ and $\iota\in\mathscr{I}(o)$, from the definition of $\Xi_i(\cdot)$ in \eqref{eqn:index_value}, we rewrite  \eqref{eqn:necessary_gamma} in the form
\begin{equation}\label{eqn:solution_existence:1}
\Bigl(1+\frac{\lambda_{\ell_{\iota}}}{\mu_{i_{\iota}}}\Bigr)\sum\limits_{j\in[J]}w_{j,i_{\iota}}\gamma_{1,j} = \Xi_{i_{\iota}}(\bm{0})-\nu_{\ell_{\iota}}(o,\pmb{\gamma}_0).
\end{equation}

Note that, from \eqref{eqn:necessary_gamma:zero}, $\gamma_{1,j}=0$ if there is no $\iota \in\mathscr{I}(o)$ with $j=j_{\iota}$.
For the remaining $|\mathscr{I}(o)|$ entries in $\pmb{\gamma}_1$ where the $\gamma_{1,j}$ are allowed to be non-zero, we construct a $|\mathscr{I}(o)|\times |\mathscr{I}(o)|$ matrix $\mathcal{M} = (m_{\iota,j})$ and write \eqref{eqn:solution_existence:1} as $\mathcal{M}\tilde{\pmb{\gamma}}_1 = \bm{y}$ with $m_{\iota,j}=(1+\lambda_{\ell_{\iota}}/\mu_{i_{\iota}})w_{j,i_{\iota}}$, $\tilde{\gamma}_{1,\iota}=\gamma_{1,j_{\iota}}$ and $y_{\iota} = \Xi_{i_{\iota}}(\bm{0})-\nu_{\ell_{\iota}}(o,\pmb{\gamma}_0)$ for all  $\iota \in\mathscr{I}(\pmb{\gamma}_0)$.

As described in \ref{case:equality:capacity}, under policy $\bar{\varphi}(o)$ produced by Algorithm~\ref{algo:varphi_gamma}, if a resource pool $j\in[J]$ becomes a critical pool of any critical pair  $\iota\in\mathscr{I}(o)$, all PS pairs  $\iota'\in[N]$ with $j\in\mathscr{J}_{i_{\iota'}}$ and lower priorities than that of $\iota$ are removed from the list of candidate PS pairs; that is, all these lower prioritized pairs $\iota'$ will not become critical pairs.
In other words, by appropriately re-ordering columns of $\mathcal{M}$ and entries of $\bm{y}$, it becomes an upper-triangular matrix with $|\mathcal{M}|=\sum_{\iota\in\mathscr{I}(o)}(1+\lambda_{\ell_{\iota}}/\mu_{i_{\iota}})w_{j_{\iota},i_{\iota}}$.
Because $w_{j_{\iota},i_{\iota}} >0$ for all $\iota\in\mathscr{I}(o)$, $|\mathcal{M}|>0$. This proves the proposition.

\endproof

\section{Proof of Proposition~\ref{prop:converge_gamma}}\label{app:prop:converge_gamma}

We start with a lemma.

\begin{lemma}\label{lemma:non-negative}
For given $\pmb{\gamma}_0\in\mathbb{R}_0^J$ and $o\in\mathscr{O}$,
if $\mathcal{T}^o(\pmb{\gamma}_0)\in\mathbb{R}_0^J$, then
the complementary slackness conditions corresponding to the capacity constraints (that is, Equation~\eqref{eqn:relaxconstraint:slack}) and Equation~\eqref{eqn:a_opt:b} are satisfied by taking $\phi=\varphi=\bar{\varphi}(o)$ and $\pmb{\gamma}=\mathcal{T}^o(\pmb{\gamma}_0)$.
\end{lemma}
The proof of Lemma~\ref{lemma:non-negative} will be given in Appendix~\ref{app:lemma:non-negative} in the e-companion. 
Further, if \eqref{eqn:a_opt:a} and \eqref{eqn:a_opt:c} are satisfied by setting $\pmb{\gamma}=\mathcal{T}^o(\pmb{\gamma}_0)$ and $\varphi=\bar{\varphi}(o)$, then $\mathcal{T}^o(\pmb{\gamma}_0)$ is a vector of decomposable multipliers  and, from the strong duality theorem, the policy $\bar{\varphi}(o)$ is an optimal policy for the relaxed problem.

\proof{Proof of Proposition~\ref{prop:converge_gamma}.}
For given $\pmb{\gamma}_0\in\mathbb{R}_0^J$ and $o\in\mathscr{O}(\pmb{\gamma}_0,\bm{0})$, if $\mathcal{T}^o(\pmb{\gamma}_0)=\pmb{\gamma}_0$, then \eqref{eqn:a_opt:a} and \eqref{eqn:a_opt:c} are satisfied  by setting $\pmb{\gamma}=\pmb{\gamma}_0$ and $\varphi=\bar{\varphi}(o)$, because in this case $\mathscr{O}(\pmb{\gamma}_0,\bm{0})=\mathscr{O}(\mathcal{T}^o(\pmb{\gamma}_0),\bm{0})$. That is, this $\mathcal{T}^o(\pmb{\gamma}_0)=\pmb{\gamma}_0$ is a vector of decomposable multipliers.

\endproof

\section{Proof of Corollary~\ref{coro:converge_gamma}}\label{app:coro:converge_gamma}
\proof{Proof.}
For all $\iota\in\mathscr{I}(o)$, from the definition of Algorithm~\ref{algo:varphi_gamma}, $\nu_{\ell_{\iota}}(o,\mathcal{T}^o(\pmb{\gamma}_0))=\Xi_{i_{\iota}}(\mathcal{T}^o(\pmb{\gamma}_0),\bm{0})$. For all $j \notin \{j_{\iota}(o)\in[J] ~|~\iota\in\mathscr{I}(o)\}$, $\mathcal{T}^o_j(\pmb{\gamma}_0)=0$ because $\mathcal{T}^o(\pmb{\gamma}_0)$ is the solution of \eqref{eqn:necessary_gamma}-\eqref{eqn:necessary_gamma:zero}.
In other words, $\mathcal{T}^o(\mathcal{T}^o(\pmb{\gamma}_0))=\mathcal{T}^o(\pmb{\gamma}_0)$ is a solution of \eqref{eqn:necessary_gamma}-\eqref{eqn:necessary_gamma:zero} by taking $\pmb{\gamma}_0=\mathcal{T}^o(\pmb{\gamma}_0)$ and, from Proposition~\ref{prop:solution_existence}, this solution is unique.
Together with Proposition~\ref{prop:converge_gamma}, this yields that $\mathcal{T}^o(\pmb{\gamma}_0)$ is a decomposable vector and $\bar{\varphi}(o)$ is optimal for the relaxed problem. 

\endproof

\section{Proof of Lemma~\ref{lemma:non-negative}}\label{app:lemma:non-negative}

\proof{Proof.}
Under policy $\bar{\varphi}(o)$, a critical resource pool $j\in\{j_{\iota}(o)\in[J] | \iota\in\mathscr{I}(o)\}$ is always fully occupied; that is,
\begin{equation}
\bm{\omega}_j\cdot\Bigl(\Pi^{\bar{\varphi}(o)}_n + \Pi^{\bar{\varphi}(o)}_a\Bigr) = C_j,~\forall j\in\{j_{\iota}(o)\in[J] | \iota\in\mathscr{I}(o)\}.
\end{equation}
Equation~\eqref{eqn:relaxconstraint:slack} is achieved for $j\in\{j_{\iota}(o)\in[J] | \iota\in\mathscr{I}(o)\}$. For the other resource pools $j\notin\{j_{\iota}(o)\in[J] | \iota\in\mathscr{I}(o)\}$, \eqref{eqn:relaxconstraint:slack} is satisfied by taking $\gamma_j=\gamma_{1,j}=0$, which follows from \eqref{eqn:necessary_gamma:zero}.
Also, for policy $\bar{\varphi}(o)$, \eqref{eqn:a_opt:b} is guaranteed by \eqref{eqn:necessary_gamma} and procedure \ref{case:equality:action} of Algorithm~\ref{algo:varphi_gamma}. This proves the proposition.
\endproof

\section{Proof of Proposition~\ref{prop:equal_opt}}\label{app:prop:equal_opt}

\IncMargin{1em}
\begin{algorithm}
\linespread{0.6}\selectfont

\SetKwFunction{FPriorityPolicyForSubProblem}{PriorityPolicyForSubProblem}
\SetKwProg{Fn}{Function}{:}{\KwRet}
\SetKwInOut{Input}{Input}\SetKwInOut{Output}{Output}
\SetAlgoLined
\DontPrintSemicolon

\Input{Real vectors $\bm{\nu}\in\mathbb{R}^L$ and $\pmb{\gamma}\in\mathbb{R}_0^J$ and a ranking of pairs $o\in\mathscr{O}(\bm{\nu},\pmb{\gamma},\bm{\beta})$ with $\bm{\beta}\in\mathbb{R}_0^{I}$.}
\Output{A policy $\varphi(o,\bm{\nu},\pmb{\gamma})$ determined by $\bm{\alpha}_i^{\varphi(o,\bm{\nu},\pmb{\gamma})}$ for sub-problems $i\in[I]$.}
\Fn{\FPriorityPolicyForSubProblem{$\bm{\nu},\pmb{\gamma},o$}}{

    $\bm{\alpha}^{\varphi}_i\gets \bm{0}$ for all $i\in[I]$\;

    Initializing a list of candidate PS pairs as the list of all PS pairs.\;

    $\iota \gets 0$ \tcc*{Iteration variable}

    \While {$\iota<N $ {\bf and} the list of candidate PS pairs is not empty}{	

	 		$\iota\gets \iota+1$\;

			\If {PS pair $\iota$ is not in the list of candidate PS pairs}{
			    {\bf continue}\;
			}

		   $a_1 \gets \inf\left\{\left\{\alpha^{\varphi}_{i_{\iota}}(n_{\iota})\in[0,1]\left|~\exists j\in[J], \bm{\omega}_j\cdot\left(\bm{\Pi}^{\varphi}_n+\bm{\Pi}^{\varphi}_a\right)=C_j\right.\right\}\cup\{1\}\right\}$\;
		   
		   $\mathscr{A} \gets \left\{\alpha^{\varphi}_{i_{\iota}}(n_{\iota}) \in [0,1] \left|~\sum_{i\in\mathscr{P}_{\ell_{\iota}}}\bm{\pi}^{\varphi}_i\cdot \bm{\alpha}^{\varphi}_i = 1\right.\right\}$\;

	       $\alpha^{\varphi}_{i_{\iota}}(n_{\iota})\gets a_1$\;
		\tcc*{Update $\alpha^{\varphi}_{i_{\iota}}(n_{\iota})$ with the maximal activating probability}
		\tcc*{without violating any capacity constraint.}
		
		\If{$\Xi_{i_{\iota}}(\pmb{\gamma},\bm{0}) = \nu_{\ell_{\iota}}$ and $\mathscr{A}\neq \emptyset$}{
			$\alpha^{\varphi}_{i_{\iota}}(n_{\iota})\gets \inf \mathscr{A}$\;
			\tcc*{When $\Xi_{i_{\iota}}(\pmb{\gamma},\bm{0}) = \nu_{\ell_{\iota}}$ and the action constraint can be satisfied by setting $\alpha^{\varphi}_{i_{\iota}}(n_{\iota})\in[0,1]$, update $\alpha^{\varphi}_{i_{\iota}}(n_{\iota})$.}
		}

		   \uIf {$\Xi_{i_{\iota}}(\pmb{\gamma},\bm{0}) < \nu_{\ell_{\iota}}$}{
		 		$\alpha^{\varphi}_{i_{\iota}}(n_{\iota})\gets 0$\;

	    		remove all PS pairs $\iota'>\iota$ with $\ell_{\iota'}=\ell_{\iota}$ from the list of candidate PS pairs.\;
			}\ElseIf {$\exists j\in[J],\ \bm{\omega}_j\cdot\left(\bm{\Pi}^{\varphi}_n+\bm{\Pi}^{\varphi}_a\right)=C_j$}{
			   
	 		   remove all PS pairs $\iota'>\iota$ with $w_{j,i_{\iota'}}>0$ from the list of candidate PS pairs.\;
			}

    }

    $\bm{\alpha}^{\varphi(o,\bm{\nu},\pmb{\gamma})}_i\gets \bm{\alpha}^{\varphi}_{i}$ for all $i\in[I]$\;
}
\caption{Priority-style policy for the sub-problems}\label{algo:varphi_nu_gamma}
\end{algorithm}
 \DecMargin{1em}

From Proposition~\ref{prop:converge_gamma}, a sufficient condition for the existence of decomposable multipliers is that there exist $\pmb{\gamma}_0\in\mathbb{R}_0^J$ and $o\in\mathscr{O}(\pmb{\gamma}_0,\bm{0})$ such that $\mathcal{T}^o(\pmb{\gamma}_0)=\pmb{\gamma}_0$.
Such a $\pmb{\gamma}_0$ is the vector of decomposable multipliers; from the definition of $\Xi_i(\cdot)$ in \eqref{eqn:index_value}, Equation \eqref{eqn:necessary_gamma} can be rewritten in the form of \eqref{eqn:solution_existence:1}. 

Similar to Section~\ref{subsec:decomposable},
for given $\bm{\nu}\in\mathbb{R}^{L}$, $\pmb{\gamma}\in\mathbb{R}^{J}_0$ and a positive vector $\bm{\beta}\in\mathbb{R}_{0}^{I}$, let states $n\in\mathscr{N}_i$, $i\in\mathscr{P}_{\ell}$, $\ell\in[L]$, be ranked according to the descending order of $\beta_{i}(\Xi_i(\pmb{\gamma},\bm{0})-\nu_{\ell})$, and let $\mathscr{O}(\bm{\nu},\pmb{\gamma},\bm{\beta})$ represent the set of all such rankings.
Note that, in Algorithm~\ref{algo:varphi_gamma} and Section~\ref{subsec:decomposable}, we ordered the PS pairs $\iota\in[N]$ according to $\Xi_{i_{\iota}}(\pmb{\gamma},\bm{0})$, but here we order them according to $\beta_{i_{\iota}}(\Xi_{i_{\iota}}(\pmb{\gamma},\bm{0})-\nu_{\ell_{\iota}})$ for given $\bm{\beta}\in\mathbb{R}_0^{I}$ and $\bm{\nu}\in\mathbb{R}^L$. 

As in Algorithm~\ref{algo:varphi_gamma}, Algorithm~\ref{algo:varphi_nu_gamma} generates a policy $\varphi(o,\bm{\nu},\pmb{\gamma})\in\tilde{\Phi}$ with PS-pair priorities defined by $o\in\mathscr{O}(\bm{\nu},\pmb{\gamma},\bm{\beta})$, satisfying 
the relaxed capacity constraints described in \eqref{eqn:constraint:relax:resources}.
The policy $\varphi(o,\bm{\nu},\pmb{\gamma})$ is dependent on the given $\bm{\beta}$ only through the ranking $o$.
The main idea of Algorithm~\ref{algo:varphi_nu_gamma}, as in Algorithm~\ref{algo:varphi_gamma}, is to activate PS pairs consecutively  according to their order in the ranking $o$ with $\bm{\alpha}^{\varphi(o,\bm{\nu},\pmb{\gamma})}_i$ initialized to be $\bm{0}$ for all $i\in[I]$. When a relaxed capacity constraint achieves equality by activating or partially activating the PS pair $\iota\in[N]$ (the $\iota$th pair in ranking $o$), all parameters are maintained in the same way as described in \ref{case:equality:capacity}.
\begin{enumerate}[label=(\Roman*)]
\setcounter{enumi}{2}
\item If the given multiplier $\nu_{\ell_{\iota}}$ for the $\iota$th pair in ranking $o$ is greater than $\Xi_{i_{\iota}}(\pmb{\gamma},\bm{0})$, then
all PS pairs $\iota'\geq \iota$ with $\ell_{\iota'}=\ell_{\iota}$ are removed from the list of candidate PS pairs awaiting activation. \label{case:nu:action}
\end{enumerate}
Note that, during the process of generating policy $\varphi(o,\bm{\nu},\pmb{\gamma})$, we do not necessarily stop activating PS pairs if the constraints in \eqref{eqn:constraint:relax:action} achieve equality; that is, here, the mechanism described in \ref{case:equality:action} is replaced by that described in \ref{case:nu:action}.
The mechanism described in \ref{case:nu:action}  guarantees that for each $\ell\in[L]$, PS pairs $\iota\in[N]$ with $\ell_{\iota}=\ell$ are passive if $\Xi_{i_{\iota}}(\pmb{\gamma},\bm{0})-\nu_{\ell}<0$; and are partially or fully active, otherwise.
Also, for the non-passive PS pairs $\iota$ satisfying  $\Xi_{i_{\iota}}(\pmb{\gamma},\bm{0})-\nu_{\ell_{\iota}} \geq 0$, the parameter $\bm{\beta}$ will affect their priorities for activation, which together with the specified ranking $o\in\mathscr{O}(\bm{\nu},\pmb{\gamma},\bm{\beta})$ determine 
whether they are fully or partially activated.


For any ranking of $\mathscr{O}(\bm{\nu},\pmb{\gamma},\bm{\beta})$, we stipulate that the PS pairs for the boundary states $|\mathscr{N}_i|-1$ of patterns $i\in[I]\backslash\{d(\ell):\ell\in[L]\}$ are the lowest ranked, so that they are always passive.

Under the policy $\varphi(o,\bm{\nu},\pmb{\gamma})$, the set of passive PS pairs $\iota$ with $\ell_{\iota}=\ell$ ($\ell\in[L]$) increases monotonically from $\{\iota\in[N]|~n_{\iota}=|\mathscr{N}_{i_{\iota}}|-1,~\ell_{\iota}=\ell\}$ to $\{\iota\in[N] | \ell_{\iota}=\ell\}$ as $\nu_{\ell}$ increases from $-\infty$ to $+\infty$. So the sum of the expected action variables, $\sum_{i\in\mathscr{P}_{\ell}}\bm{\pi}^{\varphi(o,\bm{\nu},\pmb{\gamma})}_i\cdot\bm{\alpha}^{\varphi(o,\bm{\nu},\pmb{\gamma})}_i$,
decreases with respect to $\nu_{\ell}$ from 1 to 0.



We focus on the first $N-I+L$ PS pairs represented by $\iota\in[N-I+L]$, which are allowed to be active.
Let $\mathscr{I}(o,\bm{\nu},\pmb{\gamma})$ represent the set of critical pairs and $j_{\iota}(o,\bm{\nu},\pmb{\gamma})$ represent the critical resource pool of critical pair $\iota\in\mathscr{I}(o,\bm{\nu},\pmb{\gamma})$  with respect to policy $\varphi(o,\bm{\nu},\pmb{\gamma})$. Define 
\begin{equation}
\mathscr{I}_{1}(o,\bm{\nu},\pmb{\gamma}) \coloneqq \left\{\iota \in [N-I+L] | \alpha^{\varphi(o,\bm{\nu},\pmb{\gamma})}_{i_{\iota}}(n_{\iota}) > 0 \text{ and }\iota\notin \mathscr{I}(o,\bm{\nu},\pmb{\gamma})\right\}
\end{equation}
and
\begin{equation}
\mathscr{I}_{0}(o,\bm{\nu},\pmb{\gamma}) \coloneqq \left\{\iota \in [N-I+L] | \alpha^{\varphi(o,\bm{\nu},\pmb{\gamma})}_{i_{\iota}}(n_{\iota}) = 0\right\}.
\end{equation}

For $\bm{\nu}\in\mathbb{R}^{L}$, $\pmb{\gamma}\in\mathbb{R}_0^{J}$, and $\iota\in[N-I+L]$, 
define $\bm{\Delta}\in \mathbb{R}^{[N-I+L]}$ with
\begin{equation}\label{eqn:linear_func}
\Delta_{\iota} \coloneqq \Xi_{i_{\iota}}(\bm{0},\bm{0})-\nu_{\ell_{\iota}} - \left(1+\frac{\lambda_{\ell_{\iota}}}{\mu_{i_{\iota}}}\right)\sum\limits_{j\in[J]}w_{j,i_{\iota}}\gamma_{j}.
\end{equation}
We introduce a condition.
\begin{condition}{Partial Decomposability}\label{cond:lemma:opt:for_any_nu}
We say the system is partially decomposable if and only if,
for any given $\bm{\nu}\in\mathbb{R}^L$, there exist $\pmb{\gamma} \in \mathbb{R}^{J}_0$, $\bm{\beta}\in\mathbb{R}_0^I$, $o\in\mathscr{O}(\bm{\nu},\bm{0},\bm{\beta})$ and $\bm{\Delta}\in\mathbb{R}^{[N-I+L]}$ that satisfy
\begin{equation}\label{eqn:Delta}
\Delta_{\iota}
\begin{cases}
\geq 0 & \text{if }\iota\in \mathscr{I}_1(o,\bm{\nu},\bm{0}),\\
= 0 & \text{if }\iota\in \mathscr{I}(o,\bm{\nu},\bm{0}),\\
\leq  0 & \text{if }\iota\in\mathscr{I}_0(o,\bm{\nu},\bm{0}).
\end{cases}
\end{equation}
and
\begin{equation}\label{eqn:linear_func:zero}
\gamma_{j} = 0
\end{equation} 
for all $j\in [J]\backslash \{j_{\iota}(o,\bm{\nu},\bm{0}) | \exists \iota\in\mathscr{I}(o,\bm{\nu},\bm{0})\}$.
\end{condition}

If the system is partially decomposable, then, fine-tune the $\bm{\nu}\in\mathbb{R}^{L}$ such that, for the corresponding $\pmb{\gamma}\in\mathbb{R}^J_0$, $\bm{\beta}\in\mathbb{R}^I_0$ and $o\in\mathscr{O}(\bm{\nu},\bm{0},\bm{\beta})$ satisfying \eqref{eqn:Delta} and \eqref{eqn:linear_func:zero}, the relaxed action constraints described in \eqref{eqn:constraint:relax:action} are satisfied by the policy $\varphi(o,\bm{\nu},\bm{0})$.
In this context, the policy  $\varphi(o,\bm{\nu},\bm{0})$ generated by Algorithm~\ref{algo:varphi_nu_gamma} coincides with the policy $\bar{\varphi}(o)$ generated by Algorithm~\ref{algo:varphi_gamma} with the same PS-pair ranking $o$, and the $\nu_{\ell}(o,\pmb{\gamma})$ ($\ell\in[L]$) generated by Algorithm~\ref{algo:varphi_gamma} is equal to the $\nu_{\ell}$.
The existence of the $\pmb{\gamma} $, $\bm{\beta}$ and the ranking $o\in\mathscr{O}(\bm{\nu},\bm{0},\bm{\beta})$ is guaranteed by partial decomposability condition.
Also, for the ranking $o$ and the multipliers $\nu_{\ell}(o,\pmb{\gamma})=\nu_{\ell}$ ($\ell\in[L]$), $\pmb{\gamma}$ is non-negative and a fixed point of the function $\mathcal{T}^o$ (that is, $\mathcal{T}^o(\pmb{\gamma})=\pmb{\gamma}$).
Together with Proposition~\ref{prop:converge_gamma},  
this yields that, if the system is partially decomposable, then it is  decomposable.

\begin{lemma}\label{lemma:opt:for_any_nu}
If the system is weakly coupled then it is partially decomposable.
With the vector $\pmb{\gamma}$ in the definition of partial decomposability,
for $j\in[J]$,
\begin{enumerate}[label=\roman*)]
\item 
if there is a critical PS pair $\iota\in\mathscr{I}(o,\bm{\nu},\bm{0})$ with critical resource pool $j=j_{\iota}(o,\bm{\nu},\bm{0})$, and no 
$j'\neq j$ with $j'\in\mathscr{J}_{i_{\iota}}$ is critical for any other PS pair $\iota' \in\mathscr{I}(o,\bm{\nu},\bm{0})$, then
\begin{equation}\label{eqn:opt:for_any_nu:a}
\gamma_j = \frac{\Xi_{i_{\iota}}(\bm{0},\bm{0})-\nu_{\ell_{\iota}}}{w_{j,i_{\iota}}\bigl(1+\frac{\lambda_{\ell_{\iota}}}{\mu_{i_{\iota}}}\bigr)};
\end{equation}
\item 
if there are critical PS pairs $\iota$ and $\iota'$ in $\mathscr{I}(o,\bm{\nu},\bm{0})$ with critical resource pools $j=j_{\iota}(o,\bm{\nu},\bm{0})\neq j_{\iota'}(o,\bm{\nu},\bm{0})$ and $ j_{\iota'}(o,\bm{\nu},\bm{0})\in\mathscr{J}_{i_{\iota}}$, then
\begin{equation}\label{eqn:opt:for_any_nu:b}
\gamma_j =\frac{w_{j_{\iota'}(o,\bm{\nu},\bm{0}),i_{\iota}}}{w_{j,i_{\iota}}} \Bigl(\frac{\Xi_{i_{\iota}}(\bm{0},\bm{0})-\nu_{\ell_{\iota}}}{w_{j_{\iota'}(o,\bm{\nu},\bm{0}),i_{\iota}}\bigl(1+\frac{\lambda_{\ell_{\iota}}}{\mu_{i_{\iota}}}\bigr)}
-\frac{\Xi_{i_{\iota'}}(\bm{0},\bm{0})-\nu_{\ell_{\iota'}}}{w_{j_{\iota'}(o,\bm{\nu},\bm{0}),i_{\iota'}}\bigl(1+\frac{\lambda_{\ell_{\iota'}}}{\mu_{i_{\iota'}}}\bigr)}\Bigr);
\end{equation}
\item otherwise,
\begin{equation}\label{eqn:opt:for_any_nu:c}
\gamma_j=0.
\end{equation}
\end{enumerate}

\end{lemma}

\proof{Proof.}

For 
PS pairs $\iota\in [N-I+L]$, if $\Xi_{i_{\iota}}(\bm{0},\bm{0})<\nu_{\ell_{\iota}}$, 
then $\iota\in\mathscr{I}_0(o,\bm{\nu},\bm{0})$, and for any $\pmb{\gamma}\in\mathbb{R}^J_0$,
there always exist $\Delta_{\iota}\leq 0$ satisfying \eqref{eqn:linear_func}. 

We focus, then, on equations \eqref{eqn:linear_func} corresponding to the remaining PS pairs.
Let $\mathscr{N}^0 = \{\iota\in[N-I+L] |\Xi_{i_{\iota}}(\bm{0},\bm{0})<\nu_{\ell_{\iota}}\}$, and $\sigma(\iota)$ represent the position of PS pair $\iota \in[N-I+L]\backslash\mathscr{N}^0$ among all pairs in the set $[N-I+L]\backslash\mathscr{N}^0$ according to the ranking $o\in\mathscr{O}(\bm{\nu},\bm{0},\bm{\beta})$.
The notation $\sigma(\iota)$  clarifies the difference between $\sigma(\iota)$th PS pair among $[N-I+L]\backslash \mathscr{N}^0$ and the PS pair $\iota$ defined previously.

We construct an $(N-I+L-|\mathscr{N}^0|)\times (N-I+L-|\mathscr{N}^0|)$ matrix $\mathcal{M} = (m_{i,j})$ and write \eqref{eqn:linear_func} 
as $\mathcal{M}\bm{x}=\bm{y}$ ($\bm{x},\bm{y}\in\mathbb{R}^{N-I+L-|\mathscr{N}^0|}$), where $y_{\sigma(\iota)}\coloneqq\Xi_{i_{\iota}}(\bm{0},\bm{0})-\nu_{\ell_{\iota}}$ ($\iota\in[N-I+L]\backslash\mathscr{N}^0$).
For given $\bm{\nu}\in\mathbb{R}^{L}$, $\bm{\beta}\in\mathbb{R}_0^{I}$ and $o\in\mathscr{O}(\bm{\nu},\bm{0},\bm{\beta})$, 
$\mathcal{M}$ and $\bm{x}$ are defined by:
\begin{enumerate}[label=\alph*)]
\item if PS pair $\iota\in\mathscr{I}_{1}(o,\bm{\nu},\bm{0})$, then $x_{\sigma(\iota)}$ equals $\Delta_{\iota}$, 
$m_{\sigma(\iota),\sigma(\iota)} = 1$ and $m_{\sigma(\iota'),\sigma(\iota)}=0$ for all $\iota'\neq \iota$, $\iota'\in[N-I+L]\backslash\mathscr{N}^0$; 
\label{case1}
\item if PS pair $\iota\in \mathscr{I}(o,\bm{\nu},\bm{0})$,
then $x_{\sigma(\iota)}$ equals $\gamma_{j_{\iota}(o,\bm{\nu},\bm{0})}$ and $m_{\sigma(\iota'),\sigma(\iota)}$, $\iota'\in[N-I+L]\backslash\mathscr{N}^0$, are equal to $w_{j_{\iota}(0,\bm{\nu},\bm{0}),i_{\iota'}}(1+\frac{\lambda_{\ell_{\iota'}}}{\mu_{i_{\iota'}}})$, the coefficients of $\gamma_{j_{\iota}(o,\bm{\nu},\bm{0})}$ in \eqref{eqn:linear_func}; and \label{case2}
\item for all the other PS pairs $\iota\in \mathscr{I}_0(o,\bm{\nu},\bm{0})$, 
$x_{\sigma(\iota)}$ equals $-\Delta_{\iota}$, entry $m_{\sigma(\iota),\sigma(\iota)} = -1$ and $m_{\sigma(\iota'),\sigma(\iota)}=0$ for all $\iota'\neq \iota$, $\iota'\in[N-I+L]\backslash\mathscr{N}^0$.\label{case4}
\end{enumerate}
Also, 
we set $\gamma_{j}=0$ for all $j\in [J]\backslash \{j_{\iota}(o,\bm{\nu},\bm{0})|\exists \iota\in\mathscr{I}(o,\bm{\nu},\bm{0})\}$.
Then, 
the lemma holds if and only if, for any weakly coupled system, there exists a non-negative solution to $\mathcal{M}\bm{x}=\bm{y}$, for which \eqref{eqn:opt:for_any_nu:a}-\eqref{eqn:opt:for_any_nu:c} are satisfied.

For a given $i \in [I]$, from  \eqref{eqn:linear_func}, $\Delta_{\iota} = \Delta_i$ for all PS pairs $\iota$ with $n_{\iota}\in\mathscr{N}_{i}\backslash\{|\mathscr{N}_{i}|-1\}$.
In this context, when we consider a solution to $\mathcal{M}\bm{x} = \bm{y}$ satisfying \eqref{eqn:linear_func},
we only need to consider one of these PS pairs $\iota$ with $n_{\iota}\in\mathscr{N}_{i}\backslash\{|\mathscr{N}_{i}|-1\}$ associated with one row of $\mathcal{M}$.
Thus, we remove the other $|\mathscr{N}_{i}|-2$ linear functions by removing corresponding rows and columns of $\mathcal{M}$ and elements of $\bm{x}$ and $\bm{y}$.
In particular, for all critical pairs $\iota\in\mathscr{I}(o,\bm{\nu},\bm{0})$, we keep the row and column associated with PS pair $\iota$ in $\mathcal{M}$ and remove the $|\mathscr{N}_{i_{\iota}}|-2$ others.
We represent the PS pair associated with the remaining linear function for pattern $i$ by $\iota_i$.


Removing these unnecessary rows and columns, 
we reformulate $\mathcal{M}\bm{x} = \bm{y}$ as $\tilde{\mathcal{M}}\tilde{\bm{x}} = \tilde{\bm{y}}$, where $\tilde{\mathcal{M}} \in \mathbb{R}^{K\times K}$ and $\tilde{\bm{x}},\tilde{\bm{y}}\in\mathbb{R}^{K}$ and $K\leq I$ is the number of patterns $i\in [I]$ with $\Xi_i(\bm{0},\bm{0}) \geq \nu_{\ell(i)}$.

With $m$ the number of critical pairs/pools, for $k\in [m]$, let
\begin{equation}\label{eqn:sub-matrix}
\mathcal{M}_k = \left[
\begin{array}{ M{3em}M{3em}M{3em}M{12em} M{3em}M{3em}M{3em}}
 1 &            &     & $w_{j_k,i(1)}(1+\frac{\lambda_{\ell(i(1))}}{\mu_{i(1)}})$ &&&\\
    & $\ddots$ &      &$ \vdots$ &&&\\
   &           &  1 &  $w_{j_k,i(m_k^1)}(1+\frac{\lambda_{\ell(i(m_k^1))}}{\mu_{i(m_k^1)}}) $ &&&\\
&&& $w_{j_k,i_k^*}(1+\frac{\lambda_{\ell(i_k^*)}}{\mu_{i_k^*}})$ &&& \\
&&& $w_{j_k,i(m_k^1+2)}(1+\frac{\lambda_{\ell(i(m_k^1+2))}}{\mu_{i(m_k^1+2)}}) $ & -1 &&\\
&&&  $\vdots$ &  & $\ddots$ &\\
&&&$w_{j_k,i(m_k)}(1+\frac{\lambda_{\ell(i(m_k))}}{\mu_{i(m_k)}}) $ & &&-1 
\end{array}
\right],
\end{equation}
where $j_k$ represents the $k$th critical pool with 
$i^{*}_{k}$ representing the pattern associated with the only critical pair of critical pool $j_k$,
and all the rows in $\mathcal{M}_k$ are associated with patterns $i(n)$ ($n\in[m_k]$) requiring resource pool $j_k$ (that is, $j_k\in \mathscr{J}_{i(n)}$ for all $n\in[m_k]$), of which the first $m^1_k$ ones are fully activated; that is, $\iota_{i(n)}\in\mathscr{I}_1(o,\bm{\nu},\bm{0})$ for $n\in[m^1_k]$. 
If $|\mathscr{J}_i|=1$ for all $i\in[I]$, then  there exists $m\in\mathbb{N}_0$, $m_{k}\in \mathbb{N}_{+}$, $k\in[m]$, so that with appropriately re-ordered columns, 
$\tilde{\mathcal{M}}$ has the block diagonal form $diag (\mathcal{M}_k:\ k\in[m+1])$.

The last matrix $\mathcal{M}_{m+1}$ is an identity matrix, of which the rows are related to fully activated PS pairs $\iota\in\mathscr{I}_1(o,\bm{\nu},\bm{0})$ that do not require any of the critical resource pools.

We firstly explain the correspondence of rows and columns between $\mathcal{M}_k$ ($k\in[m]$) and $\tilde{\mathcal{M}}$.
As just described, the first $m_k^1$ rows of $\mathcal{M}_k$ are associated with non-critical, activated PS pairs $\iota\in \mathscr{I}_1(o,\bm{\nu},\bm{0})$, the $(m_k^1+1)$th row corresponds to the critical pair of critical pool $j_k$, and the remaining rows correspond to passive pairs.
From the definition of the matrix $\mathcal{M}$ described in \ref{case1} and \ref{case2}, in the special case with $|\mathscr{J}_i|=1$ for all $i\in[I]$, there are in total two non-zero entries in the row for a non-critical, activated PS pair $\iota\in \mathscr{I}_1(o,\bm{\nu},\bm{0})$: one equal to one and located in the main diagonal of both $\mathcal{M}$ and $\tilde{\mathcal{M}}$; and the other in the column associated with the only resource pool required by this PS pair.
Also, from the definition described in \ref{case2} and \ref{case4}, the rows for passive PS pairs $\iota\in\mathscr{I}_0(o,\bm{\nu},\bm{0})$ have the same structure except the diagonal entries are equal to $-1$.
By re-ordering the rows and columns of $\tilde{\mathcal{M}}$, we can position the rows related to critical pool $j_k$ next to each other,
so that the first $m^1_k$ rows are related to non-critical, activated PS pairs, the $(m^1_k+1)$th row is associated with the only critical PS pair of this critical pool, and the remaining rows stand for the passive ones.

Recall that the rows of the last matrix $\mathcal{M}_{m+1}$ are related to fully activated PS pairs $\iota\in\mathscr{I}_1(o,\bm{\nu},\bm{0})$ that do not require any of the critical resource pools.
Similar to the rows of non-critical, activated PS pairs in the matrix $\mathcal{M}_k$ for $k\in[m]$, the diagonal entries for these rows are always one. Since these PS pairs require no critical resource pool, all the other entries in these rows are zero.
Accordingly, the matrix $\mathcal{M}_{m+1}$ can be constructed as an identity matrix of size $m_{k+1} \coloneqq K-\sum_{k\in[m]}m_k$, which is possibly zero. 


  

For the more general case with $|\mathscr{J}_i|>1$ for some $i\in[I]$,
in a similar way, $\tilde{\mathcal{M}}$ is a block upper triangular matrix with $\mathcal{M}_{k}$ ($k\in[m]$) as just defined as its diagonal.
When the system is weakly coupled, there is at most one critical pool $j_{k}$ with $j_k=j_{\iota_i}(o,\bm{\nu},\bm{0})$ for pattern $i$,
which is associated with a row of $\mathcal{M}_{k}$ ($k\in[m]$) corresponding to row $\bar{i}$ of matrix $\tilde{\mathcal{M}}$, and at most one other critical resource $j_{k'}$, $k\neq k'$, satisfying $j_{k'}\in\mathscr{J}_i$.

If there are such $i,k,\bar{i}$ and $k'$, we perform an elementary row operation on $[\tilde{\mathcal{M}}| \tilde{\bm{y}}]$:
replacing row $\bar{i}$ by the difference between row $\bar{i}$ and row $\bar{i}'$ where row $\bar{i}'$ of $\tilde{\mathcal{M}}$ corresponds to the only critical pair associated with critical pool $j_{k'}$.
The corresponding $\gamma_{j_{k}}$ for row $\bar{i}$ (associated with pattern~$i$) is, then, 
\begin{equation}\label{eqn:gamma_j:1}
\gamma_{j_{k}}=\frac{w_{j_{k'},i}}{w_{j_k,i}}\Bigl(\frac{y_{\sigma(\iota_i)}}{w_{j_{k'},i}(1+\lambda_{\ell_{\iota_i}}/\mu_{i})}
-\frac{y_{\sigma(\iota_{i'})}}{w_{j_{k'},i'}(1+\lambda_{\ell_{\iota_{i'}}}/\mu_{i'})}\Bigr),
\end{equation}
where $i'$ is the pattern corresponding to row $\bar{i}'$.
Note that, if the system is weakly coupled, row $j_{k'}$ is the only type-2 row for columns $i$ and $i'$ (of matrix $\mathcal{W}$) with $w_{j_{k'},i}>0$ and $w_{j_{k'},i'}>0$; that is, resource pool $j_{k'}$ is the only resource pool for patterns $i$ and $i'$ that is shared with multiple patterns. 
Let $\mathscr{M}$ represent the set of subscripts $k\in[m]$ of all critical pools $j_k$ of critical pairs $\iota_{i_k^*}$, for which there exists another critical pool $j_{k'}$ satisfying $j_k,j_{k'}\in\mathscr{J}_{i_k^*}$.

After this operation, we remove row and column $\bar{i}$ from $\tilde{\mathcal{M}}$ and the $\bar{i}$th elements from $\tilde{\bm{x}}$ and $\tilde{\bm{y}}$ correspondingly.
With some abuse of notation, we refer to the resulting matrix and vectors as 
$\tilde{\mathcal{M}}$, $\tilde{\bm{x}}$ and $\tilde{\bm{y}}$.
As a consequence, when all such rows and columns are removed from $\tilde{\mathcal{M}}$, we write the remaining matrix as $\tilde{\mathcal{M}}=diag(\mathcal{M}_{k}:\ k\in[m+1]\backslash\mathscr{M})$. 
Now, for each $\mathcal{M}_k$ ($k\in[m]\backslash\mathscr{M}$) there are $m_k^1$ rows associated with PS pairs $\iota\in\mathscr{I}_1(o,\bm{\nu},\bm{0})$, one row associated with the critical pair for critical pool $j_k$ and $m_k^0\coloneqq m_k-m_k^1-1$ rows for PS pairs $\iota\in\mathscr{I}_0(o,\bm{\nu},\bm{0})$ that are passive because of capacity constraints. 

We obtain the expression 
\begin{equation}
|\tilde{\mathcal{M}}| = \prod\limits_{k\in[m]\backslash\mathscr{M}}|\mathcal{M}_{k}|=\prod\limits_{k\in[m]\backslash\mathscr{M}}\left[w_{j_{k},i^{*}_{k}}\Bigl(\frac{\lambda_{\ell(i^{*}_{k})}}{\mu_{i^{*}_{k}}}+1\Bigr)(-1)^{m_k^0}\right],
\end{equation}
for the determinant of the square matrix $\tilde{\mathcal{M}}$.


Let $\tilde{\mathcal{M}}_{\bar{i}}$, $\bar{i}\in[\sum_{k\in[m]\backslash\mathscr{M}}m_{k}]$, be the matrix after
replacing matrix $\tilde{\mathcal{M}}$'s column $\bar{i}$ by $\tilde{\bm{y}}$.
For a column $\bar{i}$ of matrix $\tilde{\mathcal{M}}$ and pattern $i$ that is associated with row $\bar{i}$ of $\tilde{\mathcal{M}}$, there exists a unique $k\in[m]\backslash\mathscr{M}$, such that
\begin{equation}
\left|\tilde{\mathcal{M}}_{\bar{i}}\right| = 
\prod\limits_{\begin{subarray}{l}k'\in[m]\backslash\mathscr{M}\\k'\neq k\end{subarray}}|\mathcal{M}_{k'}|\times
\begin{cases}
\left|\begin{array}{ll}
y_{\sigma(\iota_i)} & w_{j_k,i}\left(\frac{\lambda_{\ell(i)}}{\mu_{i}}+1\right)\\
y_{\sigma(\iota_{i^*_k})}
& w_{j_k,i^{*}_{k}}\left(\frac{\lambda_{\ell(i^{*}_{k})}}{\mu_{i^{*}_{k}}}+1\right)
\end{array}\right|(-1)^{m_k^0}, & \text{if }\ 
\iota_i\in \mathscr{I}_1(o,\bm{\nu},\bm{0}),
\\ 
\left|\begin{array}{ll}
w_{j_k,i^{*}_{k}}\left(\frac{\lambda_{\ell(i^{*}_{k})}}{\mu_{i^{*}_{k}}}+1\right) & y_{\sigma(\iota_{i^*_k})}\\
w_{j_k,i}\left(\frac{\lambda_{\ell(i)}}{\mu_{i}}+1\right) & y_{\sigma(\iota_i)}
\end{array}\right|(-1)^{m_k^0-1}, & \text{if }
\iota_i\in \mathscr{I}_0(o,\bm{\nu},\bm{0}),
\\ 
y_{\sigma(\iota_i)}(-1)^{m_k^0}, &  \text{otherwise}.
\end{cases}
\end{equation}
By Cramer's Rule, 
we conclude that, for $\iota\in[N-I+L]\backslash \mathscr{N}^0$,
\begin{enumerate}[label=\roman*)]
\item if $\iota \in \mathscr{I}_1(o,\bm{\nu},\bm{0})$ and $\exists k\in[m]$ with $j_k \in\mathscr{J}_{i_{\iota}}$,
\begin{equation}\label{eqn:equal_opt:3}
x_{\sigma(\iota)}= w_{j_k,i_{\iota}}\bigl(1+\frac{\ell(i_{\iota})}{\mu_{i_{\iota}}}\bigr)\Biggl(\frac{y_{\sigma(\iota)}}{w_{j_k,i_{\iota}}\bigl(1+\frac{\lambda_{\ell(i_{\iota})}}{\mu_{i_{\iota}}}\bigr)}
-
\frac{y_{\sigma(\iota_{i^*_k})}}{w_{j_k,i^*_k}\bigl(1+\frac{\lambda_{\ell(i^*_k)}}{\mu_{i^*_k}}\bigr)}\Biggr);
\end{equation}
\item if $\iota \in \mathscr{I}_1(o,\bm{\nu},\bm{0})$ and $\nexists k\in[m]$ with $j_k \in\mathscr{J}_{i_{\iota}}$, 
\begin{equation}\label{eqn:equal_opt:4}
x_{\sigma(\iota)} =y_{\sigma(\iota)};
\end{equation}
\item if $\iota \in \mathscr{I}(o,\bm{\nu},\bm{0})$ and $\exists k,k'\in[m]$ with $k< k'$ and $j_k,j_{k'}\in \mathscr{J}_{i_{\iota}}$,
\begin{equation}\label{eqn:equal_opt:5}
x_{\sigma(\iota)} =\frac{w_{j_{k'},i}}{w_{j_k,i}}\Bigl(\frac{y_{\sigma(\iota_i)}}{w_{j_{k'},i}\bigl(1+\frac{\lambda_{\ell_{\iota_i}}}{\mu_{i}}\bigr)}
-\frac{y_{\sigma(\iota_{i^*_{k'}})}}{w_{j_{k'},i^*_{k'}}\bigl(1+\frac{\lambda_{\ell(i^*_{k'})}}{\mu_{i^*_{k'}}}\bigr)}\Bigr);
\end{equation}
\item if $\iota \in \mathscr{I}_0(o,\bm{\nu},\bm{0})$ and $\exists k\in[m]$ with $j_k \in\mathscr{J}_{i_{\iota}}$,
\begin{equation}\label{eqn:equal_opt:6}
x_{\sigma(\iota)} = w_{j_k,i_{\iota}}\bigl(1+\frac{\ell(i_{\iota})}{\mu_{i_{\iota}}}\bigr)\Biggl(
\frac{y_{\sigma(\iota_{i^*_k})}}{w_{j_k,i^*_k}\bigl(1+\frac{\lambda_{\ell(i^*_k)}}{\mu_{i^*_k}}\bigr)}
-\frac{y_{\sigma(\iota)}}{w_{j_k,i_{\iota}}\bigl(1+\frac{\lambda_{\ell(i_{\iota})}}{\mu_{i_{\iota}}}\bigr)}
\Biggr);
\end{equation}
\item otherwise,
\begin{equation}\label{eqn:equal_opt:7}
x_{\sigma(\iota)} = \frac{y_{\sigma(\iota)}}{w_{j_k,i_{\iota}}\bigl(1+\frac{\lambda_{\ell(i_{\iota})}}{\mu_{i_{\iota}}}\bigr)}.
\end{equation}
\end{enumerate}

Recall that $\sigma(\iota)$ represents the position of PS pair $\iota \in[N-I+L]\backslash\mathscr{N}^0$ among all pairs in the set $[N-I+L]\backslash\mathscr{N}^0$ according to the ranking $o\in\mathscr{O}(\bm{\nu},\bm{0},\bm{\beta})$, and that the PS-pair ranking $o$ is established in the descending order of $\beta_{i_{\iota}}(\Xi_{i_{\iota}}(\bm{0},\bm{0})-\nu_{\ell_{\iota}})$ ($\iota\in[N]$).
From \eqref{eqn:equal_opt:3}-\eqref{eqn:equal_opt:7}, $\bm{x}$ has a non-negative solution if 
\begin{equation}\label{eqn:choice_of_beta}
\beta_i = \frac{1}{w_{j,i}(1+\frac{\lambda_{\ell(i)}}{\mu_i})},
\end{equation}
where $j$ is the only resource pool in $\mathscr{J}_i$ that is shared with multiple patterns, if there is one; or any element of the set $\arg\min_{j'\in\mathscr{J}_i}C_{j'}/w_{j',i}$, otherwise.
In this case, $o$ can be any ranking in $\mathscr{O}(\bm{\nu},\bm{0},\bm{\beta})$.


The resulting values of $\pmb{\gamma}$, as defined in \ref{case2}, are given by \eqref{eqn:equal_opt:5} and \eqref{eqn:equal_opt:7}.
The lemma is then proved.

\proof{Proof of Proposition~\ref{prop:equal_opt}.}

The proposition can be derived from Lemma~\ref{lemma:opt:for_any_nu} by fine-tuning  values of $\bm{\nu}\in\mathbb{R}^{L}$, such that the relaxed action constraints described in \eqref{eqn:constraint:relax:action} are satisfied by the policy $\varphi(o,\bm{\nu},\bm{0})$.

\endproof

\section{Two Examples of Activating Sub-processes}\label{app:example3}
In their asymptotic optimality proof, \cite{weber1990index} allowed $h$ sub-processes simultaneously to be active in the RMABP with scaling parameter $h$. The specific birth-and-death process form of our bandit process allows us to do something different. 
We activate exactly one bandit process for each request type and 
accelerate the birth rate of active bandit processes by a factor of $h$ with
the birth rates of the passive ones remaining zero.

That is, as described in \eqref{eqn:constraint:action:h} in Section~\ref{subsec:asym_regime}, we activate exactly one sub-process $(i,k)$ ($i\in\mathscr{P}_{\ell}$, $k\in[h]$) for RT $\ell\in[L]$ regardless of the scaling parameter $h\in\mathbb{N}_+$.
The birth and death rates of this active sub-process are $h\lambda_{\ell}^0$ and $N^{\phi}_{i,k}(t)\mu_{i}^0$, respectively.


\begin{figure}[t]
\centering
\subfigure[]{\includegraphics[width=0.2\linewidth]{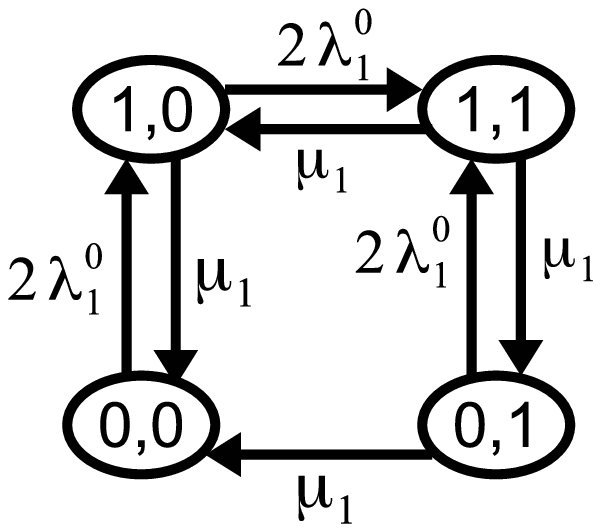}\label{fig:fixed_vs_proportional:fixed}}
\subfigure[]{\includegraphics[width=0.19\linewidth]{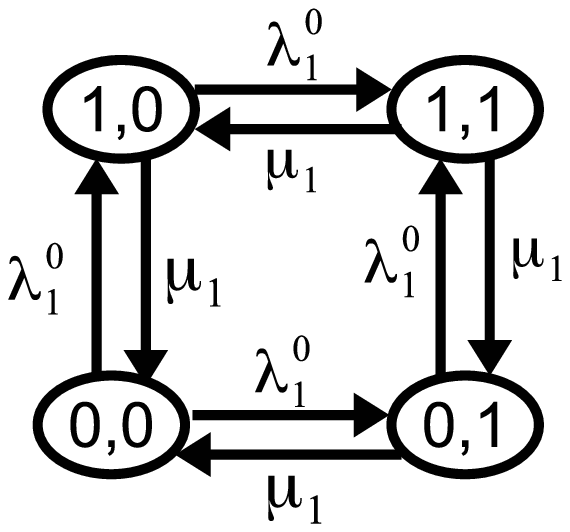}\label{fig:fixed_vs_proportional:proportional}}
\caption{Markov chains for cases with proportional and fixed numbers of active sub-processes: (a) our model; and (b) the model of \cite{weber1990index}.}\label{fig:fixed_vs_proportional}
\end{figure}


For instance, consider a simple system with only one type of request and one non-dummy pattern available for it so that $|\mathscr{P}_1|=2$. We label this pattern as pattern $1$, of which the state space is $\mathscr{N}^0_1=\{0,1\}$ with $w_{j,1}=1$ for all $j\in[J]$.
For a trivial policy that always prioritizes the only non-dummy pattern or its sub-patterns if the capacity constraints are not violated, we illustrate in Figure~\ref{fig:fixed_vs_proportional} the underlying Markov chains for our model and that of \cite{weber1990index} when $h=2$. 
For the latter, the two sub-patterns of the non-dummy pattern will be activated simultaneously all the time 
with the same Markov chain as shown in Figure~\ref{fig:fixed_vs_proportional}(b). 
For the former, its Markov chain is illustrated in Figure~\ref{fig:fixed_vs_proportional}(a),
where the first sub-pattern is activated with birth rate $2\lambda_{1}^0$, while the second sub-pattern is passive with zero birth rate. 
This has consequences for our proof of asymptotic optimality  and further discussion about the effects of fixed and proportional numbers of active sub-patterns with $h\rightarrow +\infty$ is provided in conjunction with Proposition~\ref{prop:equilibrium_vector}.

\section{Proof of Theorem~\ref{theorem:main_second}}\label{app:theorem:main_second}

Here, we first show that, for a given PS pair ranking $o\in\mathscr{O}$, 
the index policy $\varphi$ (described by \eqref{eqn:tilde_z} and \eqref{eqn:probability_active}), based on ranking $o$, achieves the same long-run average revenue as $\bar{\varphi}(o)$ (generated by Algorithm~\ref{algo:varphi_gamma}).

Let $q^{1,h}(\iota,\iota')$ and $q^{0,h}(\iota,\iota')$, $\iota,\iota'\in[N]$, represent the transition rates of a sub-process in PS pair $\iota$ transitioning to $\iota'$ with $i_{\iota'}=i_{\iota}$ and $n_{\iota'} = n_{\iota}\pm 1$, if it is active and passive, respectively.
For notational convenience, 
define $q^{1,h}(\iota,\iota')=q^{0,h}(\iota,\iota') = 0$ if $i_{\iota'}\neq i_{\iota}$ or $n_{\iota'}\neq n_{\iota}\pm 1$, $\iota,\iota'\in[N]$.
Consider a deterministic process $\bm{z}^{\varphi,h}(t)=(z^{\phi,h}_{\iota}(t):\ \iota\in[N])$ with a given initial point $\bm{z}^{\varphi,h}(0)\in\mathscr{Z}$ under the index policy $\varphi$ (described in Section~\ref{subsec:policies}),
 generated by the differential equation
\begin{multline}\label{eqn:differential}
\frac{d\ z_{\iota}^{\varphi,h}(t)}{d\ t}=\sum\limits_{
\begin{subarray}\ \iota'\in[N]\\ \iota'\neq \iota\end{subarray}}\biggl[z_{\iota'}^{\varphi,h}(t)\Bigl(\upsilon^{\varphi,h}_{\iota'}(\bm{z}^{\varphi,t}(t))q^{1,h}(\iota',\iota)+\bigl(1-\upsilon^{\varphi,h}_{\iota'}(\bm{z}^{\varphi,h}(t))\bigr)q^{0,h}(\iota',\iota)\Bigr)\\
-
z_{\iota}^{\varphi,h}(t)\Bigl(\upsilon^{\varphi,h}_{\iota}(\bm{z}^{\varphi,t}(t))q^{1,h}(\iota,\iota')+\bigl(1-\upsilon^{\varphi,h}_{\iota}(\bm{z}^{\varphi,h}(t))\bigr)q^{0,h}(\iota,\iota')\Bigr)\biggr],
\end{multline}
where the right hand side is the sum of transition rates entering PS pair $\iota$ minus the sum of transition rates leaving PS pair $\iota$.

\begin{proposition}\label{prop:convergence}
\phantom{==}

\begin{enumerate}
\item
For any $h\in\mathbb{N}_+$, $\bm{\epsilon}^h\in \mathscr{E}^h$ and $\bm{z}^{\varphi,h}(0)\in\mathscr{Z}$, there exists a unique solution $\bm{z}^{\varphi,h}(t)$, $t \geq 0$, of \eqref{eqn:differential}.
\item For any $\psi\in \Psi$ and $\delta > 0$, if $\bm{Z}^{\varphi,h}(0)=\bm{z}^{\varphi,h}(0)\in\mathscr{Z}$,
\begin{equation}\label{eqn:convergence}
\lim\limits_{h\rightarrow +\infty}\lim\limits_{t\rightarrow +\infty} \frac{1}{t}\int_{0}^{t}\mathbb{P}\Bigl\{\norm{\bm{Z}^{\varphi,h}(u)-\bm{z}^{\varphi,h}(u)}>\delta\Bigr\}du = 0,
\end{equation}
where $\norm{\cdot}$ is the Euclidean norm.
\end{enumerate}
\end{proposition}
Recall that the index policy $\varphi$, described by \eqref{eqn:tilde_z} and \eqref{eqn:probability_active}, is $\bm{\epsilon}^h$-dependent. 

The proof is given in Appendix~\ref{app:prop:convergence}.
The proof follows similar ideas and methods to those of  \cite{weber1990index} by invoking Picard's existence theorem and \cite[Chapter 7, Theorem 2.1]{freidlin2012random}.

We define a special point
\begin{equation}\label{eqn:equilibrium_vector}
\bm{z}^{\bar{\varphi}(o),h}= \lim\limits_{t\rightarrow +\infty}\mathbb{E}\left[\bm{Z}^{\bar{\varphi}(o),h}(t)\right],
\end{equation}
which exists since,  with given $h\in\mathbb{N}_+$, the underlying Markov chain of the process $\bm{N}^{\bar{\varphi}(o)}_h(t)$ is irreducible with finitely many states and the process $\bm{Z}^{\bar{\varphi}(o),h}(t)$ is a function of the process $\bm{N}^{\bar{\varphi}(o)}_h(t)$ 
 as defined in Section~\ref{subsec:asym_regime}.

The key to our proof of Theorem~\ref{theorem:main_second} is the following result.
\begin{lemma}\label{lemma:zero_derivative}
The limit $ \bm{z}^{\bar{\varphi}(o)}=\lim_{h\rightarrow +\infty} \bm{z}^{\bar{\varphi}(o),h}$ exists, and for any $\bm{z}^{\varphi,h}(0)\in\mathscr{Z}$,
\begin{equation}\label{eqn:zero_derivative}
\lim\limits_{\norm{ \bm{\epsilon} }\to 0}\lim\limits_{h\rightarrow +\infty}\lim\limits_{t\rightarrow+\infty} \bm{z}^{\varphi,h}(t) =  \bm{z}^{\bar{\varphi}(o)}.
\end{equation}
\end{lemma}
The proof is given in Appendix~\ref{app:prop:equilibrium_vector}.

{\bf Remark}
A condition similar to Lemma~\ref{lemma:zero_derivative}, trapping the process $\bm{z}^{\varphi,h}(t)$ in the neighborhood of $\bm{z}^{\bar{\varphi}(o)}$, was an important assumption in \cite{weber1990index}: there exists a unique equilibrium point of the constructed deterministic process with any initial point.
In our problem, the underlying Markov chain for each sub-process is a birth-and-death process with state-dependent service rates that are monotonically non-increasing in their priorities. This simplifies the analysis of $\bm{z}^{\varphi,h}(t)$ in the asymptotic regime.  

Intuitively, under an index policy, the process related to the pattern with the highest priority will keep occupying or releasing resource units until it reaches a stable point: a capacity constraint is about to be violated or the arrival rate of the related RT is balanced by the total service rate. Then, another process with the second priority will behave similarly: it cannot occupy any resource units already occupied by the first process but
can take resource units from other processes with lower priorities. The final equilibrium point can be obtained by calculating these stable points in turn, moving from the pattern with highest priority to the lowest one. 
A detailed analysis is provided in the proof of Lemma~\ref{lemma:zero_derivative} in Appendix~\ref{app:prop:equilibrium_vector}.


\begin{proposition}\label{prop:equilibrium_vector}
For any $\sigma >0$,
\begin{equation}\label{eqn:convergence:2}
\lim\limits_{\norm{ \bm{\epsilon} }\to 0}\lim\limits_{h\rightarrow +\infty}\lim\limits_{t\rightarrow+\infty} \frac{1}{t}\int_0^t \mathbb{P}\Bigl\{\norm{\bm{Z}^{\varphi,h}(u)-\bm{z}^{\bar{\varphi}(o)}} > \sigma\Bigr\} du =0
\end{equation}
\end{proposition}
\proof{Proof.}
Proposition~\ref{prop:convergence} shows that for any $\bm{\epsilon}$, the stochastic process $\bm{Z}^{\varphi,h}(t)$ will not leave a neighbourhood of $\bm{z}^{\varphi,h}(t)$ for any substantially positive proportion of time when $h\rightarrow +\infty$.
Lemma~\ref{lemma:zero_derivative} shows that the deterministic process $\bm{z}^{\varphi,h}(t)$ will not leave a neighbourhood of $\bm{z}^{\bar{\varphi}(o)}$ for any substantially positive proportion of time when $\norm{ \bm{\epsilon} }\to 0$ and $h\rightarrow +\infty$.

This proposition is a consequence of Proposition~\ref{prop:convergence} and Lemma~\ref{lemma:zero_derivative}.

\endproof

{\bf Remark}
Recall the long-run average revenue normalized by $h$ of the resource allocation problem under policy $\phi$, $R^{\phi,h}$, as defined in Section~\ref{sec:asym_opt}.
For given $o\in\mathscr{O}$, 
if policy $\bar{\varphi}(o)\in\tilde{\Phi}$ maximizes the long-run average revenue of the relaxed problem for a given $h$, then $R^{\bar{\varphi}(o),h} \geq \max_{\phi\in\Phi}R^{\phi,h} \geq R^{\varphi,h}$.
From Lemma~\ref{lemma:zero_derivative} and Proposition~\ref{prop:equilibrium_vector}, since $R^{\varphi,h}$ approaches $R^{\bar{\varphi}(o),h}$ as $\norm{ \bm{\epsilon} }\to 0$ and $h\rightarrow +\infty$, $R^{\varphi,h}$ asymptotically approaches optimality of the original problem.


\proof{Proof of Theorem~\ref{theorem:main_second}.}
If \eqref{eqn:main_second} holds for $\bar{\varphi}(o)$, then, from Proposition~\ref{prop:equilibrium_vector}, $R^{\varphi,h}$ approaches $\max_{\phi\in\Phi}R^{\phi,h}$ as $\lVert \bm{\epsilon} \rVert\to 0$ and $h\rightarrow+\infty$. 
If $R^{\varphi,h}$ asymptotically approaches $\max_{\phi\in\Phi}R^{\phi,h}$ as $\lVert \bm{\epsilon} \rVert\to 0$ and $h\rightarrow+\infty$, then, from Proposition~\ref{prop:equilibrium_vector}, $R^{\bar{\varphi}(o),h}$ also approaches $\max_{\phi\in\Phi}R^{\phi,h}$
as $h\rightarrow+\infty$; that is, \eqref{eqn:main_second} is satisfied.

\endproof

\section{Proof of Proposition~\ref{prop:convergence}.}\label{app:prop:convergence}

\proof{Proof of Proposition~\ref{prop:convergence}.}
This proof follows the ideas and methods of \cite{weber1990index} and \cite{freidlin2012random}.

We firstly construct a stochastic process that matches the hypothesis of \cite[Chapter 7, Theorem 2.1]{freidlin2012random}.

Let $t^{h}_{\ell,s}$ be the times of the $s$th arrival of a request of type $\ell\in[L]$ in a system with scaling parameter $h$.
For convenience, define a one-to-one mapping $\theta$: $\Bigl([I]\backslash \{d(\ell):\ell\in[L]\}\Bigr)\times [h] \to[h(I-L)]$, so that each non-dummy sub-pattern or sub-process $(i,k)$, $i\in[I]\backslash \{d(\ell):\ell\in[L]\}$, $k\in[h]$, can be labeled by an integer $\theta = \theta(i,k) \in [h(I-L)]$.
Let $(i^{\theta},k^{\theta})$ represent the sub-pattern or sub-process labeled by $\theta\in[h(I-L)]$, and $t^{h}_{L+\theta,s}(n)$ be the times of the $s$th \emph{potential departure} of a completed request with sub-process
$(i^{\theta}, k^{\theta})$ in state $n\in\mathscr{N}_{i^{\theta}}^0$, $s\in\mathbb{N}_+$.
In particular, let the $h$ sub-processes for the same pattern be labeled successively: the mapping $\theta$ satisfies that, for any $i\in[I]\backslash\{d(\ell):\ell\in[L]\}$ and $k\in [h-1]$,  $\theta(i,k+1) = \theta(i,k)+1$.

The time intervals between successive arrivals and potential departures of requests are independently and exponentially distributed random variables.
The word \emph{potential} is used because, for each pattern/sub-pattern, a potential departure is regarded as a real departure if there is at least one request being served by this pattern; and it is ignored, otherwise.


Let $\tau(t)$ represent the occurrence time of the latest event, either arrival or departure, before time $t$.
Define a random vector $\bm{\xi}^{h}_{t}=(\xi^{h}_{\ell,t}: \ell\in[L],\ \xi^{h}_{L+\theta,t}(n):\theta\in[h(I-L)],\ n\in\mathscr{N}^0_{i^{\theta}})$ as follows:
for $\ell\in[L]$,
\begin{equation}\label{eqn:app:convergence:1}
\xi^{h}_{\ell,t}=
\begin{cases}
\frac{1}{t^{h}_{\ell,s^*}-\tau(t^h_{\ell,s^*})}, & \text{if }
\tau(t^{h}_{\ell,s^*}) \leq t < t^{h}_{\ell,s^*} \text{, where } 
t^{h}_{\ell,s^*}=\min\limits_{\begin{subarray}~k''\in\mathbb{N}_+\end{subarray}}\{t^{h}_{\ell,k''}|t^{h}_{\ell,k''}> t\},
\\
0, & \text{otherwise},
\end{cases}
\end{equation}
for $\theta\in[h(I-L)]$, $n\in\mathscr{N}^0_{i^{\theta}}$,
\begin{equation}\label{eqn:app:convergence:2}
\xi^{h}_{L+\theta,t} (n)=
\begin{cases}
\frac{1}{t^{h}_{\theta,s^*}(n)-\tau(t^{h}_{\theta,s^*}(n))}, & \text{if }
\tau(t^{h}_{\theta,s^*}(n)) \leq t < t^{h}_{\theta,s^*}(n) \text{, where }
 t^{h}_{\theta,s^*}(n)=\min\limits_{\begin{subarray}~k''\in\mathbb{N}_+\end{subarray}}\{t^{h}_{\theta,k''}(n)|t^{h}_{\theta,k''}(n)> t\},
\\
0, & \text{otherwise}.
\end{cases}
\end{equation}

From \eqref{eqn:app:convergence:1}, at any time $t>0$, there is at most one $\ell \in [L]$ for which $\xi^h_{\ell,t} > 0$ and all the others are zero. Also, for $\ell\in[L]$, if the next event is not an arrival of request type $\ell$ (but a potential departure or an arrival of other request types), $\xi^h_{\ell,t}=0$. 
If the next event to occur after time $t$ is the $s$th arrival of type $\ell$, then  $\xi^h_{\ell,t}$ is the reciprocal of the time between the last and the next event. Otherwise it is zero.
Similar results apply to $\xi^h_{L+\theta}(n)$ ($\theta\in[h(I-L)]$, $n\in\mathscr{N}^0_{i^{\theta}}$) defined in \eqref{eqn:app:convergence:2} associated with potential departures of sub-process $(i^{\theta},k^{\theta})$ in state $n$.

Then, for $\ell\in[L]$, the number $\lfloor \int_0^T \xi^{h}_{\ell,t} dt\rfloor$ represents the number of arrivals of RT $\ell$ by time $T$; and for $\theta\in[h(I-L)]$ and $n\in\mathscr{N}^0_{i^{\theta}}$, the number $\lfloor \int_0^T \xi^{h}_{L+\theta,t}(n)  dt \rfloor$ is the number of potential departures associated with sub-pattern
$(i^{\theta},k^{\theta})$ in its state $n$ by time $T$, when $T$ is large.


We define a function, $Q^{h}(\iota,\iota',\bm{x},\bm{\xi}^{h})$, on $\iota,\iota'\in[N]$, $\bm{x}\in\mathbb{R}^{N}$, $\bm{\xi}^{h}\in\mathbb{R}^{L+h(N-L)}$.
For $\bm{z}\in\mathscr{Z}$, $\iota\in[N]$, let 
$a_{\iota}(\bm{z}) = \upsilon^{\varphi,h}_{\iota}(\bm{z})z_{\iota}hI$, which, from the definition of $\upsilon^{\varphi,h}_{\iota}(\bm{z})$ in \eqref{eqn:probability_active}, takes values in $\{0,1\}$ and becomes algebraically independent of $h$ with given $\bm{z}$, when $h$ is sufficiently large. 
For given $\bm{\xi}^{h}\in \mathbb{R}^{L+h(N-L)}$ and $\bm{x}\in \mathbb{R}^{N}_{0}$, 
$Q^{h}(\iota,\iota',\bm{x},\bm{\xi}^{h})$ is defined by
\begin{equation}\label{eqn:definition_Q}
Q^{h}(\iota,\iota',\bm{x},\bm{\xi}^{h})=
\begin{cases}
a_{\iota}(\bm{x}/hI)\ 
\xi^{h}_{\ell_{\iota}}+f^{0,h}_{\iota,a}(\bm{x},\bm{\xi}^{h}),&\text{if }i_{\iota}=i_{\iota'},\ n_{\iota'}=n_{\iota}+1,\\
\sum\limits_{\theta=\theta(i_{\iota},\lceil x^{-}_{\iota-1}\rceil+1)}^{\theta(i_{\iota},\lceil x^{-}_{\iota}\rceil)}\xi^{h}_{L+\theta}(n_{\iota})
 + f^{h}_{\iota,a}(\bm{x},\bm{\xi}^{h}),&\text{if }i_{\iota}=i_{\iota'},\ n_{\iota'}=n_{\iota}-1,\\
0, &\text{otherwise},
\end{cases}
\end{equation}
where $x^{-}_{\iota}=\sum_{\iota'=1,i_{\iota'}=i_{\iota}}^{\iota}x_{\iota'}$ with $x^{-}_{0}=0$, and 
$f^{0,h}_{\iota,a}(\bm{z},\bm{\xi}^h)$ and $f^{h}_{\iota,a}(\bm{z},\bm{\xi}^h)$ are appropriate functions to make $Q^{h}(\iota,\iota',\bm{x},\bm{\xi}^{h})$ Lipschitz continuous in $\bm{x}$ for all $\bm{\xi}^{h}$ and $0<a<1$, see the discussion below.  For $\bm{x}\in \mathbb{R}^{N}\backslash\mathbb{R}^{N}_0$, let $(\bm{x})^+=(\max\{x_{\iota},0\}:\ \iota\in[N])$, and define $Q^{h}(\iota,\iota',\bm{x},\bm{\xi}^{h})=Q^{h}(\iota,\iota',(\bm{x})^+,\bm{\xi}^{h})$.

If we set $f^{0,h}_{\iota,a}(\bm{z},\bm{\xi}^h)$ and $f^{h}_{\iota,a}(\bm{z},\bm{\xi}^h)$ to zero all the time, then $Q^{h}(\cdot,\cdot,\bm{x},\cdot)$ is a discontinuous function of $\bm{x}$.
The idea to make the function smooth in $\bm{x}$ is to utilize the fact that the Dirac delta function can be considered as a limit 
of normal distributions of mean $0$.

With $0<a<1$, define
\begin{equation}
y_a(u)=\begin{cases}
\int_{-\infty}^{\varrho(u)} \frac{1}{a\sqrt{\pi}}e^{-(v-\frac{1}{a})^{2}/a^{2}}d  v & \text{if }u\in(0,1)\\
0 & \text{if }u=0\\
1 & \text{if }u=1,
\end{cases}
\end{equation}
where $\varrho(u)$ is a function of $u\in (0,1)$ which is $-\infty$ at $0$, $+\infty$ at $1$ and suitably smooth. 
It can possibly be specified as a negative cotangent function, $\varrho(u)=-\cot(u\pi)$.
The function $y_a(u)$ is continuous at $u\in(0,1)$, right-continuous at $u=0$ and left continuous at $u=1$. Let $\tilde{\mathscr{N}}$ be  the set of PS pairs associated with dummy patterns, and for any $\iota\in\tilde{\mathscr{N}}$, define $f^{0,h}_{\iota,\cdot}(\cdot,\cdot)=f^h_{\iota,\cdot}(\cdot,\cdot)=0$.

For $\iota\in[N]$, let 
\begin{equation}
\Gamma_{\iota}(\bm{z}) \coloneqq \lim_{h\rightarrow +\infty}\min\Biggl\{z_{\iota},\max\biggl\{0,\min\limits_{j\in\mathscr{J}_i} \frac{1}{w_{j,i}I}\Bigl( C_{j}^{0}(1-\epsilon^h_{j,\iota})-\sum\limits_{\iota'=1}^{N}w_{j,i_{\iota'}}n_{\iota'}z_{\iota'}I-\sum\limits_{\iota'\in \mathscr{N}_{\iota}^{+}}w_{j,i_{\iota'}}\upsilon^{\varphi,h}_{\iota'}(\bm{z})z_{\iota'}I\Bigr)\biggr\}\Biggr\},
\end{equation}
where, from the definition of $\upsilon^{\varphi,h}_{\iota}(\bm{z})$ in \eqref{eqn:probability_active}, the limit exists.
Note that $\Gamma_{\iota}(\bm{z})$ is used for definitions of $f^{0,h}_{\iota,a}(\bm{x},\bm{\xi}^h)$ and $f^{h}_{\iota,a}(\bm{x},\bm{\xi}^h)$ with $\iota\in[N]\backslash \tilde{\mathscr{N}}$, and is similar to $\lim\limits_{h\rightarrow +\infty}\zeta^{\varphi,h}_{\iota}(\bm{z})$ defined in \eqref{eqn:tilde_z}; $\Gamma_{\iota}(\bm{z})$ is Lipschitz continuous in $\bm{z}\in\mathscr{Z}$ and the latter is not.
Also, let $\chi^{h}_{\iota}(\bm{x})\coloneqq \sum_{\iota''\in\mathscr{N}^{+}_{\iota},\ \ell_{\iota''}=\ell_{\iota}}\Gamma_{\iota''}(\bm{x}/hI)hI$,
so that $\chi^{h}_{\iota}(\bm{x})$ is continuous in $\bm{x}\in \mathbb{R}^{N}_0$.
For $\bm{x}\in\mathbb{R}^N_0$, $\iota\in[N]\backslash \tilde{\mathscr{N}}$ and $h\in\mathbb{N}_+$, let
$\iota'=\max\{\iota''\in\mathscr{N}^{+}_{\iota}\ |\ \ell_{\iota''}=\ell_{\iota},\ \Gamma_{\iota}(\bm{x}/hI)>0\}$.
We then define
\begin{equation}
f^{0,h}_{\iota,a}(\bm{x},\bm{\xi}^h) = 
\begin{cases}
\xi^{h}_{\ell_{\iota}} y_a(1-\chi^{h}_{\iota'}(\bm{x})), &
\text{if }\ \Gamma_{\iota}(\bm{x}/hI)>0,\ \text{and } 0<\chi^{h}_{\iota'}(\bm{x})<1, \\
-\xi^{h}_{\ell_{\iota}} y_a(1-\chi^{h}_{\iota}(\bm{x})), & \text{if }\ \chi^{h}_{\iota'}(\bm{x})=0,\ \text{and }0<\chi^{h}_{\iota}(\bm{x})<1,\\
0,&\text{otherwise}.
\end{cases}
\end{equation}
For $\bm{x}\in\mathbb{R}^N_0$ and $\iota\in[N]\backslash\tilde{\mathscr{N}}$ with $x_{\iota-1}^{-}>0$ and $\lceil x_{\iota}^{-}\rceil<h$, define
\begin{multline}
f^{h}_{\iota,a}(\bm{x},\bm{\xi}^{h})=
-\xi^h_{L+\theta\bigl(i_{\iota},\lceil x^{-}_{\iota-1}\rceil+1\bigr)}(n_{\iota})y_a(1-\lceil x^{-}_{\iota-1}\rceil+x^{-}_{\iota-1})
+{\xi^h_{L+\theta\bigl(i_{\iota},\lceil x_{\iota-1}^{-}\rceil\bigr)}(n_{\iota})}y_a(\lceil x_{\iota-1}^{-}\rceil - x_{\iota-1}^{-})\\
-{\xi^h_{L+\theta\bigl(i_{\iota},\lceil x_{\iota}^{-}\rceil\bigr)}(n_{\iota})}y_a(\lceil x_{\iota}^{-}\rceil - x_{\iota}^{-})
+{\xi^h_{L+\theta\bigl(i_{\iota},\lceil x^{-}_{\iota} \rceil+1\bigr)}(n_{\iota})}y_a(1-\lceil x^{-}_{\iota}\rceil+x^{-}_{\iota}).
\end{multline}
For $\bm{x}\in\mathbb{R}^N_0$ and $\iota\in[N]\backslash\tilde{\mathscr{N}}$ with $x_{\iota-1}^{-}=0$ and $\lceil x_{\iota}^{-}\rceil<h$, define
\begin{multline}
f^{h}_{\iota,a}(\bm{x},\bm{\xi}^{h})=-{\xi^h_{L+\theta\bigl(i_{\iota},\lceil x^{-}_{\iota-1}\rceil+1\bigr)}(n_{\iota})}y_a(1-\lceil x^{-}_{\iota-1}\rceil+x^{-}_{\iota-1}) \\
-{\xi^h_{L+\theta\bigl(i_{\iota},\lceil x_{\iota}^{-}\rceil\bigr)}(n_{\iota})}y_a(\lceil x_{\iota}^{-}\rceil - x_{\iota}^{-})
+{\xi^h_{L+\theta\bigl(i_{\iota},\lceil x^{-}_{\iota} \rceil+1\bigr)}(n_{\iota})}y_a(1-\lceil x^{-}_{\iota}\rceil+x^{-}_{\iota}),
\end{multline}
For $\bm{x}\in\mathbb{R}^N_0$ and $\iota\in[N]\backslash\tilde{\mathscr{N}}$ with $x_{\iota-1}^{-}>0$ and $\lceil x_{\iota}^{-}\rceil=h$, define
\begin{multline}
f^{h}_{\iota,a}(\bm{x},\bm{\xi}^{h})=
-{\xi^h_{L+\theta\bigl(i_{\iota},\lceil x^{-}_{\iota-1}\rceil+1\bigr)}(n_{\iota})}y_a(1-\lceil x^{-}_{\iota-1}\rceil+x^{-}_{\iota-1}) \\
+{\xi^h_{L+\theta\bigl(i_{\iota},\lceil x_{l-1}^{-}\rceil\bigr)}(n_{\iota})}y_a(\lceil x_{l-1}^{-}\rceil - x_{l-1}^{-})
-\xi^h_{L+\theta\bigl(i_{\iota},\lceil x_{\iota}^{-}\rceil\bigr)}(n_{\iota})y_a(\lceil x_{\iota}^{-}\rceil - x_{\iota}^{-})
\end{multline}
For $\bm{x}\in\mathbb{R}^N_0$ and $\iota\in[N]\backslash\tilde{\mathscr{N}}$ with $x_{\iota-1}^{-}=0$ and $\lceil x_{\iota}^{-}\rceil=h$, define
\begin{equation}
f^{h}_{\iota,a}(\bm{x},\bm{\xi}^{h})=-{\xi^h_{L+\lceil x^{-}_{\iota-1}\rceil+1}(n_{\iota})}y_a(1-\lceil x^{-}_{\iota-1}\rceil+x^{-}_{\iota-1}) 
-{\xi^h_{L+\lceil x_{\iota}^{-}\rceil}(n_{\iota})}y_a(\lceil x_{\iota}^{-}\rceil - x_{\iota}^{-}).
\end{equation}
With these $f^{0,h}_{\iota,a}(\bm{z},\bm{\xi}^h)$ and $f^{h}_{\iota,a}(\bm{z},\bm{\xi}^h)$,  the function $Q^{h}(\iota,\iota',\bm{x},\bm{\xi}^{h})$ is Lipschitz continuous in $\bm{x}\in\mathbb{R}^{N}$ for all given $\bm{\xi}^{h} \in \mathbb{R}^{L+h(N-L)}$ and $0<a<1$.


For the special case with $h=1$, any given $0<a<1$, 
and $\sigma>0$,
we define $\bm{X}_{t}^{\sigma}$ such that
\begin{equation}\label{eqn:differential_equation_x}
\dot{\bm{X}}_{t}^{\sigma}\coloneqq
b(\bm{X}_{t}^{\sigma}, \bm{\xi}^{1}_{t/\sigma})\coloneqq\sum\limits_{\iota'\in[N]}Q^{1}(\iota',\iota,\bm{X}^{\sigma}_{t},\bm{\xi}_{t/\sigma}^{1})-Q^{1}(\iota,\iota',\bm{X}^{\sigma}_{t},\bm{\xi}_{t/\sigma}^{1}).
\end{equation}
It follows that $b(\bm{X}^{\sigma}_{t}, \bm{\xi}^{1}_{t/\sigma})$ satisfies a Lipschitz condition over $\bm{X}^{\sigma}_{t}$ and $\bm{\xi}^{1}_{t/\sigma}$. 

Also, from the definition of the function $Q^h$ in \eqref{eqn:definition_Q}, for $\bm{x}\in \mathbb{R}^{N}$, $t\geq 0$ and $h=1$, there exists a matrix $\tilde{\mathcal{Q}}(\bm{x})$, such that
$
b(\bm{x},\bm{\xi}^1_t) = \tilde{\mathcal{Q}}(\bm{x})\bm{\xi}^1_t.
$
For any $\bm{x}\in \mathbb{R}^{N}$, $h=1$, $\delta>0$ and any function $\overline{b}(\bm{x})\in\mathbb{R}^{N}$ of $\bm{x}\in\mathbb{R}^N$, we obtain
\begin{multline}
\mathbb{P}\biggl\{
\Bigl\lVert \frac{1}{T}\int_{t}^{t+T}b(\bm{x},\bm{\xi}^{1}_{v})dv-\overline{b}(\bm{x})\Bigr\rVert>\delta\biggr\}\\
\leq 
\mathbb{P}\biggl\{\Bigl\lVert\tilde{\mathcal{Q}}(\bm{x})\frac{1}{T}\bigl\lfloor\int_{t}^{t+T}\bm{\xi}^1_v dv\bigr\rfloor - \overline{b}(\bm{x})\Bigr\rVert + \frac{1}{T}\Bigl\lVert\tilde{\mathcal{Q}}(\bm{x})\Bigl(\int_{t}^{t+T}\bm{\xi}^1_v dv - \bigl \lfloor\int_{t}^{t+T}\bm{\xi}^1_v dv\bigr\rfloor \Bigr)\Bigr\rVert > \delta\biggr\}\\
\leq 
\mathbb{P}\biggl\{\Bigl\lVert\tilde{\mathcal{Q}}(\bm{x})\frac{1}{T}\bigl\lfloor\int_{t}^{t+T}\bm{\xi}^1_v dv\bigr\rfloor - \overline{b}(\bm{x})\Bigr\rVert + \frac{o(T)}{T} > \delta\biggr\},
\end{multline}
where $\lfloor \bm{\xi} \rfloor$ for a vector $\bm{\xi}$ takes the floor operation for each of its element.
Recall that, for $\ell\in[L]$, the $\ell$th element of the vector $\lfloor \int_t^{t+T} \bm{\xi}^1_v dv\rfloor$ is a Poisson distributed random variable with rate $\lambda_{\ell}$, representing the number of arrivals of requests of type $\ell$; and, for $\theta\in [I-L]$ and $n\in\mathscr{N}^0_{i^{\theta}}$, its $\theta+L$th element is a Poisson distributed random variable with rate $n\mu_{i^{\theta}}$, representing the number of potential departures of requests for a sub-process in state $n$ of pattern $i^{\theta}$. 
Thus, by the Law of Large Numbers, for any $\bm{x}\in \mathbb{R}^{N}$, $h=1$, $\delta>0$, there exists $\overline{b}(\bm{x})=\mathbb{E}[b(\bm{x},\bm{\xi}^{1}_{t})]$, which is independent from $t$, satisfying
\begin{equation}\label{eqn:equilibrium_b}
\lim\limits_{T\rightarrow +\infty} \mathbb{P}\left\{
\Bigl\lVert \frac{1}{T}\int_{t}^{t+T}b(\bm{x},\bm{\xi}^{1}_{v})d\ v-\overline{b}(\bm{x})\Bigr\rVert>\delta\right\} = 0,
\end{equation}
uniformly in $t>0$.

Let $\overline{\bm{x}}(t)$ be the solution of 
\begin{equation}\label{eqn:prop:convergence:differential_b}
\dot{\overline{\bm{x}}}(t)=\overline{b}(\overline{\bm{x}}(t)),
\end{equation}
with given $\overline{\bm{x}}_{0}=\bm{X}^{\sigma}_{0}$.
By Picard's Existence Theorem (\cite{coddington1955theory}), there exists a unique solution $\overline{\bm{x}}_{t}$, $t\geq 0$, satisfying \eqref{eqn:prop:convergence:differential_b} with given $\overline{\bm{x}}_{0}$. 


Now we invoke \cite[Chapter 7, Theorem 2.1]{freidlin2012random}: if (\ref{eqn:equilibrium_b}) holds true,
and $\mathbb{E} \norm{b(\bm{x},\bm{\xi}^{1}_{t})}^2 < +\infty$ for all $\bm{x}\in R^{N}$, 
then,
for any $T>0$, $\delta >0 $,
\begin{equation}\label{eqn:sup_equilibrium_1}
\lim\limits_{\sigma\rightarrow 0}\mathbb{P}\left\{ \sup\limits_{0\leq t\leq T}\norm{\bm{X}^{\sigma}_{t}-\overline{\bm{x}}(t)}>\delta\right\} = 0.
\end{equation}


We interpret the scalar $\sigma$ and the scaling effects in another way.
For $\bm{x}\in R^{N}$ and $\bm{\xi}^{h}\in R^{L+h(N-L)}$, we define 
\begin{equation}\label{eqn:differential_equation_xn}
b^{h}(\bm{x}, \bm{\xi}^{h})\coloneqq\sum\limits_{\iota'\in[N]}Q^{h}(\iota',\iota,\bm{x},\bm{\xi}^{h})-Q^{h}(\iota,\iota',\bm{x},\bm{\xi}^{h}).
\end{equation}
Following the technique of \cite{fu2016asymptotic}, we set $\sigma = 1/h$, and observe that, for any $\bm{x}\in \mathbb{R}^{N}$, $h>0$ and $T>0$, $\int_{0}^{T}b(\bm{x},\xi^{1}_{t/\sigma})dt$ and $\int_{0}^{T}(b^{h}(h\bm{x},\xi^{h}_{t})/h)dt$ are identically distributed.
With $\bm{Z}^{\sigma}_{0}= \bm{Z}^{h}_{0}= \bm{x}_{0}/I$, we define
$\dot{\bm{Z}}^{h}_{t} \coloneqq \frac{1}{hI}b^{h}(hI\bm{Z}^{h}_{t},\bm{\xi}^{h}_{t})$
and 
$\dot{\bm{Z}}^{\sigma}_{t} \coloneqq \frac{1}{I}b(I\bm{Z}^{\sigma}_{t},\bm{\xi}^{1}_{t/\sigma})$.
From \eqref{eqn:sup_equilibrium_1}, for any $T>0$ and $\delta>0$, we obtain
\begin{equation}\label{eqn:sup_equilibrium_2}
\lim\limits_{h\rightarrow +\infty}\mathbb{P}\left\{ \sup\limits_{0\leq t\leq T}\norm{\bm{Z}^{h}_{t}-\overline{\bm{x}}(t)/I}>\delta\right\} = 0.
\end{equation}
Effectively then, scaling time by $\sigma = \frac{1}{h}$ is equivalent to scaling the system size by $h$.

Note that $\dot{\bm{Z}}^{h}_t$ and $\dot{\overline{\bm{x}}}(t)$ are dependent on the parameter $a\in(0,1)$ through functions $f^{0,h}_{\iota,a}(\bm{x},\bm{\xi}^h)$ and $f^{h}_{\iota,a}(\bm{x},\bm{\xi}^h)$ that are used in definition \eqref{eqn:definition_Q}.
Equation \eqref{eqn:sup_equilibrium_2} holds for any given $0<a<1$. 
Because of the Lipschitz behavior of $\dot{\bm{Z}}^{h}_{t}$ and $\dot{\overline{\bm{x}}}(t)$ on $0<a<1$, $\lim_{a\downarrow 0}d \dot{\bm{Z}}^{h}_{t}/d a = 0$ and $\lim_{a\downarrow 0}d \dot{\overline{\bm{x}}}(t)/d a = 0$,  equation \eqref{eqn:sup_equilibrium_2} holds in the limit as $a\rightarrow 0$. 
Also, if $\bm{Z}^{h}_{0}=\bm{Z}^{\varphi,h}(0)$, and $\overline{\bm{x}}(0)/I=\bm{z}^{\varphi,h}(0)$, then $\lim_{a\downarrow 0}\overline{\bm{x}}(t)/I = \lim_{h\rightarrow +\infty}\bm{z}^{\varphi,h}(t)$
and $\lim_{h\rightarrow +\infty}\lim_{a \downarrow 0}\bm{Z}^{h}_{t} = \lim_{h\rightarrow +\infty}\bm{Z}^{\varphi,h}(t)$.

For any $h\in\mathbb{N}_+$, the existence of $\bm{z}^{\varphi,h}(t)$, $t\geq 0$, satisfying \eqref{eqn:differential} with given $\bm{z}^{\varphi,h}(0)$ can be proved along similar lines as $\overline{\bm{x}}(t)$ by introducing a function $y_a(u)$ and then invoking Picard's Existence Theorem for initial value problems.

Recall that $\bm{Z}^{\varphi,h}(t)$ is the vector of proportions of sub-processes under policy $\varphi$ at time $t$ and $\bm{z}^{\varphi,h}(t)$ is given by \eqref{eqn:differential} (as defined in Section~\ref{sec:asym_opt}).  
Then, for any $T>0$ and $\delta>0$,
\begin{equation}
\lim\limits_{h\rightarrow +\infty}\mathbb{P}\left\{ \sup\limits_{0\leq t\leq T}\norm{\bm{Z}^{\varphi,h}(t)-\bm{z}^{\varphi,h}(t)}>\delta\right\} = 0,
\end{equation}
which leads to \eqref{eqn:convergence}.
This proves the proposition.

\endproof

\section{Proof of Lemma~\ref{lemma:zero_derivative}}\label{app:prop:equilibrium_vector}



We start with showing the existence of $\bm{z}^{\bar{\varphi}(o)}$.
\begin{lemma}
The limit $ \bm{z}^{\bar{\varphi}(o)}=\lim_{h\rightarrow +\infty} \bm{z}^{\bar{\varphi}(o),h}$ exists
\end{lemma}
\proof{Proof.}
Recall that the policy $\bar{\varphi}(o)$ is generated by Algorithm~\ref{algo:varphi_gamma} in Section~\ref{subsec:decomposable}, where $o$ is a given PS pair ranking.
Note that, as mentioned in Section~\ref{subsec:asym_regime},
the optimization problem consisting of the $hI$ sub-processes associated with $hI$ sub-patterns, coupled through constraints described in \eqref{eqn:constraint:action:h}-\eqref{eqn:constraint:resources:zero:h} can be analyzed and relaxed along the same lines as in Section~\ref{sec:relaxation}. 
Algorithm~\ref{algo:varphi_gamma} can also be applied directly to the problem scaled by $h$. To clarify, for the scaled problem described by \eqref{eqn:objective:h}-\eqref{eqn:constraint:resources:zero:h}, 
the action and capacity constraints can be relaxed to
\begin{equation}\label{eqn:relax:constraint:action:h}
\lim\limits_{t\rightarrow +\infty}\mathbb{E}\left[\sum\limits_{i\in\mathscr{P}_{\ell}}\sum\limits_{k\in[h]}a^{\phi}_{i,k}(\bm{N}^{\phi}_h(t))\right] = 1,~\forall \ell\in[L],
\end{equation}
and
\begin{equation}\label{eqn:relax:constraint:resources:h}
\lim\limits_{t\rightarrow +\infty}\mathbb{E}\left[
\sum\limits_{i\in[I]}\frac{w_{j,i}}{h}\sum\limits_{k\in[h]}\biggl(N^{\phi}_{i,k}(t)+a^{\phi}_{i,k}(\bm{N}^{\phi}_h(t))\biggr) \right]\leq C^0_j,~\forall j\in[J],
\end{equation}
which correspond to the relaxed constraints in \eqref{eqn:constraint:relax:action} and \eqref{eqn:constraint:relax:resources}, respectively.
Here, we first recall and rewrite the procedures of generating the policy $\bar{\varphi}(o)$ (in Algorithm~\ref{algo:varphi_gamma}):  initialize all PS pairs to be passive and then
sequentially activate the PS pairs according to their priorities defined by $o\in\mathscr{O}$.
In particular, for each PS pair $\iota\in[N]$, the action variables for the $h$ sub-patterns, $\alpha^{\bar{\varphi}(o)}_{i,k}(n_{\iota})$ ($k\in[h]$), are sequentially activated according to any permutation of $1,2,\ldots,h$, because all of them correspond to the same index $\Xi_{i_{\iota}}(\bm{\gamma},\bm{0})$ and the same priority in the PS pair ranking $o$. We assume without loss of generality that the $\alpha^{\bar{\varphi}(o)}_{i_{\iota},k}(n_{\iota})$ are activated in the order of $k=1,2,\ldots,h$.
The activating process continues until either a relaxed action or capacity constraint described in \eqref{eqn:relax:constraint:action:h} and \eqref{eqn:relax:constraint:resources:h}, respectively, achieves equality.  
Also, the mechanism described in \ref{case:equality:action} and \ref{case:equality:capacity} is directly applicable to this procedure by replacing \eqref{eqn:constraint:relax:action} and \eqref{eqn:constraint:relax:resources} with \eqref{eqn:relax:constraint:action:h} and \eqref{eqn:relax:constraint:resources:h}, respectively.

In this context, we still use Algorithm~\ref{algo:varphi_gamma} (including the mechanisms described in \ref{case:equality:action} and \ref{case:equality:capacity}) to indicate the procedure of generating the policy $\bar{\varphi}(o)$ in the scaled system.

For a given PS pair ranking $o$, the action variables $\alpha^{\bar{\varphi}(o)}_{i_{\iota},k}(n_{\iota})$ under the policy $\bar{\varphi}(o)$ are initialized to zero and activated sequentially from $\iota=1$ to $N$ (from the PS pair with the highest priority to the lowest), as described above. 
For clarity, we write $\kappa\in\{0\}\cup[N]$ to indicate the initial condition and the $N$ iterations in Algorithm~\ref{algo:varphi_gamma}, and define a sequence of intermediate policies, referred to as $\varphi(o,\kappa)$, for which $\alpha^{\varphi(o,0)}_{i,k}(\cdot) = 0$ ($i\in[I]$, $k\in[h]$) and $\alpha^{\varphi(o,N)}_{i,k}(\cdot) = \alpha^{\bar{\varphi}(o)}_{i,k}(\cdot)$.
In particular, for the $\kappa$th iteration and PS pairs $\iota\in[N]$, $\alpha^{\varphi(o,\kappa)}_{i_{\iota},k}(n_{\iota}) = \alpha^{\bar{\varphi}(o)}_{i_{\iota},k}(n_{\iota})$ ($k\in[h]$) if $\iota \leq \kappa$; and $\alpha^{\varphi(o,\kappa)}_{i_{\iota},k}(n_{\iota}) = 0$ if $\iota > \kappa$.

In the scaled system, each sub-process is a birth-and-death Markov process with finitely many states.
For any $o$, 
the Markov chain for the sub-process has a stationary distribution which is limiting as $t\rightarrow +\infty$ in the sense that the time-dependent distribution for any initial state converges to it. For a sub-process $(i,k)$ ($i\in[I]\backslash\{d(\ell):\ell\in[L]\}$, $k\in[h]$) and $\kappa\in[N]$, the stationary distribution of state $n\in\mathscr{N}_i^0\backslash \{0\}$ under policy $\varphi(o,\kappa)$ is 
\begin{equation}\label{eqn:stationary_distribution} 
\pi_{i,k}^{o,\kappa,h}(n) = \pi_{i,k}^{o,\kappa,h}(0)\prod_{n'=1}^n\frac{\alpha^{\varphi(o,\kappa)}_{i,k}(n'-1)h\lambda_{\ell(i)}}{n'\mu_i},
\end{equation}
with $\pi_{i,k}^{o,\kappa,h}(0)$ the stationary distribution for state $0$ normalized by $\sum_{n\in\mathscr{N}_i^0}\pi_{i,k}^{o,\kappa,h}(n) =1$. 

For $\kappa\in[N]$ and $h\in\mathbb{N}_+$, 
define the expected sum of the action variables with respect to the stationary distribution \eqref{eqn:stationary_distribution} associated with pattern $i\in[I]$
\begin{equation}
S_A^{o,h}(i,\kappa)\coloneqq \sum_{n\in\mathscr{N}^0_i}\sum\limits_{k\in[h]}\pi_{i,k}^{o,\kappa,h}(n)\alpha^{\varphi(o,\kappa)}_{i,k}(n)
\end{equation} 
 and, for $j\in[J]$, the expected sum of occupied capacities 
\begin{equation}
S_C^{o,h}(i,\kappa,j)\coloneqq \sum_{n\in\mathscr{N}^0_i}\frac{1}{h}\sum\limits_{k\in[h]}\pi_{i,k}^{o,\kappa,h}(n)w_{j,i}n,
\end{equation} 
under policy $\varphi(o,\kappa)$.

For $\kappa=1$ and the action variables for PS pairs $\iota>1$ initialized to zero, 
all the action variables $\alpha^{\varphi(o,\kappa)}_{i_{\kappa},k}(n_{\kappa})$ ($k\in[h]$) remain constant if no constraint in \eqref{eqn:relax:constraint:action:h} and \eqref{eqn:relax:constraint:resources:h} is violated as $h\rightarrow+\infty$; or some of them decrease in $h$ for sufficiently large $h$ if a constraint in \eqref{eqn:relax:constraint:action:h} or \eqref{eqn:relax:constraint:resources:h} achieves equality.
Slightly abusing notation, when $\kappa\in[N]$, we write $i_{\kappa}$, $n_{\kappa}$ and $\ell_{\kappa}$ to indicate $i_{\iota}$, $n_{\iota}$ and $\ell_{\kappa}$ with $\iota = \kappa$, respectively.
In the former case, the $\kappa$th element of $\bm{z}^{\bar{\varphi}(o)}$, $z^{\bar{\varphi}(o)}_{\kappa} = \lim\limits_{h\rightarrow+\infty} \lim\limits_{t\rightarrow +\infty}\mathbb{E}[Z^{\varphi(o,\kappa),h}_{\kappa}(t)]=\lim\limits_{h\rightarrow+\infty}\frac{1}{h}\sum\limits_{k\in[h]}\pi_{i_{\kappa},k}^{o,h,\kappa}(n_{\kappa})$,  $\lim\limits_{h\rightarrow+\infty}S_A^{o,h}(i,\kappa)$ and $\lim\limits_{h\rightarrow+\infty}S_C^{o,h}(i,\kappa,j)$ ($i\in[I]$, $j\in[J]$) exist.  
In the latter case, 
for any $H<+\infty$, there exists $h>H$ such that 
an equality is  achieved in \eqref{eqn:relax:constraint:action:h} or \eqref{eqn:relax:constraint:resources:h}; that is,
\begin{equation}\label{eqn:limit:action}
S_A^{o,h}(i_{\kappa},\kappa) = \sum\limits_{k\in[h]}\pi^{o,\kappa,h}_{i_{\kappa},k}(n_{\kappa})\alpha^{\varphi(o,\kappa),h}_{i_{\kappa},k}(n_{\kappa})=1,
\end{equation}
or, there is a $j\in[J]$ such that
\begin{equation}\label{eqn:limit:resources}
S_C^{o,h}(i_{\kappa},\kappa,j) = \frac{w_{j,i_{\kappa}}}{h}\sum\limits_{k\in[h]}\Bigl(\pi_{i_{\kappa},k}^{o,\kappa,h}(n_{\kappa}) n_{\kappa}+\pi_{i_{\kappa},k}^{o,\kappa,h}(n_{\kappa}+1)\bigl(n_{\kappa}+1\bigr)\Bigr)= C_j^0 - \frac{o(h)}{h},
\end{equation}
where $\pi_{i_{\kappa},k}^{o,\kappa,h}(\cdot)$ is the stationary distribution of sub-pattern $(i_{\kappa},k)$ that is a solution of \eqref{eqn:stationary_distribution} with given $\alpha^{\varphi(o,\kappa)}_{i_{\kappa},k}(\cdot)$.
In particular, the $o(h)$ on the right hand side of \eqref{eqn:limit:resources} corresponds to the second term on the left hand side of \eqref{eqn:relax:constraint:resources:h},
\begin{equation}
\lim\limits_{t\rightarrow+\infty}\mathbb{E}\Bigl[\sum\limits_{i\in[I]}w_{j,i}\sum\limits_{k\in[h]}a^{\varphi(o,\kappa)}_{i,k}\bigl(\bm{N}_h^{\varphi(o,\kappa)}(t)\bigr)\Bigr],
\end{equation}
which is equal to 
\begin{equation}\label{eqn:reserved_part}
R^h_{\kappa}(j)\coloneqq \sum\limits_{i\in[I]}w_{j,i}\sum\limits_{k\in[h]}\sum\limits_{n\in\mathscr{N}^0_{i}}\pi^{o,\kappa,h}_{i,k}(n)\alpha^{\varphi(o,\kappa)}_{i,k}(n)
\end{equation}
and  is bounded for any $h\in\mathbb{N}_+\cup\{+\infty\}$ because of the relaxed action constraints described in \eqref{eqn:relax:constraint:action:h}.

\begin{definition}\label{definition:saturated}
We say that the $\kappa$th iteration of Algorithm~\ref{algo:varphi_gamma} is \emph{saturated}, if, for any $H<+\infty$, there exists $h>H$ such that equality is achieved in \eqref{eqn:relax:constraint:action:h} or \eqref{eqn:relax:constraint:resources:h}.
\end{definition}
We now show the existence of $z^{\bar{\varphi}(o)}_{\kappa} =\lim\limits_{h\rightarrow+\infty}\frac{1}{h}\sum\limits_{k\in[h]}\pi_{i_{\kappa},k}^{o,h,\kappa}(n_{\kappa})$,  $\lim\limits_{h\rightarrow+\infty}S_A^{o,h}(i,\kappa)$ and $\lim\limits_{h\rightarrow+\infty}S_C^{o,h}(i,\kappa,j)$ ($i\in[I]$, $j\in[J]$)  in the saturated case.

Because the action variables for the same PS pair are  sequentially activated from $k=1,2,\ldots,h$, there exists a $k^*_{\kappa}(h)\in[h]$ such that, for all $k < k^*_{\kappa}(h)$, $\alpha^{\varphi(o,\kappa)}_{i_{\kappa},k}(n_{\kappa}) = 1$; for $k=k^*_{\kappa}(h)$, $\alpha^{\varphi(o,\kappa)}_{i_{\kappa},k}(n_{\kappa})\in(0,1]$; and, for $k > k^*_{\kappa}(h)$, $\alpha^{\varphi(o,\kappa)}_{i_{\kappa},k}(n_{\kappa}) = 0$. 
Let $\pi^{o,\kappa,h}_{i_{\kappa},k}(\cdot) = \pi_{i_{\kappa}}^{o,\kappa,h,+}(\cdot)$ for  $k<k^*_{\kappa}(h)-1$; and $\pi_{i_{\kappa}}^{o,\kappa,h,-}(\cdot)$ for $k > k^*_{\kappa}(h)$. For clarity, for given $h$ and $\kappa$, since the action variables $\alpha^{\varphi(o,\kappa),h}_{i_{\kappa},k}(n_{\kappa})$ are uniquely determined by the value
\begin{equation}
\rho^h_{\kappa} \coloneqq \Bigl(k^*_{\kappa}(h)+\alpha^{\varphi(o,\kappa),h}_{i_{\kappa},k^*_{\kappa}(h)}(n_{\kappa})\Bigr)/h,
\end{equation}
we refer to  $\rho^h_{\kappa}\in[0,1]$ as the \emph{decision} in the $\kappa$th iteration of Algorithm~\ref{algo:varphi_gamma}.

For $\rho\in[0,1]$, let 
$\pi^{o,\kappa,h}_{\rho}(n)$ be the solution of equations: for $n\in\mathscr{N}^0_{i_{\kappa}}$,
if $n \leq n_{\kappa}$,
\begin{equation}\label{eqn:stationary_distribution_rho:1} 
\pi_{\rho}^{o,\kappa,h}(n) = \pi_{\rho}^{o,\kappa,h}(0)\prod_{n'=1}^{n}\frac{h\lambda_{\ell_{\kappa}}}{n'\mu_{i_{\kappa}}};
\end{equation}
if $n = n_{\kappa}+1$,
\begin{equation}\label{eqn:stationary_distribution_rho:2} 
\pi_{\rho}^{o,\kappa,h}(n) = \pi_{\rho}^{o,\kappa,h}(0)\frac{(\rho h -\lfloor \rho h \rfloor)h\lambda_{\ell_{\kappa}}}{n \mu_{i_{\kappa}}}\prod_{n'=1}^{n-1}\frac{h\lambda_{\ell_{\kappa}}}{n'\mu_{i_{\kappa}}};
\end{equation}
if $n > n_{\kappa}+1$,
\begin{equation}\label{eqn:stationary_distribution_rho:3} 
\pi_{\rho}^{o,\kappa,h}(n) = 0;
\end{equation}
and 
\begin{equation}\label{eqn:stationary_distribution_rho:4} 
\sum_{n'\in\mathscr{N}_{i_{\kappa}}^0}\pi_{\rho}^{o,\kappa,h}(n') =1.
\end{equation}
The $\pi^{o,\kappa,h}_{\rho}(\cdot)$ represents the stationary distribution of a sub-process with respect to pattern $i_{\kappa}$ under a policy determined by $\rho\in[0,1]$ in the $\kappa$th iteration of Algorithm~\ref{algo:varphi_gamma}.
We refer to this policy as $\varphi_{\rho}(o,\kappa)$.
In particular, if $\rho = \rho^h_{\kappa}$, $\pi^{o,\kappa,h}_{\rho}(\cdot) = \pi^{o,\kappa,h}_{i_{\kappa},k^*_{\kappa}(h)}(\cdot)$.

Next we define, for $\rho\in[0,1]$,
\begin{equation}
f^{o,\kappa,h}_C(\rho,j) \coloneqq  
w_{j,i_{\kappa}}\sum\limits_{n=0}^{n_{\kappa}+1}\biggl(\rho\pi^{o,\kappa,h,+}_{i_{\kappa},k}(n)   + \frac{1}{h}\pi^{o,\kappa,h}_{\rho}(n)\biggr)n
\end{equation}
which is the expected  number of RUs of resource pool $j$ occupied by sub-patterns with respect to pattern $i_{\kappa}$ under policy $\varphi_{\rho}(o,\kappa)$. Also, if $\rho = \rho^h_{\kappa}$, $f^{o,\kappa,h}_C(\rho,j) = S^{o,h}_C(i_{\kappa},\kappa,j)$.
Since, for given $\rho\in[0,1]$, the action variables under policy $\varphi_{\rho}(o,\kappa)$ remain constant for any $h\in\mathbb{N}_+$, by solving the stationary distribution in \eqref{eqn:stationary_distribution_rho:1}-\eqref{eqn:stationary_distribution_rho:4},
we can show the existence of $\lim\limits_{h\rightarrow +\infty} f^{o,\kappa,h}_C(\rho,j)$ ($j\in[J]$). Similarly, define 
\begin{equation}\label{eqn:define:fA}
f^{o,\kappa,h}_A(\rho) \coloneqq \sum_{n=0}^{n_{\kappa}}\biggl(\lfloor\rho h \rfloor\pi_{i_{\kappa}}^{o,\kappa,h,+}(n) + \pi_{\rho}^{o,\kappa,h}(n)\mathds{1}_{n<n_{\kappa}} + \pi_{\rho}^{o,\kappa,h}(n) \Bigl(\rho h - \lfloor \rho h \rfloor\Bigr)\mathds{1}_{n=n_{\kappa}}\biggr),
\end{equation}
which is the expected sum of action variables with respect to pattern $i_{\kappa}$ under policy $\varphi_{\rho}(o,\kappa)$, so that $\lim\limits_{h\rightarrow +\infty} f^{o,\kappa,h}_A(\rho)$ exists.
Insertion of the stationary distribution in \eqref{eqn:stationary_distribution_rho:1}-\eqref{eqn:stationary_distribution_rho:4} shows that, for any $h\in\mathbb{N}_+\cup\{+\infty\}$, both $f^{o,\kappa,h}_A(\rho)$ and $f^{o,\kappa,h}_C(\rho,\cdot)$ are strictly increasing and Lipschitz continuous in $\rho\in[0,1]$. 
In particular, 
by taking limit of all terms in \eqref{eqn:define:fA}, we obtain $\lim\limits_{h\rightarrow +\infty} f^{o,\kappa,h}_A(\rho) = \rho n_{\kappa}\mu_{i_{\kappa}}/\lambda_{\ell_{\kappa}}$, and thus $f^{o,\kappa,h}_A(\rho)$ is Lipschitz continuous in $\rho$ for $h=+\infty$.

Define the inverse functions of $f^{o,\kappa,h}_A(\rho)$ and $f^{o,\kappa,h}_C(\rho,j)$ ($j\in[J]$) as $\rho = \bar{f}^{o,\kappa,h}_A(s)$ and $\rho=\bar{f}^{o,\kappa,h}_C(s,j)$, respectively. These are also strictly increasing and Lipschitz continuous in $s$.

We consider two cases:
\begin{enumerate}[label=$<\roman*>$]
\item $\lim\limits_{h\rightarrow +\infty}f^{o,\kappa,h}_A(1) < 1$ and, for all $j\in[J]$, $\lim\limits_{h\rightarrow +\infty}f^{o,\kappa,h}_C(1,j) < C_j^0$; and \label{case:not_saturated}
\item $\lim\limits_{h\rightarrow +\infty}f^{o,\kappa,h}_A(1) \geq 1$ or there is a $j\in[J]$ such that $\lim\limits_{h\rightarrow +\infty}f^{o,\kappa,h}_C(1,j) \geq C_j^0$. \label{case:saturated}
\end{enumerate}
Define $\mathscr{J}^{\kappa}$ to be the set of $j$ satisfying $\lim\limits_{h\rightarrow +\infty}f^{o,\kappa,h}_C(1,j) \geq C_j^0$, and 
 $\mathscr{J}^{\kappa}_1$ to be the subset of $\mathscr{J}^{\kappa}$ where for any $j\in\mathscr{J}^{\kappa}_1$, $\lim\limits_{h\rightarrow+\infty} f^{o,\kappa,h}_C(1,j) = 1$.

In Case~\ref{case:not_saturated}, there exists $H<+\infty$ such that, for all $h>H$ the constraints in \eqref{eqn:relax:constraint:action:h} and \eqref{eqn:relax:constraint:resources:h} are not violated --- this iteration is not saturated. In the non-saturated case, $z^{\bar{\varphi}(o)}_{\kappa} =\lim\limits_{h\rightarrow+\infty}\frac{1}{h}\sum\limits_{k\in[h]}\pi_{i_{\kappa},k}^{o,h,\kappa}(n_{\kappa})$,  $\lim\limits_{h\rightarrow+\infty}S_A^{o,h}(i,\kappa)$ and $\lim\limits_{h\rightarrow+\infty}S_C^{o,h}(i,\kappa,j)$ ($i\in[I]$, $j\in[J]$) exist as discussed earlier.

For Case~\ref{case:saturated}, 
if $\lim\limits_{h\rightarrow+\infty}f^{o,\kappa,h}_A(1) > 1$ or $\lim\limits_{h\rightarrow +\infty}f^{o,\kappa,h}_C(1,j) > C_j^0$ for some $j\in [J]$ (that is, $\mathscr{J}^{\kappa}\neq \emptyset$ and $\mathscr{J}^{\kappa}_1 \neq \mathscr{J}^{\kappa}$), then, because of the constraints in \eqref{eqn:relax:constraint:action:h} and \eqref{eqn:relax:constraint:resources:h}, there exists $H<+\infty$ such that, for all $h>H$, \eqref{eqn:limit:action} or \eqref{eqn:limit:resources} (for some $j\in\mathscr{J}^{\kappa}$) holds and thus $\lim\limits_{h\rightarrow +\infty} S^{o,h}_A(i_{\kappa},\kappa)$ or $\lim\limits_{h\rightarrow+\infty}S^{o,h}_C(i_{\kappa},\kappa,j)$ (for some $j\in\mathscr{J}^{\kappa}$) exists. 
Recall that $f^{o,\kappa,h}_A(\rho)$ and $f^{o,\kappa,h}_C(\rho,j)$ ($j\in[J]$) are strictly increasing and Lipschitz continuous in $\rho$, and so $\lim\limits_{h\rightarrow +\infty} \rho^h_{\kappa}$ exists.

If $\lim\limits_{h\rightarrow+\infty}f^{o,\kappa,h}_A(1) = 1$ or $\lim\limits_{h\rightarrow +\infty}f^{o,\kappa,h}_C(1,j) = C_j^0$ for all $j\in\mathscr{J}^{\kappa}$ (that is, $\mathscr{J}^{\kappa}=\mathscr{J}^{\kappa}_1$) and $\mathscr{J}^{\kappa}\neq \emptyset$,
then we write $f^{o,\kappa,h}_A(1) = 1 - \Delta^{\kappa}_A(h)$ and $f^{o,\kappa,h}_C(1,j) = C^0_j- \Delta^{\kappa}_C(h,j)$ ($j\in\mathscr{J}^{\kappa}$) with $\lim\limits_{h\rightarrow+\infty}\Delta^{\kappa}_A(h) = 0$ and $\lim\limits_{h\rightarrow +\infty}\Delta^{\kappa}_C(h,j)=0$.
In this context, recalling the constraints stated in \eqref{eqn:relax:constraint:action:h} and \eqref{eqn:relax:constraint:resources:h}, we obtain that 
\begin{equation}\label{eqn:add_for_peter:1}
S^{o,h}_A(i_{\kappa},\kappa)=f^{o,\kappa,h}_A(\rho^h_{\kappa}) \leq \min\{1-\Delta^{\kappa}_A(h),1\}
\end{equation}
and, for $j\in\mathscr{J}^{\kappa}$, 
\begin{equation}\label{eqn:add_for_peter:2}
S^{o,h}_C(i_{\kappa},\kappa,j)=f^{o,\kappa,h}_A(\rho^h_{\kappa},j)\leq \min\{C_j^0 -\frac{R^{h}_{\kappa}(j)}{h},C_j^0-\frac{R^h_{\kappa}(j)}{h}-\Delta^{\kappa}_C(h,j)\}.
\end{equation} 
Note that at least one of the 
equalities in equations~\eqref{eqn:add_for_peter:1} and \eqref{eqn:add_for_peter:2} is achieved.
Because $f^{o,\kappa,h}_A(\rho)$ and $f^{o,\kappa,h}_C(\rho,j)$ are strictly increasing and Lipschitz continuous in $\rho$,
\begin{equation}\label{eqn:rho}
\rho^h_{\kappa} = \min\Biggl\{\biggl\{\bar{f}^{o,\kappa,h}_A(1),\bar{f}^{o,\kappa,h}_A\Bigl(1-\Delta^{\kappa}_A(h)\Bigr)\biggr\}\cup \bigcup\limits_{j\in\mathscr{J}^{\kappa}}\biggl\{
\bar{f}^{o,\kappa,h}_C\Bigl(C_j^0-\frac{R^{h}_{\kappa}(j)}{h}\Bigr),\bar{f}^{o,\kappa,h}_C\Bigl(C_j^0-\frac{R^h_{\kappa}(j)}{h}-\Delta_C^{\kappa}(h,j)\Bigr)\biggr\}\Biggr\}.
\end{equation}
As explained in \eqref{eqn:reserved_part}, $R^h_{\kappa}(j)=o(h)$, and the functions $\bar{f}^{o,\kappa,h}_A(s)$ and $\bar{f}^{o,\kappa,h}_C(s,j)$ ($j\in[J]$) are Lipschitz continuous in $s$, so the limits of the bracketed arguments on the right hand side of \eqref{eqn:rho} all exist, and $\lim\limits_{h\rightarrow +\infty} \rho^h_{\kappa} $ exists. 
In summary, in the saturated case, $\lim\limits_{h\rightarrow +\infty}\rho^h_{\kappa}$ exists.

Since $\rho^h_{\kappa}$ exists, 
\begin{multline}
z^{\bar{\varphi}(o)}_{\kappa}=  \lim\limits_{h\rightarrow +\infty}\lim\limits_{t\rightarrow+\infty}\mathbb{E}[Z^{\bar{\varphi}(o),h}_{\kappa}(t)] 
=\lim\limits_{h\rightarrow +\infty}\frac{1}{h}\sum\limits_{k\in[h]}\pi_{i_{\kappa},k}^{o,\kappa,h}(n_{\kappa})\\
=\lim\limits_{h\rightarrow +\infty}\Biggl(\frac{1}{h}\biggl(\sum\limits_{k\in[k^*_{\kappa}(h)-1]}\pi_{i_{\kappa},k}^{o,\kappa,h}(n_{\kappa}) + \sum\limits_{k\in [h]/[k^*_{\kappa}(h)]}\pi_{i_{\kappa},k}^{o,\kappa,h}(n_{\kappa})\biggr)  + \frac{o(h)}{h} \Biggr)\\
= \lim\limits_{h\rightarrow +\infty}\Biggl(\rho^h_{\kappa}\pi_{i_{\kappa}}^{o,\kappa,h,+}(n_{\kappa}) + \Bigl(1-\rho^h_{\kappa}\Bigr)\pi_{i_{\kappa}}^{o,\kappa,h,-}(n_{\kappa}) + \frac{o(h)}{h}\Biggr),
\end{multline}
exists. 
Because  $f^{o,\kappa,h}_A(\rho)$ and  $f^{o,\kappa,h}_C(\rho,j)$ ($j\in[J]$) are Lipschitz continuous in $\rho\in[0,1]$, if $\lim\limits_{h\rightarrow +\infty}\rho^h_{\kappa}$ exists, then both
$\lim\limits_{h\rightarrow +\infty}S^{o,h}_A(i_{\kappa},\kappa) = \lim\limits_{h\rightarrow +\infty}f^{o,\kappa,h}_A(\rho^h_{\kappa})$ and $\lim\limits_{h\rightarrow +\infty}S^{o,h}_C(i_{\kappa},\kappa,j) = \lim\limits_{h\rightarrow +\infty}f^{o,\kappa,h}_C(\rho^h_{\kappa},j)$ ($j\in[J]$)  exist.
In addition, for $i\neq i_{\kappa}$, $\lim\limits_{h\rightarrow+\infty}S_A^{o,h}(i_{\kappa},\kappa) = 0$ and $\lim\limits_{h\rightarrow+\infty}S_C^{o,h}(i,\kappa,j) = 0$ ($j\in[J]$).

For the iteration $\kappa > 1$, we assume that $z^{\bar{\varphi}(o)}_{\kappa'} = \lim\limits_{h\rightarrow+\infty} \lim\limits_{t\rightarrow +\infty}\mathbb{E}[Z^{\bar{\varphi}(o),h}_{\kappa'}(t)] = \lim\limits_{h\rightarrow+\infty}\frac{1}{h}\sum\limits_{k\in[h]}\pi_{i_{\kappa'},k}^{o,\kappa',h}(n_{\kappa'})$, $\lim\limits_{h\rightarrow+\infty}S_A^{o,h}(i,\kappa')$ and $\lim\limits_{h\rightarrow+\infty}S_C^{o,h}(i,\kappa',j)$ ($i\in[I]$, $j\in[J]$) exist for all $\kappa' \in[N]$ with $\kappa' < \kappa$.
Along similar lines,  the action variables $\alpha^{\varphi(o,\kappa)}_{i_{\kappa},k}(n_{\kappa})$ ($k\in[h]$) are determined when the action variables of all PS pairs $\iota<\kappa$ are determined and other PS pairs $\iota'>\kappa$ remain passive. 
If, as $h\rightarrow +\infty$, no constraint in \eqref{eqn:relax:constraint:action:h} and \eqref{eqn:relax:constraint:resources:h} is violated or one of the constraints achieves equality in a PS pair $\iota < \kappa$ with $i_{\iota}=i_{\kappa}$, then 
there exists $H>0$ such that for all $h>H$, $\alpha^{\varphi(o,\kappa)}_{i_{\kappa},k}(n_{\kappa})$ ($k\in[h]$) is algebraically independent from $h$, and $z^{\bar{\varphi}(o)}_{\kappa}$, $\lim\limits_{h\rightarrow +\infty}S_A^{o,h}(i,\kappa)$ and $\lim\limits_{h\rightarrow +\infty}S_C^{o,h}(i,\kappa,j)$ ($i\in[I]$, $j\in[J]$) exist.

We focus on the saturated case where, for any $H<+\infty$, there is $h>H$ such that a constraint in \eqref{eqn:relax:constraint:action:h} or \eqref{eqn:relax:constraint:resources:h} achieves equality on the  $\kappa$th iteration; that is, either
\begin{multline}\label{eqn:limit3}
S_A^{o,h}(i_{\kappa},\kappa)=\sum\limits_{n\in\mathscr{N}_{i_{\iota}}^0}\sum\limits_{k\in[h]}\pi_{i_{\kappa},k}^{o,\kappa,h}(n)\alpha^{\varphi(o,\kappa),h}_{i_{\kappa},k}(n) = 1-\sum\limits_{\begin{subarray}~i\in\mathscr{P}_{\ell_{\kappa}}\\i\neq i_{\kappa}\end{subarray}}\sum\limits_{n\in\mathscr{N}^0_i}\sum\limits_{k\in[h]}\pi_{i,k}^{o,\kappa,h}(n)\alpha^{\varphi(o,\kappa),h}_{i,k}(n)\\
=1 - \sum\limits_{\begin{subarray}~i\in\mathscr{P}_{\ell_{\kappa}}\\i\neq i_{\kappa}\end{subarray}}S_A^{o,h}(i,\kappa-1),
\end{multline}
or,  for some $j\in[J]$,
\begin{multline}\label{eqn:limit4}
S_C^{o,h}(i_{\kappa},\kappa,j) = \sum\limits_{n\in\mathscr{N}_{i_{\kappa}}^0}\frac{w_{j,i_{\kappa}}}{h}\sum\limits_{k\in[h]}\pi_{i_{\kappa},k}^{o,\kappa,h}(n)n\\
= C_j^0 - \sum\limits_{\begin{subarray}~i\in[I]\\i\neq i_{\kappa}\end{subarray}}\sum\limits_{n\in\mathscr{N}^0_i}\frac{w_{j,i}}{h}\sum\limits_{k\in[h]}\pi_{i,k}^{o,\kappa,h}(n)n - \frac{R^h_{\kappa}}{h}
=C_j^0 - \sum\limits_{\begin{subarray}~i\in[I]\\i\neq i_{\kappa}\end{subarray}}S_C^{o,h}(i,\kappa-1,j) - \frac{o(h)}{h}.
\end{multline}

Along similar lines as discussed for Cases~\ref{case:not_saturated} and \ref{case:saturated} for $\kappa=1$, in the saturated case with respect to the $\kappa$th iteration, either
\begin{equation}\label{eqn:limit:fA}
\lim\limits_{h\rightarrow +\infty}f^{o,\kappa,h}_A(1) \geq 1-\lim\limits_{h\rightarrow+\infty}\biggl[\sum\limits_{\begin{subarray}~i\in\mathscr{P}_{\ell_{\kappa}}\\i\neq i_{\kappa}\end{subarray}}S_A^h(i,\kappa-1)\biggr],
\end{equation}
or, for some $j\in[J]$, 
\begin{equation}\label{eqn:limit:fC}
\lim\limits_{h\rightarrow +\infty}f^{o,\kappa,h}_C(1,j) \geq C_j^0 - \lim\limits_{h\rightarrow+\infty}\biggl[\sum\limits_{\begin{subarray}~i\in[I]\\i\neq i_{\kappa}\end{subarray}}S_C^h(i,\kappa-1,j)\biggr].
\end{equation} 
By slightly abusing notation, we also rewrite $\mathscr{J}^{\kappa}$ as the subset of $[J]$, where, for any $j\in\mathscr{J}^{\kappa}$ \eqref{eqn:limit:fC} is satisfied; and rewrite $\mathscr{J}^{\kappa}_1$, the subset of $\mathscr{J}^{\kappa}$, such that, for any $j\in\mathscr{J}^{\kappa}_1$, \eqref{eqn:limit:fC} achieves equality.

If \eqref{eqn:limit:fA} holds with strict inequality, or $\mathscr{J}^{\kappa} \neq \emptyset$ and $\mathscr{J}^{\kappa}_1 \neq \mathscr{J}^{\kappa}$, then there exists $H<+\infty$ such that for all $h>H$, \eqref{eqn:limit3} or \eqref{eqn:limit4} holds for a $j\in\mathscr{J}^{\kappa}$; by taking limit for both side of \eqref{eqn:limit3} or \eqref{eqn:limit4}, $\lim\limits_{h\rightarrow +\infty} S^{o,h}_A(i_{\kappa},\kappa)$ or $\lim\limits_{h\rightarrow +\infty} S^{o,h}_C(i_{\kappa},\kappa,j)$ exists.
This proves the existence of 
\begin{multline}\label{eqn:limit6}
\lim\limits_{h\rightarrow +\infty} \rho^h_{\kappa} = \lim\limits_{h\rightarrow +\infty}  \bar{f}^{o,\kappa,h}_A\Bigl(S^{o,h}_A(i_{\kappa},\kappa)\Bigr)=\lim\limits_{h\rightarrow +\infty} \bar{f}^{o,\kappa,h}_C\Bigl(S^{o,h}_C(i_{\kappa},\kappa,j)\Bigr)\\
= \bar{f}^{o,\kappa,h}_A\Bigl(\lim\limits_{h\rightarrow +\infty} S^{o,h}_A(i_{\kappa},\kappa)\Bigr) = \bar{f}^{o,\kappa,h}_C\Bigl(\lim\limits_{h\rightarrow +\infty}  S^{o,h}_C(i_{\kappa},\kappa,j)\Bigr),
\end{multline}
where the equality holds in the second line of \eqref{eqn:limit6} because of the Lipschitz continuity of the inverse functions, $\bar{f}^{o,\kappa,h}_A$ and $\bar{f}^{o,\kappa,h}_C$, of $f^{o,\kappa,h}_A$ and $f^{o,\kappa,h}_C$, respectively.

If \eqref{eqn:limit:fA} holds with equality or $\mathscr{J}^{\kappa}_1 = \mathscr{J}^{\kappa} \neq \emptyset$, then let
\begin{equation}\label{eqn:delat_kappa_A}
\Delta_A^{\kappa}(h) = 1 - \sum\limits_{\begin{subarray}~i\in\mathscr{P}_{\ell_{\kappa}}\\i\neq i_{\kappa}\end{subarray}}S_A^{o,h}(i,\kappa-1) - f^{o,\kappa,h}_A(1),
\end{equation}
and, for $j\in\mathscr{J}^{\kappa}_1$,
\begin{equation}\label{eqn:delat_kappa_C}
\Delta_C^{\kappa}(h,j) = C_j^0 - \sum\limits_{\begin{subarray}~i\in[I]\\i\neq i_{\kappa}\end{subarray}}S_C^{o,h}(i,\kappa-1,j) - f^{o,\kappa,h}_C(1,j),
\end{equation}
and so $\lim\limits_{h\rightarrow +\infty} \Delta_A^{\kappa}(h) = \lim\limits_{h\rightarrow +\infty} \Delta_C^{\kappa}(h,j) = 0$. 

Similar to the discussion for \eqref{eqn:rho}, because of the strict monotonicity of the functions $f^{o,\kappa,h}_A(\rho)$ and $f^{o,\kappa,h}_C(\rho,j)$ ($j\in[J]$) ($h\in\mathbb{N}_+$) $\rho^{\kappa}_h$ takes the minimum value such that at least one equality in 
\begin{equation}
S^{o,h}_A(i_{\kappa},\kappa) = f^{o,\kappa,h}_A(\rho^h_{\kappa}) \leq \min\Biggl\{1 - \sum\limits_{\begin{subarray}~i\in\mathscr{P}_{\ell_{\kappa}}\\i\neq i_{\kappa}\end{subarray}}S_A^{o,h}(i,\kappa-1)-\Delta_A^{\kappa}(h), 1 - \sum\limits_{\begin{subarray}~i\in\mathscr{P}_{\ell_{\kappa}}\\i\neq i_{\kappa}\end{subarray}}S_A^{o,h}(i,\kappa-1)\Biggr\},
\end{equation}
and, for $j\in\mathscr{J}^{\kappa}_1$,
\begin{multline}
S^{o,h}_C(i_{\kappa},\kappa,j) = f^{o,\kappa,h}_C(\rho^h_{\kappa},j)\\
\leq \min \Biggl\{ C_j^0 -\frac{R^{h}_{\kappa}(j)}{h} - \sum\limits_{\begin{subarray}~i\in[I]\\i\neq i_{\kappa}\end{subarray}}S_C^{o,h}(i,\kappa-1,j)-\Delta^{\kappa}_C(h,j),  C_j^0 -\frac{R^{h}_{\kappa}(j)}{h} - \sum\limits_{\begin{subarray}~i\in[I]\\i\neq i_{\kappa}\end{subarray}}S_C^{o,h}(i,\kappa-1,j) \Biggr\}
\end{multline}
holds; that is,
\begin{multline}\label{eqn:limit5}
\rho^h_{\kappa} = \min\Biggl\{\biggl\{\bar{f}^{o,\kappa,h}_A\biggl(1 - \sum\limits_{\begin{subarray}~i\in\mathscr{P}_{\ell_{\kappa}}\\i\neq i_{\kappa}\end{subarray}}S_A^{o,h}(i,\kappa-1) - \Delta^{\kappa}_A(h)\biggr),\bar{f}^{o,\kappa,h}_A\biggl(1 - \sum\limits_{\begin{subarray}~i\in\mathscr{P}_{\ell_{\kappa}}\\i\neq i_{\kappa}\end{subarray}}S_A^{o,h}(i,\kappa-1)\biggr)\biggr\}\cup\\
\bigcup\limits_{j\in\mathscr{J}^{\kappa}_1}\biggl\{\bar{f}^{o,\kappa,h}_C\Bigl(C_j^0- \frac{R^{h}_{\kappa}(j)}{h} - \sum\limits_{\begin{subarray}~i\in[I]\\i\neq i_{\kappa}\end{subarray}}S_C^{o,h}(i,\kappa-1,j)-\Delta^{\kappa}_C(h,j),j\Bigr), \\
\bar{f}^{o,\kappa,h}_C\Bigl(C_j^0  -\frac{R^{h}_{\kappa}(j)}{h}- \sum\limits_{\begin{subarray}~i\in[I]\\i\neq i_{\kappa}\end{subarray}}S_C^{o,h}(i,\kappa-1,j),j\Bigr)\biggr\}\Biggr\}.
\end{multline}
Because of the Lipschitz continuity of the functions $\bar{f}^{o,\kappa,h}_A$ and $\bar{f}^{o,\kappa,h}_C$, the limits of all the bracketed arguments on the right hand side of Equation~\eqref{eqn:limit5} exist and so $\lim\limits_{h\rightarrow +\infty} \rho^h_{\kappa}$ exists.
In summary, for the $\kappa$th iteration with $\kappa>1$, $\lim\limits_{h\rightarrow +\infty}\rho^h_{\kappa}$ always exists.

Similar to the analysis for the case with $\kappa=1$, it follows that
\begin{equation}
z^{\bar{\varphi}(o)}_{\kappa} = \lim\limits_{h\rightarrow +\infty} \biggr(\sum\limits_{n\in \mathscr{N}_{i_{\kappa}}^0}
\Bigl(\rho^h_{\kappa}\pi^{o,\kappa,h,+}_{i_{\kappa}}(n) + \Bigl(1-\rho^h_{\kappa}\Bigr)\pi^{o,\kappa,h,-}_{i_{\kappa}}(n)\Bigr) + \frac{o(h)}{h}\biggr)
\end{equation}
exists. Also, the Lipschitz continuity of $f^{o,\kappa,h}_A$ and $f^{o,\kappa,h}_C$ and the existence of $\lim\limits_{h\rightarrow +\infty} \rho^h_{\kappa}$ lead to the existence of $\lim\limits_{h\rightarrow+\infty}S^{o,h}_A(i_{\kappa},\kappa) = \lim\limits_{h\rightarrow +\infty}f^{o,\kappa,h}_A(\rho^h_{\kappa})$ and $\lim\limits_{h\rightarrow+\infty}S^{o,h}_C(i_{\kappa},\kappa,j)=\lim\limits_{h\rightarrow+\infty}f^{o,\kappa,h}_C(\rho^h_{\kappa},j)$ ($j\in[J]$).

For $i\neq i_{\kappa}$, $\lim\limits_{h\rightarrow +\infty}S_A^{o,h}(i,\kappa)$ and $\lim\limits_{h\rightarrow +\infty}S_C^{o,h}(i,\kappa,j)$ ($j\in[J]$) exist because $S_A^{o,h}(i,\kappa) = S_A^{o,h}(i,\kappa-1)$ and $S_C^{o,h}(i,\kappa,j)=S_C^{o,h}(i,\kappa-1,j)$ ($j\in[J]$).

In summary, $z^{\bar{\varphi}(o)}_{\kappa}$, $\lim\limits_{h\rightarrow +\infty}S_A^{o,h}(i,\kappa)$ and $\lim\limits_{h\rightarrow +\infty}S_A^{o,h}(i,\kappa,j)$ ($i\in[I]$, $j\in[J]$) still exist for $\kappa>1$.
The proof is completed by iteratively showing the existence of $z^{\bar{\varphi}(o)}_{\kappa}$, $\lim\limits_{h\rightarrow +\infty}S_A^{o,h}(i,\kappa)$ and $\lim\limits_{h\rightarrow +\infty}S_A^{o,h}(i,\kappa,j)$ ($i\in[I]$, $j\in[J]$) from $\kappa=1$ to $N$.

\endproof

\proof{Proof of Lemma~\ref{lemma:zero_derivative}.}

We now consider \eqref{eqn:zero_derivative}.
We will show that, when $h$ is sufficiently large, the value of $z^{\varphi,h}_{\iota}(t)$ as $t\rightarrow +\infty$ is independent of those of PS pairs $\iota'$ with $\iota' > \iota$ (PS pairs with lower priorities) at any time.

For $\iota\in\mathscr{N}_i^0$, $\iota'\in\mathscr{N}_{i'}^0$, $i\in\mathscr{P}_{\ell}$, $i'\in\mathscr{P}_{\ell'}$, and $\ell,\ell'\in[L]$,
the rates at which an active and passive sub-process in PS pair $\iota$ transition to PS pair $\iota'$ are, 
\begin{equation}
q^{1,h}(\iota,\iota')=\begin{cases}
h\lambda^0_{\ell}, & \text{if } i=i',\ n_{\iota'}=n_{\iota}+1,\\
n_{\iota}\mu_{i}, & \text{if } i=i',\ n_{\iota'}=n_{\iota}-1,\\
0,& \text{otherwise},
\end{cases}
\end{equation}
and
\begin{equation}
q^{0,h}(\iota,\iota')=\begin{cases}
n_{\iota}\mu_{i}, & \text{if } i=i',\ n_{\iota'}=n_{\iota}-1,\\
0,& \text{otherwise},
\end{cases}
\end{equation}
respectively. 
By substituting for $q^{1,h}(\cdot,\cdot)$ and $q^{0,h}(\cdot,\cdot)$ in \eqref{eqn:differential}, we obtain
\begin{multline}\label{eqn:trans_rate:1}
\frac{d\ z_{\iota}^{\varphi,h}(t)}{d\ t}=\\
\begin{cases}
\bigl(h\lambda^0_{\ell}\upsilon^{\varphi,h}_{\iota^-}(\bm{z}^{\varphi,h}(t))z^{\varphi,h}_{\iota^-}(t)+ n_{\iota^+}\mu_i z^{\varphi,h}_{\iota^+}(t)- h\lambda^0_{\ell}\upsilon^{\varphi,h}_{\iota}(\bm{z}^{\varphi,h}(t))z^{\varphi,h}_{\iota}(t)- n_{\iota}\mu_{i} z^{\varphi,h}_{\iota}(t)\bigr)Ih, & \text{if }
0 < n_{\iota} < |\mathscr{N}^0_{i}|-1,\\
- h\lambda^0_{\ell}\upsilon^{\varphi,h}_{\iota}(\bm{z}^{\varphi,h}(t))z^{\varphi,h}_{\iota}(t)Ih, & \text{if } n_{\iota} = 0,\\
- n_{\iota}\mu_{i} z^{\varphi,h}_{\iota}(t)Ih, & \text{otherwise},
\end{cases}
\end{multline}
where $\ell = \ell_{\iota}$, $i=i_{\iota}$, and $\iota^-$ and $\iota^+$ are the PS pairs with $i_{\iota^-}=i_{\iota^+} = i$ and satisfying $n_{\iota^-} = n_{\iota}-1$ and $n_{\iota^+}=n_{\iota}+1$, respectively.


We define the proportion of RUs of resource pool $j$ \emph{occupied} by pattern $i$ related to process $\bm{z}^{\varphi,h}(t)$ to be
\begin{equation}
\overline{\zeta}_{j,i}(\bm{z}^{\varphi,h}(t)) = \frac{w_{j,i} I}{C_j^0} \sum\limits_{\begin{subarray}~\iota'\in[N]\\ i_{\iota'} = i\end{subarray}}
n_{\iota'}z^{\varphi,h}_{\iota'}(t).
\end{equation}

Consider the PS pair $\iota$ with the highest priority among all pairs in $[N]$.
This must have $n_{\iota}=0$, because pair $(i,0)$ has the highest priority among all pairs $(i,n)$, $n\in\mathscr{N}_i^0$, for any $i\in[I]$.
From \eqref{eqn:trans_rate:1},

\begin{enumerate}[label=(\roman*)]
\item if  $\lambda_{\ell_{\iota}}^0 > u_{\iota}(\bm{z}^{\varphi,h}(t_0))$, where
\begin{equation}\label{eqn:lemma:1}
u_{\iota}(\bm{z}) \coloneqq n_{\iota^+}\mu_i z_{\iota^+},~t\geq 0, 
\end{equation}
then $z^{\varphi,h}_{\iota^+}(t)$ ($t\geq t_0$) will increase and $z^{\varphi,h}_{\iota}(t)$ will decrease, until $z^{\varphi,h}_{\iota}(t)$ becomes zero or a capacity constraint is about to be violated; that is, there exists a $j\in[J]$,\label{case:1}
\begin{equation}\label{eqn:lemma:2}
\sum\limits_{i'\in[I]}\overline{\zeta}_{j,i'}(z^{\varphi,h}_{\iota}(t))/(1-\epsilon^h_{j,\iota})=1;
\end{equation}
\item if $\lambda^0_{\ell_{\iota}}= u_{\iota}(\bm{z}^{\varphi,h}(t_0))$ and 
$z^{\varphi,h}_{\iota'}(t)=0$ for all  $\iota'>\iota^+$ with $i_{\iota'}=i_{\iota}$, then $z^{\varphi,h}_{\iota}(t)=z^{\varphi,h}_{\iota}(t_0)$ for any $t \geq t_0$;
\label{case:2}
\item if $u_{\iota}(\bm{z}^{\varphi,h}(t_0))\geq \lambda^0_{\ell_{\iota}} $  and $n_{\iota} = 0$, then 
$z^{\varphi,h}_{\iota'}(t)$ for all  $\iota'>\iota^+$ with $i_{\iota'}=i_{\iota}$ will decrease until $\lambda^0_{\ell_{\iota}}= u_{\iota}(\bm{z}^{\varphi,h}(t))$ and $z^{\varphi,h}_{\iota'}(t)=0$ for all such $\iota'$
 (invoking Case~\ref{case:2}).\label{case:3}
\end{enumerate}

In all of the cases \ref{case:1}-\ref{case:3},  for the PS pair $\iota$ with the highest priority, we consider a process starting from time $t_0$ and finishing at time $t=t_1$, at which either \eqref{eqn:lemma:2} holds or $z^{\varphi,h}_{\iota'}(t)=0$ for all $\iota'>\iota^+$ with $i_{\iota'}=i_{\iota}$.
In other words, this process ends in Case~\ref{case:1} or \ref{case:2}.
We refer to the process as the \emph{I-process} and time $t_1$ as its stopping time.
Note that the I-process may continue indefinitely with $t_1= +\infty$.



Define a \emph{stopping pair} of pattern $i$ at time $t$, denoted by $p_i(t)$, as the PS pair satisfying $z^{\varphi,h}_{\iota}(t)>0$ and $z^{\varphi,h}_{\iota'}(t)=0$ for all $\iota'>\iota$ with $i_{\iota'}=i_{\iota}$.


If the I-process ends in Case~\ref{case:1}, 
we then consider another process starting from its stopping time $t_1$.
If \eqref{eqn:lemma:2} holds for some $t_1\geq t_0$ and $j\in[J]$, $\lambda_{\ell_{\iota}}^{0} > u_{\iota}(\bm{z}^{\varphi,h}(t_1))$ and $\overline{\zeta}_{j,i'}(\bm{z}^{\varphi,h}(t_1)) > 0$ for some $i'\in[I]$ with $w_{j,i'}>0$ and $p_{i'}(t_1) > \iota$ (that is, resource pool $j$ is partially occupied by another pattern $i'$ of which the stopping pair $p_{i'}(t_1)$ has less priority than pair $\iota$), then 
\begin{enumerate}[label=(\roman*)]
\setcounter{enumi}{3}
\item $z^{\varphi,h}_{\iota^+}(t)$ and $u_{\iota}(\bm{z}^{\varphi,h}(t))$ ($t>t_1$)  will increase until   $\lambda_{\ell_{\iota}}^{0} = u_{\iota}(\bm{z}^{\varphi,h}(t))$ (invoking Case~\ref{case:3}); \label{case:4}
\item $z^{\varphi,h}_{\iota^+}(t)$ and $u_{\iota}(\bm{z}^{\varphi,h}(t))$ ($t>t_1$)  will increase until   $z^{\varphi,h}_{\iota}(t)=0$ (stopping pair $p_{i_{\iota}}(t)$ for pattern $i_{\iota}$ is no longer $\iota$); or \label{case:5}
\item $z^{\varphi,h}_{\iota^+}(t)$ and $u_{\iota}(\bm{z}^{\varphi,h}(t))$ ($t>t_1$)  will increase until   $\overline{\zeta}_{j,i'}(\bm{z}^{\varphi,h}(t)) = 0$ for all $i'$ with $w_{j,i'}>0$ and $p_{i'}(t) > \iota$, and $z^{\varphi,h}_{\iota'}(t)=0$ for all $\iota'>\iota^+$ with $i_{\iota'}=i_{\iota}$.\label{case:6}
\end{enumerate}

We refer to such a process starting from time $t_1$ until any of above three conditions is satisfied as the \emph{II-process}, and refer to its stopping time as $t_2$. 
The II-process may continue indefinitely with $t_2 =  +\infty$.
The II-process describes the case where the RUs occupied by pattern $i'$ can be taken by pattern $i_{\iota}$ if the stopping pair of $i_{\iota}$ (that is, PS pair $\iota$) has higher priority than that of $i'$ and $\lambda_{\ell_{\iota}}^{0} > u_{\iota}(\bm{z}^{\varphi,h}(t))$.

We now generalize above discussions to a PS pair $\iota\in[N]$ with $z^{\varphi,h}_{\iota'}(t_0) = 0$ for all $\iota' < \iota$; that is, PS pair $\iota$ is the one with the highest priority among $\iota'\in[N]$ with $z^{\varphi,h}_{\iota'}(t_0)>0$.
We can go through a similar argument by generalizing the definition of $u_{\iota}(\bm{z}^{\varphi,h}(t))$ so that
\begin{equation}\label{eqn:u}
u_{\iota}(\bm{z}^{\varphi,h}(t))\coloneqq \sum_{\iota'\leq\iota, \ell_{\iota'}=\ell}n_{\iota'}\mu_{i_{\iota'}}z^{\varphi,h}_{\iota'}(t) + n_{\iota^+}\mu_i z^{\varphi,h}_{\iota^+}(t).
\end{equation}
If $n_{\iota} > 0$, the description of the I-process is completed by adding
\begin{enumerate}[label=(\roman*)]
\setcounter{enumi}{6}
\item if $u_{\iota}(\bm{z}^{\varphi,h}(t_0))\geq \lambda^0_{\ell_{\iota}} $  and $n_{\iota} > 0$, then 
$z^{\varphi,h}_{\iota'}(t)$ for all  $\iota'>\iota^+$ with $i_{\iota'}=i_{\iota}$ will decrease until \label{case:7}
\begin{enumerate}[label=\ref{case:7}.\arabic*]
\item $\lambda^0_{\ell_{\iota}}= u_{\iota}(\bm{z}^{\varphi,h}(t))$ and $z^{\varphi,h}_{\iota'}(t)=0$ for all  $\iota'>\iota^+$ with $i_{\iota'}=i_{\iota}$ (invoking Case~\ref{case:2}); or \label{case:7:1}
\item 
$\lambda^0_{\ell_{\iota}} > u_{\iota^-}(\bm{z}^{\varphi,h}(t))$ and $z^{\varphi,h}_{\iota^-}(t) > 0$, so the stopping pair $p_{i_{\iota}}(t)$ of pattern $i_{\iota}$ becomes $\iota^-$.\label{case:7:2}
\end{enumerate}
\end{enumerate}

If the I-process for PS pair $\iota$ ends in Case~\ref{case:7:2} at time $t_1$, since $\iota^-$ always has higher priority than $\iota$, $\iota^-$ will be considered as the PS pair with the highest priority and $z^{\varphi,h}_{\iota^-}(t_1)>0$ at time $t_1$. 
In this case, because $\lambda^0_{\ell_{\iota}} > u_{\iota^-}(\bm{z}^{\varphi,h}(t_1))$ is guaranteed in Case~\ref{case:7:2}, the value $z^{\varphi,h}_{\iota^-}(t)$ for $t\geq t_1$ will be considered in Cases \ref{case:3} or \ref{case:7}. Iteratively, once $z^{\varphi,h}_{\iota}(t)$ goes into Case~\ref{case:7}, there will be an $\iota' \leq \iota$ with $i_{\iota'}=i_{\iota'}$, such that $z^{\varphi,h}_{\iota'}(t)$ ends in Case~\ref{case:2}; or continues indefinitely  in Case~\ref{case:3} or \ref{case:7}.



We firstly consider the situation with $t_1,t_2 < +\infty$.
From above description of the I and II-processes, the process $z^{\varphi,h}_{\iota}(t)$ ($t\geq t_0$) will enter one of two regimes:
\begin{itemize}
\item $\lambda^0_{\ell_{\iota}} = u_{\iota}(z^{\varphi,h}_{\iota}(t))$ and $z^{\varphi,h}_{\iota'}(t)=0$ for all $\iota'>\iota^+$ with $i_{\iota'}=i_{\iota}$;  or
\item $\lambda^0_{\ell_{\iota}} > u_{\iota}(z^{\varphi,h}_{\iota}(t))$, Equation \eqref{eqn:lemma:2} holds, there exists a $j\in[J]$ such that
$\overline{\zeta}_{j,i'}(\bm{z}^{\varphi,h}(t)) = 0$ for all $i'$ with $w_{j,i'}>0$ and $p_{i_{\iota'}}(t) > \iota$, and $z^{\varphi,h}_{\iota'}(t)=0$ for all $\iota'>\iota^+$ with $i_{\iota'}=i_{\iota}$.
\end{itemize}
Intuitively, in the first regime, the arrival rate of RT $\ell_{\iota}$ is balanced by the total service rate supported by all PS pairs with higher priorities than $\iota^+$ and PS pair $\iota^+$; and, in the second regime, the capacity constraints are about to be violated, which helps balance the arrival rate by forcing the sub-processes to be passive (with zero arrival rate) when a resource pool is fully occupied.

From the differential equation stated in \eqref{eqn:trans_rate:1}, the value of $z^{\varphi,h}_{\iota}(t)$ will not change once it enters any of the two regimes.
We refer to them as the \emph{stability regimes}, and say the value of $z^{\varphi,h}_{\iota}(t)$ for the stopping pair $\iota$ of pattern $i_{\iota}$ at time $t$ becomes \emph{stable} if $z^{\varphi}_{\iota}(t)=z^{\varphi}_{\iota}(t')$ for any $t'\geq t$.

If the I or II-process continues indefinitely, then 
it will continue indefinitely in one of Cases~\ref{case:1}, \ref{case:3}-\ref{case:7}.
We recall the $\upsilon^{\varphi,h}_{\iota'}(\bm{z}^{\varphi,h}(t))$ defined in \eqref{eqn:probability_active}, representing
the proportion of activated sub-processes in PS pair $\iota'\in[N]$.
From its definition, if $z^{\varphi,h}_{\iota}(t)>0$, then, for all $\iota'>\iota^+$ with $i_{\iota'}=i_{\iota}$,
$\upsilon^{\varphi,h}_{\iota'}(\bm{z}^{\varphi,h}(t))=0$.
If the process $z^{\varphi,h}_{\iota}(t)$ ($\iota\in[N]$) continues indefinitely for time period $t\in[T,+\infty)$ in any of Cases~\ref{case:1}, \ref{case:3}-\ref{case:7}, then, for all $\iota'>\iota^+$ with $i_{\iota'}=i_{\iota}$ and $t\in[T,\infty)$, $\upsilon^{\varphi,h}_{\iota'}(\bm{z}^{\varphi,h}(t))=0$.
Thus, for these $\iota'$, $\lim\limits_{t\rightarrow+\infty}z^{\varphi,h}_{\iota'}(t)=0$. 
In other words, $\lim\limits_{t\rightarrow +\infty} (z^{\varphi,h}_{\iota}(t)+z^{\varphi,h}_{\iota^+}(t)) = 1$:
all the sub-processes of pattern $i_{\iota}$ are transitioning between two states $n_{\iota}$ and $n_{\iota^+}$ indefinitely. In this case, from the differential equation stated in \eqref{eqn:trans_rate:1}, $\lim\limits_{t\rightarrow +\infty} z^{\varphi,h}_{\iota}(t)$ exists, and is a point that satisfies one of the stability regimes.

Accordingly, the stable value of $z^{\varphi,h}_{\iota}(t)$ is independent of the values of $z^{\varphi,h}_{\iota'}(t')$ for any $t'\geq 0$ and $\iota'>\iota$ (PS pairs with lower priorities).
Note that if $\iota$ becomes stable without delivering the stopping pair role to $\iota'>\iota$ with $i_{\iota'}=i_{\iota}$, then $\iota'$ 
is removed from future consideration, because the value $z^{\varphi,h}_{\iota'}(t)$ remains constant in the future.

From the definition of the stability regimes, for any $\bm{\epsilon}\in\mathscr{E}^{+\infty}$ and $\bm{z}^{\varphi,h}(0)\in\mathscr{Z}$, there exists $H\in\mathbb{N}_+$ such that, for all $h>H$, the stable values of $\bm{z}^{\varphi,h}(t)$ (that is, $\lim\limits_{t\rightarrow +\infty}\bm{z}^{\varphi,h}(t)$) can be iteratively calculated from PS pair $\iota=1$ to $N$, which are independent of $h$.
Thus, $\bm{z}^{\varphi}\coloneqq\lim\limits_{h\rightarrow +\infty} \lim\limits_{t\rightarrow +\infty}\bm{z}^{\varphi,h}(t)$ exists with any initial point in $\mathscr{Z}$.
Note that the existence of positive elements of $\bm{\epsilon}$ is used to prioritize PS pairs in the asymptotic regime, which is crucial for the stable value of $z^{\varphi,h}_{\iota}(t)$ to be independent from those of other stopping pairs with lower priorities than $\iota$.

Recall that the $\bm{z}^{\varphi}$ is continuous in $\bm{\epsilon}\in\mathscr{E}^{+\infty}$ and is a bounded vector in the probability simplex $\mathscr{Z}$.
For the PS pair $\iota=1$,  $z^{\varphi}_{\iota}$ is continuous and increasing in $\epsilon_{j,\iota}\in(0,1]$ ($j\in\mathscr{J}_{i_{\iota}}$), then $\lim_{\lVert \bm{\epsilon} \rVert\to 0}z^{\varphi}_{\iota}$ exists.
If $\lim_{\lVert \bm{\epsilon} \rVert\to 0}z^{\varphi}_{\iota'}$ exists for any $\iota' < \iota$, then, since $z^{\varphi}_{\iota}$ is continuous and increasing in $\epsilon_{j,\iota}\in(0,1]$ ($j\in\mathscr{J}_{i_{\iota}}$), $\lim_{\lVert \bm{\epsilon} \rVert\to 0}z^{\varphi}_{\iota}$ also exists.
Accordingly, $\lim_{\lVert \bm{\epsilon} \rVert\to 0}\bm{z}^{\varphi}$ exists,  for which  the calculating procedure is the same as calculating $\bm{z}^{\bar{\varphi}(o)}$.
The lemma has been proved.



\endproof


\section{Proof of Theorem~\ref{theorem:main}}\label{app:theorem:main}
\proof{Proof.}
If the capacity constraints described in \eqref{eqn:constraint:resources:h} (or equivalently \eqref{eqn:constraint:resources}) are decomposable with decomposable values $\pmb{\gamma}\in\mathbb{R}^J_0$ of  the multipliers in the asymptotic regime, then there exist $\bm{\nu}\in\mathbb{R}^L$ and a PS-pair ranking $o\in\mathscr{O}(\pmb{\gamma},\bm{\nu})$ such that the complementary slackness is satisfied by the policy $\bar{\varphi}(o)$ and the multipliers $\pmb{\gamma}$.
The policy $\bar{\varphi}(o)$
is optimal for the relaxed problem.
Together with Theorem~\ref{theorem:main_second}, the index policy $\varphi$ derived from the same ranking $o$ is also optimal in the asymptotic regime.
The theorem is proved.

\endproof

\section{Proof of Corollary~\ref{coro:main}}\label{app:coro:main}

\proof{Proof.}
If the system is weakly coupled,
then, from \eqref{eqn:equal_opt:3}-\eqref{eqn:equal_opt:7} in Appendix~\ref{app:prop:equal_opt}, 
the ranking of PS pairs $o$ following the descending order of $\Xi^*_{\iota}$ defined in \eqref{eqn:main:index} for $\iota\in[N]$ leads to an optimal policy $\bar{\varphi}(o)$ of the relaxed problem.
Together with Theorem~\ref{theorem:main}, the index policy $\varphi$ derived from such a PS-pair ranking $o$ is asymptotically optimal.

\endproof

\section{Settings of Simulations}
\label{app:simulation:opt}

\subsection{Settings of Simulations in Figure~\ref{fig1:opt_v1}}
The simulations whose results are exhibited in Figure~\ref{fig1:opt_v1} have four request types and fourteen resource pools, with capacities and cost rates per RU given in Table~\ref{table:fig1:opt_v1}.

\begin{table}[t]
\caption{Capacities and Cost Rates for RUs of Different Pools.}\label{table:fig1:opt_v1}
\renewcommand{\arraystretch}{1.2}
\centering
\begin{tabular}{|c|c|c|c|c|c|c|c|}
\hline
Resource Pool ($j$) & $1$ & $2$ & $3$ & $4$ & $5$ & $6$ & $7$  \\\hline
Capacity ($C_j^0$) & $5$ & $7$ & $6$ & $8$ & $6$ & $7$ & $6$\\\hline
Cost Rate ($\varepsilon_j$) & $9.046$ &
$4.995$ & $0.679$ & $2.761$ & $9.010$ & $4.775$ & $3.033$   \\\hline
Resource Pool ($j$) & $8$ & $9$ & $10$ & $11$ & $12$ & $13$ & $14$ \\\hline
Capacity ($C_j^0$) & $6$& $9$& $8$& $8$ & $5$& $8$& $5$
\\\hline
Cost Rate ($\varepsilon_j$) & $5.033$ & $3.318$ & $4.686$ & $3.302$ & $0.938$ & $6.770$ & $7.775$
\\\hline
\end{tabular}
\end{table}

\begin{itemize}
\item $\mathscr{P}_1=\{1,2,\ldots,12,13\}$,  $\lambda^0_1=1.275$, $\mu_{1}=\mu_{2}=\ldots=\mu_{12}=0.255$, $R_1=4026.22$;
\item $\mathscr{P}_2=\{14,15,\ldots,21,22\}$,   $\lambda^0_2=1.203$, $\mu_{14}=\mu_{15}=\ldots=\mu_{21}=0.241$,  $R_2=3871.33$;
\item $\mathscr{P}_3=\{23,24,\ldots,32,33\}$,  $\lambda^0_3=1.134$, $\mu_{23}=\mu_{24}=\ldots=\mu_{32}=0.227$,  $R_3=3731.69$;
\item $\mathscr{P}_4=\{34,35,\ldots,40,41\}$, $\lambda^0_4=1.314$, $\mu_{34}=\mu_{35}=\ldots=\mu_{40}=0.263$, $R_4=3242.58$;
\end{itemize}
where patterns $13,22,33,41$ are dummy patterns for blocking requests.
The weight vectors $\bm{w}_i$ of patterns  $i\in[41]\backslash \{13,22,33,41\}$ are given in Table~\ref{table:fig1:opt_v1:weights}, where $\bm{e}^{k}\in \{0,1\}^{J}$, $k\in[J]$, is a vector with all zero entries except its $k$th entry $e^k_k = 1$.
\begin{table}[t]
\caption{Weight Vectors of Patterns for Four Different RTs.}\label{table:fig1:opt_v1:weights}
\renewcommand{\arraystretch}{1.2}
\centering
\begin{tabular}{|c|c|c|c|c|c|c|c|}
\hline
$\bm{w}_1$ & $2\bm{e}^{12}+3\bm{e}^6+4\bm{e}^{13}$ & $\bm{w}_{14}$ & $\bm{e}^{4}+3\bm{e}^8+\bm{e}^{12}$  & $\bm{w}_{23}$ & $\bm{e}^{11}+\bm{e}^{13}+2\bm{e}^{6}$ & $\bm{w}_{34}$ & $2\bm{e}^{8}+2\bm{e}^{4}+2\bm{e}^{3}$
\\\hline
$\bm{w}_2$ & $2\bm{e}^{12}+3\bm{e}^1+4\bm{e}^{13}$  & $\bm{w}_{15}$ & $\bm{e}^{4}+3\bm{e}^8+\bm{e}^{7}$  & $\bm{w}_{24}$ & $\bm{e}^{11}+\bm{e}^{13}+2\bm{e}^{1}$ & $\bm{w}_{35}$ & $2\bm{e}^{8}+2\bm{e}^{4}+2\bm{e}^{12}$
\\\hline
$\bm{w}_3$ & $2\bm{e}^{7}+3\bm{e}^6+4\bm{e}^{13}$  & $\bm{w}_{16}$ & $\bm{e}^{4}+3\bm{e}^8+\bm{e}^{11}$  & $\bm{w}_{25}$ & $\bm{e}^{9}+\bm{e}^{13}+2\bm{e}^{6}$ & $\bm{w}_{36}$ & $2\bm{e}^{8}+2\bm{e}^{4}+2\bm{e}^{7}$
\\\hline
$\bm{w}_4$ & $2\bm{e}^{7}+3\bm{e}^1+4\bm{e}^{13}$  & $\bm{w}_{17}$ & $\bm{e}^{4}+3\bm{e}^8+\bm{e}^{9}$  & $\bm{w}_{26}$ & $\bm{e}^{11}+\bm{e}^{13}+2\bm{e}^{1}$ & $\bm{w}_{37}$ & $2\bm{e}^{8}+2\bm{e}^{4}+2\bm{e}^{11}$
\\\hline
$\bm{w}_5$ & $2\bm{e}^{11}+3\bm{e}^6+4\bm{e}^{13}$  & $\bm{w}_{18}$ & $\bm{e}^{4}+3\bm{e}^8+\bm{e}^{10}$  & $\bm{w}_{27}$ & $\bm{e}^{10}+\bm{e}^{13}+2\bm{e}^{6}$ & $\bm{w}_{38}$ & $2\bm{e}^{8}+2\bm{e}^{4}+2\bm{e}^{2}$
\\\hline
$\bm{w}_6$ & $2\bm{e}^{11}+3\bm{e}^1+4\bm{e}^{13}$  & $\bm{w}_{19}$ & $\bm{e}^{4}+3\bm{e}^8+\bm{e}^{2}$  & $\bm{w}_{28}$ & $\bm{e}^{10}+\bm{e}^{13}+2\bm{e}^{1}$ & $\bm{w}_{39}$ & $2\bm{e}^{8}+2\bm{e}^{4}+2\bm{e}^{14}$
\\\hline
$\bm{w}_7$ & $2\bm{e}^{2}+3\bm{e}^6+4\bm{e}^{13}$  & $\bm{w}_{20}$ & $\bm{e}^{4}+3\bm{e}^8+\bm{e}^{14}$  & $\bm{w}_{29}$ & $\bm{e}^{2}+\bm{e}^{13}+2\bm{e}^{6}$ & $\bm{w}_{40}$ & $2\bm{e}^{8}+2\bm{e}^{4}+2\bm{e}^{5}$
\\\hline
$\bm{w}_8$ & $2\bm{e}^{2}+3\bm{e}^1+4\bm{e}^{13}$  & $\bm{w}_{21}$ & $\bm{e}^{4}+3\bm{e}^8+\bm{e}^{5}$  & $\bm{w}_{30}$ & $\bm{e}^{2}+\bm{e}^{13}+2\bm{e}^{1}$ &&
\\\hline
$\bm{w}_9$ & $2\bm{e}^{14}+3\bm{e}^6+4\bm{e}^{13}$  && & $\bm{w}_{31}$ & $\bm{e}^{5}+\bm{e}^{13}+2\bm{e}^{6}$ &&
\\\hline
$\bm{w}_{10}$ & $2\bm{e}^{14}+3\bm{e}^1+4\bm{e}^{13}$  && & $\bm{w}_{32}$ & $\bm{e}^{5}+\bm{e}^{13}+2\bm{e}^{1}$&&
\\\hline
$\bm{w}_{11}$ & $2\bm{e}^{5}+3\bm{e}^6+4\bm{e}^{13}$  &&&&&&
\\\hline
$\bm{w}_{12}$ & $2\bm{e}^{5}+3\bm{e}^1+4\bm{e}^{13}$  &&&&&&
\\\hline
\end{tabular}
\end{table}

\subsection{Settings of Simulations in Figure~\ref{fig1:opt_v3}}
The simulations whose results are exhibited in Figure~\ref{fig1:opt_v3} have two request types and six resource pools, with capacities and cost rates per RU given in Table~\ref{table:fig1:opt_v3}.

\begin{table}[t]
\caption{Capacities and Cost Rates for RUs of Different Pools.}\label{table:fig1:opt_v3}
\renewcommand{\arraystretch}{1.2}
\centering
\begin{tabular}{|c|c|c|c|c|c|c|}
\hline
Resource Pool ($j$) & $1$ & $2$ & $3$ & $4$ & $5$ & $6$  \\\hline
Capacity ($C_j^0$) & $8$ & $8$ & $9$ & $9$ & $9$ & $6$ \\\hline
Cost Rate ($\varepsilon_j$) & $5.684$ &
$7.249$ & $0.224$ & $4.969$ & $6.874$ & $8.539$
\\\hline
\end{tabular}
\end{table}
\begin{itemize}
\item $\mathscr{P}_1=\{1,2,3\}$,  $\lambda^0_1=1.385$, $\mu_{1}=\mu_{2}=0.277$, $R_1=3635.69$;
\item $\mathscr{P}_2=\{4,5,6,7,8\}$,   $\lambda^0_2=1.317$, $\mu_{4}=\mu_{5}=\mu_6=\mu_{7}=0.263$,  $R_2=3758.85$;
\end{itemize}
where patterns $3,8$ are dummy patterns for blocking requests.
The weight vectors $\bm{w}_i$ of patterns  $i\in[8]\backslash \{3,8\}$ are given in Table~\ref{table:fig1:opt_v3:weights}. 
\begin{table}[t]
\caption{Weight Vectors of Patterns for Two Different RTs.}\label{table:fig1:opt_v3:weights}
\renewcommand{\arraystretch}{1.2}
\centering
\begin{tabular}{|c|c|c|c|}
\hline
$\bm{w}_1$ & $\bm{e}^{1}+\bm{e}^4+\bm{e}^{3}$ & $\bm{w}_4$ & $2\bm{e}^{4}+2\bm{e}^6+4\bm{e}^{1}$
\\\hline
$\bm{w}_2$ & $\bm{e}^{1}+\bm{e}^5+\bm{e}^{3}$& $\bm{w}_5$ & $2\bm{e}^{5}+2\bm{e}^6+4\bm{e}^{1}$
\\\hline
&& $\bm{w}_6$ & $2\bm{e}^{4}+2\bm{e}^6+4\bm{e}^{2}$
\\\hline
&& $\bm{w}_7$ & $2\bm{e}^{5}+2\bm{e}^6+4\bm{e}^{2}$
\\\hline
\end{tabular}
\end{table}

\subsection{Settings of Simulations in Figure~\ref{fig:fig2}}
The simulations whose results are exhibited in Figure~\ref{fig:fig2} have three request types and fifteen resource pools, with capacities and cost rates per RU given in Table~\ref{table:fig2:nonopt_v1}.


\begin{table}[t]
\caption{Capacities and Cost Rates for RUs of Different Pools.}\label{table:fig2:nonopt_v1}
\renewcommand{\arraystretch}{1.2}
\centering
\begin{tabular}{|c|c|c|c|c|c|c|c|c|}
\hline
Resource Pool ($j$) & $1$ & $2$ & $3$ & $4$ & $5$ & $6$ & $7$ & $8$  \\\hline
Capacity ($C_j^0$) & $9$ & $7$ & $6$ & $7$ & $8$ & $7$ & $7$ & $6$ \\\hline
Cost Rate ($\varepsilon_j$) & $9.995$ &
$2.707$ & $2.237$ & $4.656$ & $0.624$ & $5.705$ & $0.385$ & $6.065$   \\\hline
Resource Pool ($j$) & $9$ & $10$ & $11$ & $12$ & $13$ & $14$ & $15$ & \\\hline
Capacity ($C_j^0$) & $9$& $7$& $5$& $6$ & $6$& $5$& $6$ &
\\\hline
Cost Rate ($\varepsilon_j$) & $7.492$ & $6.584$ & $1.085$ & $7.332$ & $5.862$ & $1.938$ & $8.411$ &
\\\hline
\end{tabular}
\end{table}
\begin{itemize}
\item $\mathscr{P}_1=\{1,2,\ldots,30,31\}$,  $\lambda^0_1=1.177$, $\mu_{1}=\mu_{2}=\ldots=\mu_{30}=0.235$, $R_1=3710.05$;
\item $\mathscr{P}_2=\{32,33,34,35\}$,   $\lambda^0_2=1.108$, $\mu_{32}=\mu_{33}=\mu_{34}=0.222$,  $R_2=3712.66$;
\item $\mathscr{P}_3=\{36,37,\ldots,41,42\}$,  $\lambda^0_3=1.286$, $\mu_{36}=\mu_{37}=\ldots=\mu_{41}=0.257$,  $R_3=3821.33$;
\end{itemize}
where patterns $31,35,42$ are dummy patterns for blocking requests.
The weight vectors $\bm{w}_i$ of patterns  $i\in[42]\backslash \{31,35,42\}$ are given in Table~\ref{table:fig2:nonopt_v1:weights}.
\begin{table}[t]
\caption{Weight Vectors of Patterns for Four Different RTs.}\label{table:fig2:nonopt_v1:weights}
\renewcommand{\arraystretch}{1.2}
\centering
\begin{tabular}{|c|c|c|c|c|c|c|c|}
\hline
$\bm{w}_1$ & $4\bm{e}^{3}+3\bm{e}^7+\bm{e}^{14}$ & $\bm{w}_{11}$ & $4\bm{e}^{3}+3\bm{e}^1+\bm{e}^{14}$ & $\bm{w}_{21}$ & $4\bm{e}^{2}+3\bm{e}^7+\bm{e}^{9}$ & $\bm{w}_{32}$ & $\bm{e}^{10}+3\bm{e}^5+2\bm{e}^{14}$
\\\hline
$\bm{w}_2$ & $4\bm{e}^{6}+3\bm{e}^7+\bm{e}^{14}$ & $\bm{w}_{12}$ & $4\bm{e}^{3}+3\bm{e}^1+\bm{e}^{2}$ & $\bm{w}_{22}$ & $4\bm{e}^{6}+3\bm{e}^7+\bm{e}^{9}$ & $\bm{w}_{33}$ & $\bm{e}^{10}+3\bm{e}^5+2\bm{e}^{2}$
\\\hline
$\bm{w}_3$ & $4\bm{e}^{6}+3\bm{e}^1+\bm{e}^{14}$ & $\bm{w}_{13}$ & $4\bm{e}^{12}+3\bm{e}^7+\bm{e}^{8}$ & $\bm{w}_{23}$ & $4\bm{e}^{3}+3\bm{e}^1+\bm{e}^{9}$ & $\bm{w}_{34}$ & $\bm{e}^{10}+3\bm{e}^5+2\bm{e}^{8}$
\\\hline
$\bm{w}_4$ & $4\bm{e}^{6}+3\bm{e}^1+\bm{e}^{2}$ & $\bm{w}_{14}$ & $4\bm{e}^{3}+3\bm{e}^7+\bm{e}^{8}$ & $\bm{w}_{24}$ & $4\bm{e}^{6}+3\bm{e}^1+\bm{e}^{9}$ &&
\\\hline
$\bm{w}_5$ & $4\bm{e}^{12}+3\bm{e}^1+\bm{e}^{14}$ & $\bm{w}_{15}$ & $4\bm{e}^{3}+3\bm{e}^1+\bm{e}^{8}$ & $\bm{w}_{25}$ & $4\bm{e}^{3}+3\bm{e}^7+\bm{e}^{5}$ & $\bm{w}_{36}$ & $2\bm{e}^{15}+\bm{e}^{11}+3\bm{e}^{3}$
\\\hline
$\bm{w}_6$ & $4\bm{e}^{12}+3\bm{e}^1+\bm{e}^{2}$ & $\bm{w}_{16}$ & $4\bm{e}^{6}+3\bm{e}^1+\bm{e}^{8}$ & $\bm{w}_{26}$ & $4\bm{e}^{6}+3\bm{e}^7+\bm{e}^{5}$ & $\bm{w}_{37}$ & $2\bm{e}^{15}+\bm{e}^{11}+3\bm{e}^{12}$
\\\hline
$\bm{w}_7$ & $4\bm{e}^{3}+3\bm{e}^7+\bm{e}^{2}$ & $\bm{w}_{17}$ & $4\bm{e}^{12}+3\bm{e}^1+\bm{e}^{8}$ & $\bm{w}_{27}$ & $4\bm{e}^{12}+3\bm{e}^7+\bm{e}^{5}$& $\bm{w}_{38}$ & $2\bm{e}^{15}+\bm{e}^{4}+3\bm{e}^{3}$
\\\hline
$\bm{w}_8$ & $4\bm{e}^{6}+3\bm{e}^7+\bm{e}^{2}$  & $\bm{w}_{18}$ & $4\bm{e}^{6}+3\bm{e}^7+\bm{e}^{8}$ & $\bm{w}_{28}$ & $4\bm{e}^{3}+3\bm{e}^1+\bm{e}^{5}$ & $\bm{w}_{39}$ & $2\bm{e}^{15}+\bm{e}^{4}+3\bm{e}^{12}$
\\\hline
$\bm{w}_9$ & $4\bm{e}^{12}+3\bm{e}^7+\bm{e}^{14}$ & $\bm{w}_{19}$ & $4\bm{e}^{12}+3\bm{e}^7+\bm{e}^{9}$ & $\bm{w}_{29}$ & $4\bm{e}^{6}+3\bm{e}^1+\bm{e}^{5}$ & $\bm{w}_{40}$ & $2\bm{e}^{15}+\bm{e}^{13}+3\bm{e}^{3}$
\\\hline
$\bm{w}_{10}$ & $4\bm{e}^{12}+3\bm{e}^7+\bm{e}^{2}$ & $\bm{w}_{20}$ & $4\bm{e}^{12}+3\bm{e}^1+\bm{e}^{9}$ & $\bm{w}_{30}$ & $4\bm{e}^{12}+3\bm{e}^1+\bm{e}^{5}$ & $\bm{w}_{41}$ & $2\bm{e}^{15}+\bm{e}^{13}+3\bm{e}^{12}$
\\\hline
\end{tabular}
\end{table}





\end{document}